\documentclass[a4paper,reqno]{amsart}

\usepackage{hyperref}
\usepackage{amssymb}
\usepackage{amsmath}
\usepackage{amsfonts}
\usepackage{mathtools}
\usepackage{mathrsfs}
\usepackage{enumitem}
\usepackage{amsthm}
\usepackage{mathrsfs}
\usepackage{color}
\usepackage{dsfont}
\usepackage{a4wide}
\usepackage{tikz-cd}

\usepackage{euscript} 

\numberwithin{equation}{section}
\usepackage{parskip}

\makeatletter
\def\thm@space@setup{%
  \thm@preskip=\parskip \thm@postskip=0pt
}
\makeatother

\allowdisplaybreaks

\AddToHook{env/abstract/begin}{%
  \edef\restoreabstractparagraphs{%
    \parindent=\the\parindent\relax
    \parskip=\the\parskip\relax
  }%
}
\AddToHook{cmd/abstract/after}{\restoreabstractparagraphs}

\usepackage{etoolbox}

\makeatletter

\patchcmd{\@setauthors}
  {\MakeUppercase{\authors}}
  {\large\authors}
  {}
  {\PackageWarning{mystyle}{Could not change author capitalization}}

\patchcmd{\ps@firstpage}
  {\thepage}
  {}
  {}
  {\PackageWarning{mystyle}{Could not remove first-page number}}

\makeatother

\makeatletter

\patchcmd{\@maketitle}
  {\global\topskip42\p@\relax}
  {\global\topskip28\p@\relax}
  {}
  {\PackageWarning{mystyle}{Could not change title spacing}}

\makeatother

\usepackage{hyperref}
\hypersetup{
    colorlinks = true,
    citecolor  = blue!75!black,  
    urlcolor   = blue!75!black,   
    linkcolor  = red!70!black 
}

\newtheorem{theorem}{Theorem}[section]
\newtheorem{corollary}[theorem]{Corollary}

\newtheorem{lemma}[theorem]{Lemma}
\newtheorem{proposition}[theorem]{Proposition}
\newtheorem{remark}[theorem]{Remark}

\newtheorem{ltheorem}{Theorem}

\newtheorem{lproposition}[ltheorem]{Proposition}

\theoremstyle{definition}
\newtheorem*{nnoproblem}{Problem}

\newcommand{\cA}{\mathcal{A}}
\newcommand{\cB}{\mathcal{B}}

\newcommand{\cD}{\mathcal{D}}
\newcommand{\cE}{\mathcal{E}}
\newcommand{\cF}{\mathcal{F}}

\newcommand{\cH}{\mathcal{H}}

\newcommand{\cL}{\mathcal{L}}
\newcommand{\cM}{\mathcal{M}}

\newcommand{\cP}{\mathcal{P}}

\newcommand{\cR}{\mathcal{R}}

\newcommand{\cW}{\mathcal{W}}

\newcommand{\1}{\mathds{1}}
\newcommand\imagi{\mathrm{i}}
\newcommand{\id}{\mathrm{id}}
\newcommand{\cb}{\mathrm{cb}}

\newcommand{\into}{\hookrightarrow}
\newcommand{\acts}{\curvearrowright}

\newcommand{\CC}{\mathbb{C}}

\newcommand{\EE}{\mathbb{E}}
\newcommand{\FF}{\mathbb{F}}

\newcommand{\NN}{\mathbb{N}}

\newcommand{\ZZ}{\mathbb{Z}}

\newcommand{\rad}{\mathrm{rad}}
\newcommand{\ang}{\mathrm{ang}}

\newcommand{\RR}{\mathbb{R}}

\newcommand{\Fr}{\mathbb{F}_r}
\newcommand{\Erad}{\mathbb{E}_{\mathrm{rad}}}

\newcommand{\cLFr}{\mathcal{L} \mathbb{F}_r }
\newcommand{\LFr}{\mathcal{L}(\mathbb{F}_r)}
\newcommand{\CFr}{\mathbb{C}[\mathbb{F}_r]}
\newcommand{\LFrrad}{\cLFr^\rad}

\newcommand{\CFrrad}{\mathbb{C}[\mathbb{F}_r]^{\rad}}

\makeatletter
\newcommand{\avgint@overlay}[3]{%
  {%
    \setbox0=\hbox{$#1{#2#3}{\int}$}%
    \vcenter{\hbox{$#2#3$}}%
    \kern-.5\wd0
  }%
}

\newcommand{\avgint@build}[1]{%
  \mathchoice
    {\avgint@overlay\displaystyle\textstyle{#1}}%
    {\avgint@overlay\textstyle\scriptstyle{#1}}%
    {\avgint@overlay\scriptstyle\scriptscriptstyle{#1}}%
    {\avgint@overlay\scriptscriptstyle\scriptscriptstyle{#1}}%
  \!\int
}

\newcommand{\avgint}{\avgint@build{-}}
\makeatother

\begin{document}

\title[\emph{Mixed-norm spaces over free groups}]{Radial multipliers on mixed-norm spaces over free groups}

\author[]{Adrián González-Pérez}
\author[]{Binbing Lyu}
\author[]{Javier Parcet}
\author[]{Jorge Pérez-García}
\author[]{\\Simeng Wang}
\author[]{Tiantian Zhao}

\begin{abstract}
    The problem of Fourier summation on free groups remains open and challenging. More precisely, if $B_N(e)$ are the balls of radius $N$ on $\mathbb{F}_r$, it is still open whether they induce a uniformly bounded family of Fourier multipliers on the non-commutative $L_p$-spaces associated with their regular von Neumann algebras. Nevertheless, progress  has been made. In particular, Kuhn and later Bozejko/Fendler showed that the ball multipliers are non-uniformly bounded when restricted to the radial subalgebra of the free group if $p \not\in (\tfrac{3}{2}, 3)$. 

We will produce a converge to their results by introducing mixed $L_p(L_2)$-norms given by 
\[
    \| x \|_{p, \mathrm{mix}} := \big\| \mathbb{E}_{\mathrm{rad}} \big[ x^\ast \, x \big]^\frac12 \big\|_p,
\]
where $\mathbb{E}_\mathrm{rad}$ is the expectation onto the radial subalgebra. These spaces can be seen as free analogues of the mixed $L_p$-radial/$L_2$-angular norms studied in the Euclidean case by Mockenhaupt, Rubio de Francia and A. Cordoba. We show that the ball multipliers are uniformly bounded on these mixed-norm spaces exactly when $p \in (\tfrac{3}{2}, 3)$, thus confirming intuitions going back to Biane about the dual of the free group behaving as a $3$ dimensional object.

Lastly, using Fourier-Schur comparison, we show that the balls associated with the $q$-convex deformations of the word length function of $\mathbb{F}_r$ have divergent completely bounded norms as multipliers on $L_p$ spaces if $q > 1$. When $q = 1$, recent results on radial Schur multipliers show that the strategy cannot possibly work.
\end{abstract}

\maketitle

\section{Introduction}

The study of the convergence of Fourier series and integrals is a natural problem in harmonic analysis that naturally leads to the study of the boundedness of Fourier multipliers. 
In Euclidean spaces, the ball multiplier is defined by 
\[
    \widehat{T(f)}(\xi)=\chi_{B}(\xi)\widehat{f}(\xi),\quad f\in L_p(\mathbb{R}^n),
\] 
where $B$ is the unit ball $B(0,1) \subset \mathbb R^n$ and $\widehat f$ denotes the Fourier transform. The boundedness of $T$ on $L_p(\mathbb R^n)$ is equivalent to the convergence of Fourier integrals over Euclidean balls and thus it has received considerable attention in the past. 
In dimension \(n=1\), the boundedness of the ball multiplier reduces to that of the Hilbert transform. As such, its boundedness in $L_p$ has been well understood for close to a century. 
In sharp contrast, in dimensions \(n>1\), Herz \cite{He1954} proved that \(T\) is a bounded operator on the subspace of radial functions in \(L_p(\mathbb{R}^n)\) if and only if \(2n/(n+1)<p<2n/(n-1)\). Following the work of Herz, it was conjectured that \(T\) would be bounded on \(L_p(\mathbb{R}^n)\) throughout the same range. 
In the early 1970s, C. Fefferman disproved this conjecture in his famous work \cite{Fe1971} showing that $T$ is unbounded in $L_p(\mathbb R^n)$ for every $p\neq2$.
%
%
These early results and techniques would grow into a small subfield of harmonic analysis exhibiting deep connections with topics such as the Bochner--Riesz conjecture  \cite{Fe1972,CS1972,Co1979,Zy2002}, Kakeya maximal functions \cite{Fe1971,Co1974,Co1979,W1995}, and Fourier restriction theory \cite{T1975,Ve1988,S1993}, that has grown in parallel to techniques such as weighted and vector-valued inequalities.

The study of the ball multiplier has progressed in parallel with techniques like geometric maximal functions, weight theory and vector-valued inequalities.
Indeed, Fefferman's proof of the unboundedness of the ball multiplier on  \(L_p\) is based on Meyer's observation that the \(L_p\)-boundedness of the ball multiplier would imply a square-function estimate for directional Hilbert transforms. 
Fefferman then built his counterexample by introducing a Kakeya configuration contradicting Meyer's square-function inequality after the action of the Hilbert transforms.
This eventually led to the guiding question of finding geometrically defined maximal operators $\cM$ controlling the behavior of the ball multiplier
\begin{equation}
    \label{T2wM}
	\int|T f|^2g \leq C\, \int|f|^2\cM(g).
\end{equation}
In \cite{St1979}, see also \cite{Co1977}, it was conjectured that $\cM$ could be taken as a $L_q$-convexified universal maximal over rectangles. While the conjecture is still open, much more is currently known, see \cite{Co1984, CRS1992} as well as \cite{RdF1982Weighted, Car1985} for the analogous inequality for Bochner-Riesz multipliers. In particular, by \cite{CHS1991}, Stein's conjectured bound holds when $g$ is a radial function. Through duality, these weighted \(L^2\) estimates allow multiplier bounds to be reduced to maximal function estimates.

\subsection*{Mixed radial-angular spaces}
For our own purposes, results about mixed-norm inequalities for the ball multiplier and for Bochner-Riesz multipliers are of utmost interest. Indeed, passing to polar coordinates, we can decompose $\RR^n$ as $\RR_+ \times S^{n-1}$. The mixed radial-angular norms are then defined by taking the $L_p$-norm on the radial component with values in $L_2(S^{n-1})$ for the \emph{spherical} or \textit{angular} component. This gives the mixed-norm
\begin{align*}
    \| f \|_{L_p^\rad(L_2^\ang)}
    :=
    \bigg(
    \int_{0}^\infty
        \Big(\int_{S^{n-1}}|f(r \, u)|^2 \, d\sigma(u) \Big)^{\frac{p}{2}}r^{n-1}dr
    \bigg)^{\frac{1}{p}},
\end{align*}
where $d \sigma$ is the invariant probability measure on $S^{n-1}$. The literature on Fourier multiplier bounds on mixed radial-angular spaces starts with the pioneering work of Rubio de Francia \cite{Ru1989}, who showed that over $L_p^\rad(L_2^\ang)$ on even-dimensional spaces, the Bochner-Riesz multipliers are bounded in the range predicted by the Bochner-Riesz conjecture. Independently, Córdoba \cite{Co1989} and Mockenhaupt \cite{Mo1990} extended the approach to the limit case of the ball multiplier. Their results state that over $\RR^n$
\begin{equation}
    \label{eq:MNtheoren}
    \tag{MN}
    T :L_p^\rad(L_2^\ang) \to L_p^\rad(L_2^\ang) \; \text{ boundedly }
    \;\;
    \Longleftrightarrow
    \;\;
    \Big| \frac1{2} - \frac1{p} \Big| < \frac{1}{2n},
\end{equation}
see \cite{Co1989, CHS1991, CRS1992, DMOS2008} as well. At this point, it is worth expanding on the techniques used to prove the above theorem. First, when we restrict $T$ to radial functions, we obtain, through the use of Bessel harmonics, a singular integral operator acting on the positive half-line $\RR_+$ with the weight given by $w(r) dr = r^{{n-1}-p\frac{n-1}{2}} d r$. This weight sits in the Muckenhoupt's $A_p$ class with $p$ in the range described in \eqref{eq:MNtheoren}. This is enough to obtain the radial result of Herz. In order to pass from the radial subspace to mixed-norm spaces, the key is decomposing the angular part as
\[
    L_2(S^{n-1}) = \bigoplus_{\ell \geq 0} \cH_\ell \; \implies \; L_p^\rad(L_2^\ang) \cong L_p\Big((0,\infty), r^{n-1}dr; \, \bigoplus_{\ell \geq 0}  \cH_\ell\Big),
\]
where each of the summands $\cH_\ell$ denotes the space of spherical harmonics of degree $\ell$. The rotational symmetry of the ball multiplier then implies that $T$ splits as a direct sum of singular integral operators $T_\ell: L_p(\RR_+,r^{n-1} dr; \cH_\ell) \to L_p(\RR_+,r^{n-1} dr; \cH_\ell)$ with scalar kernels and the integral kernels of each of the operators $T_\ell$ are explicitly computable through special functions. This reduces the bound in \eqref{eq:MNtheoren} to a weighted square function estimate. The quantitative part of the theorem is to obtain that square function estimate via duality arguments, weighted one dimensional singular integrals, and control by maximal operators; see \cite{CRS1992} for the full details. Lastly, a common intuition regarding mixed radial-angular spaces is that $L_p^\rad(L_2^\ang)$ is blind to the angular components and this nullifies the potential effects of Kakeya sets that appear in dimensions $n \geq 2$.

\medskip
\subsection*{Analysis on free group algebras}
The topic of this article is the extension of the mixed radial-angular results above, from the context of Euclidean harmonic analysis to free group algebras. The study of noncommutative harmonic analysis over free groups goes back to the pioneering work of Haagerup \cite{Haa1978} who proved that the reduced $C^\ast$-algebra of the free group $\Fr$ has Grothendieck's metric approximation property. His result used a clever combination of the conditional negativity of the word length, see \cite{CheCowJoJulgVal2001}, and what later became known as Haagerup's inequality. Later on, de Canni\`ere, Cowling, and Haagerup \cite{DH1985,CH1989} extensively studied completely bounded Fourier multipliers on Lie group von Neumann algebras and the corresponding completely bounded approximation property, providing one of the first examples of non-isomorphism between nonamenable group von Neumann algebras. Nowadays, the theory of noncommutative Fourier multipliers has become a useful tool in deformation/rigidity theory, strong solidity, and absorption properties \cite{OP2010,PV2014,Po2007,HI2024}. Moreover, in recent decades, the study of the ${L_p}$-boundedness of Fourier multipliers on group von Neumann algebras has been an active and challenging topic. 
One of the first breakthroughs was obtained by Lafforgue and de la Salle \cite{LD2011}, who proved that the group von Neumann algebra of higher rank lattices and their associated noncommutative $L_p$-spaces fail to have the completely bounded approximation property when $p$ is large, using Fourier and Schur multiplier techniques. Further advances on Riesz transforms, Hilbert transforms, and H\"ormander--Mikhlin theory have recently been obtained; see \cite{JMP2014,GJP2017,MR2017,JMP2018,PRD2022}.


Despite the centrality of free groups in noncommutative harmonic analysis, establishing $L_p$-boundedness for radial Fourier multipliers remains challenging. For instance, a spectral Hörmander-Mikhlin theorem for the word length of the form $\gamma \mapsto m(|\gamma|)$ with finite smoothness order is currently unavailable. Similarly, the range of exponents $p$ for which the ball multipliers are uniformly bounded is also unknown. Let $\Fr$ be the free group on $r$ generators, equipped with the natural word-length function $g \mapsto |g|$, and denote by $\cL \Fr$ its von Neumann algebra and by $L_p(\cL \Fr)$ the noncommutative $L_p$-space associated to $\cL \Fr$ with the Plancherel trace. The question of the $L^p$-boundedness of the ball multipliers goes back to the celebrated work of Junge, Nielsen, Ruan, and Xu \cite{JNRX2004}, who proved that ${L_p}(\cL \Fr)$ admits a (completely bounded) Schauder basis for $1<p<\infty$. This is a much stronger property than merely having the metric approximation property or the completely bounded approximation property. However, explicit constructions of such bases have remained elusive. In fact, in \cite{JNRX2004}, the authors asked whether the group elements, arranged according to the word length on $\Fr$, could constitute a Schauder basis of $L_p$ for $1<p<\infty$. This question amounts precisely to asking for the uniform boundedness of the ball multipliers on the free group. 
\begin{nnoproblem}[{The free ball multiplier problem}]
    \label{prb:BallProblem}
    Let $S_N(\lambda_\gamma) = \chi_{\{|\gamma|\leq N\}} \, \lambda_\gamma$ be the ball multiplier on the free group. For which $1 < p < \infty$ does one have that
    \[
        \sup_{N} \big\| S_N: L_p(\cL \Fr) \to L_p(\cL \Fr) \big\| < \infty \, \text{?}
    \]
\end{nnoproblem}
This motivated the work of Bo\.zejko and Fendler \cite{BF2006}, who studied the $L_p$-boundedness of $S_N$ when restricted to the commutative \emph{radial} subalgebra. Indeed, we will say that a function $\varphi: \Fr \to \CC$ is \emph{radial} if its values $\varphi(g)$ depend only on the word length of $g$. It is a well-known fact that the radial functions form an Abelian $\ast$-algebra under convolution, see \cite{Pytlik1981, Cohen1982Operator, FP1983}. We will refer to its weak-$\ast$ closure as the radial subalgebra $\cL \Fr^{ \rad} \subset \cL \Fr$. This radial subalgebra can be conceptualized as a frequency analogue of the radial functions on the Euclidean space. Their properties,  however, differ substantially from the Euclidean situation. Most crucially $\cL \Fr^\rad$ is not a priori the fixed point algebra of a group action, nor a tensor factor. In fact, a common intuition is that $\cL \Fr^\rad$ behaves analogously to the algebra of $K$-bi-invariant functions of $SL_2(\RR)$. Using the spherical transform of the radial subalgebra, Bo\.zejko and Fendler proved that $S_N$ is not uniformly bounded on $L_p(\cL \Fr^\rad)$ for $p<\tfrac{3}{2}$ and $p>3$ (\'Eric Ricard communicated to us that it is also not uniformly bounded for $p=\tfrac{3}{2}$ and $p=3$, by refining the same arguments), but is uniformly bounded for $\tfrac{3}{2}<p<3$, thereby answering the aforementioned free ball problem in the negative for $p\leq \tfrac{3}{2}$ and for $p\geq 3$. However, since the radial part of $L_p(\cL \Fr)$ is a proper subspace of $L_p(\cL \Fr)$, the uniform boundedness over the range $p \in (\tfrac{3}{2}, 3)$ on the radial part does not provide an answer in this smaller range. Similar results were also obtained by Kuhn in an earlier work \cite{K1984}. These results should be regarded as a noncommutative analogue of the work of Herz in the Euclidean setting. As in the Euclidean setting, the study of $L_p$-boundedness on the non-radial part of ${L_p}(\cL \Fr)$ appears to be much more delicate, and remain underexplored.

\subsection*{Our results}
This article takes a first step towards understanding mixed-norm estimates of ball multipliers on free groups, especially on the non-radial part, in close analogy with the classical mixed-norm results of C\'ordoba \cite{Co1989}, Mockenhaupt \cite{Mo1990}, and Carbery et al.~\cite{CRS1992} in the Euclidean case. These mixed-norm estimates in the noncommutative setting are not only motivated by the analogous classical theory, but also arise naturally from noncommutative analysis and operator algebras. Indeed, Parcet \cite{P2025} recently noted that similar mixed radial-angular Schatten-$p$-norm estimates for Fourier--Schur multipliers on a simple Lie group $G$ would imply completely bounded approximation by Schur multipliers on \(S_p(G)\); the latter is currently well known to be a core open problem in noncommutative analysis, deeply connected to Connes' W*-rigidity conjecture in von Neumann algebra theory and the Kakeya conjecture in modern harmonic analysis \cite{PRD2022,De2023,P2025,De2026}.

More precisely, let $\Erad\colon \cL \Fr \to \cL \Fr^{\rad}$ be the unique trace-preserving conditional expectation between these two von Neumann algebras. Then we can define the free version of the classical radial-angular $L_p$-norm as follows:
\begin{equation}
    \label{mixed-norms}
		\left\|x\right\|_{{L_p}(\cLFr; \Erad^c)}
		\coloneqq
		\left\|
		\Erad
		\big[x^*x\big]^{\frac12}
		\right\|_{p},\quad 
		\left\|x\right\|_{{L_p}(\mathcal L(\mathbb F_r);
			\mathbb E_{\mathrm{rad}}^r)}
		\coloneqq
		\left\|
		\Erad
		\big[xx^*\big]^{\frac12}
		\right\|_{p}.
\end{equation}
These norms give rise to two non-equivalent spaces called \textit{the column and the row mixed-norm spaces}, respectively, and were systematically studied by Junge and Parcet \cite{JP2010} in a more general setting. With this notation, our main result is as follows.

\begin{ltheorem}
    \label{Tm-pp}  
    Let \(r\geq 2\)  and \(\frac{3}{2}< p<3\). Then there exists a constant $C_{p}>0$ depending only on $p$ such that
	\[ 
        \sup_{N\ge0 }\left \| S_N(x) \right \|_{{L_p}(\cLFr;\Erad^c)}\leq C_{p} \left \| x \right \|_{{L_p}(\cLFr;\Erad^c)},
    \]
	for all \(x\in L_p(\cL \Fr;\mathbb{E}_{\mathrm{rad}}^c)\). The analogous assertion also holds for
	\(L_p(\mathcal{L}(\mathbb{F}_r);\mathbb{E}_{\mathrm{rad}}^r)\).
\end{ltheorem}
Together with the aforementioned negative result by Kuhn and Bo\.zejko--Fendler for $p\leq \tfrac{3}{2}$ and $p\geq 3$, this theorem yields a sharp range of $p$ for the boundedness of ball multipliers in the mixed norms spaces of the free group. This result establishes the analogue of the aforementioned theory of Rubio de Francia, C\'ordoba, Mockenhaupt, and Carbery--Romera--Soria \cite{CRS1992} on free groups. To the best of our knowledge, this is the first sharp mixed-norm theorem of this type for a genuinely noncommutative group von Neumann algebra.

As in the Euclidean setting, the above result follows from the noncommutative weighted $L_2$ estimates below. Identifying $\cL \Fr^\rad$ with the commutative von Neumann algebra of an interval $I$ depending on $r$, $L^\infty(I)$ in a standard way, see Section \ref{subsec: radial subalgebra definition}, we consider the Hardy--Littlewood maximal operator $\cM$ on $\cLFr^\rad$ and write $\mathcal M_\alpha(f) = \mathcal M (f^\alpha)^{\frac1{\alpha}}$. The following establishes a noncommutative analogue of inequality \eqref{T2wM}, which was extensively studied in  \cite{CRS1992,Be2014,Co2016} in the Euclidean setting.
\begin{ltheorem}
    \label{mainthm}
	Let \(r\ge 2\) and \(\alpha>2\).  There exists a constant  \(C_\alpha>0\) depending only on \(\alpha\) such that for every $N\geq 0$, \(x\in \cLFr\) and  \(w\in (\cLFr^{\rad})_{+}\),  we have
	\begin{equation*} 
		\tau\Big(|S_N(x)|^2\cdot w\Big)\le C_{\alpha}\tau\Big(|x|^2\cdot \mathcal M_\alpha(w)\Big).
	\end{equation*}
\end{ltheorem}

Note that, in the Euclidean setting, the restriction of the Hardy--Littlewood maximal operator to radial functions is controlled by the Kakeya maximal operator---this is precisely why Kakeya maximal functions appear in C\'ordoba's version of inequality \eqref{T2wM}. Given the analogous geometry of Cayley graphs, our result involving a Hardy--Littlewood maximal function on the free radial part naturally suggests the intuition of a ``noncommutative Kakeya maximal function'' on free groups, which is of independent interest and will be helpful for further developments in harmonic analysis on free groups. As it is not required for the proof of our main results, we will pursue this direction elsewhere.

Lastly, we deal with an analogue of the free ball multiplier problem described above. Recall that if \(1 \leq q < \infty\) is an exponent, the following expression
\begin{equation}
    \label{def:qlength}
    \big| s_1^{n_1} \, s_2^{n_2} \cdots s_m^{n_m} \big|_q
    \; = \;
    \bigg( \sum_{j = 1}^m |n_j|^q \bigg)^\frac1q
\end{equation}
gives a length function on the free group $| \cdot |_q: \Fr \to \RR_+$. Furthermore, it is conditionally negative precisely when $1 \leq q \leq 2$, see \cite{Bozejko1986, ChuahLiuMei2023}. We will show that the balls associated to the $q$-length above are non-uniformly completely bounded in $L_p(\mathcal{L} \Fr)$. In the statement below, the $\cb$ represents the completely bounded norm of an operator.

\begin{ltheorem}
    \label{thm:qBall}
    Let $B_N^q(e) = \big\{ \gamma \in \Fr : |\gamma|_q \leq N \big\}$ and $S_N^q(\lambda_\gamma) = \chi_{B_N^q(e)}(\gamma) \lambda_\gamma$ be the Fourier multipliers associated to the characteristic function of $B_N^q(e)$, then, 
    \[
        \sup_{N \geq 0} \big\| S_N^q: L_p(\mathcal{L} \Fr) \to L_p(\mathcal{L} \Fr) \big\|_\cb = \infty,
        \qquad 
        \text{ for every } 
        \; q > 1 \; \text{ and } p \in [1,\infty] \setminus \{2\}.
    \]
\end{ltheorem}
The proof of the above theorem uses that the completely bounded operator norm on $L_p$ of a Fourier multiplier is bounded below by that of its associated Herz-Schur multiplier acting on its $p$-Schatten class $S_p$, see \cite{NeuRic2011, CasSall2015}. Then, restricting to a two-parameter subspace of indices and using an approximation argument gives the boundedness of an explicit idempotent Schur multiplier on $\RR^2$, which is easily seen to be unbounded due to the results in \cite{Parcet2024Idempotentschur}. It is also interesting to note that, while these $q$-lengths approximate the usual one as $q \to 1^+$, the proof schema fails for $q=1$ since, by \cite{GPPW2026}, the usual free balls $B_N(e)$ induce uniformly (completely) bounded Herz-Schur multipliers in $S_p$.

\subsection*{The techniques}
Let us briefly outline the new ideas and ingredients in the proof of the main results since they carry intrinsic value beyond their use for the proofs of Theorems \ref{Tm-pp} and \ref{mainthm}. We will start by reviewing some well-known facts on spherical analysis on the free groups. Let 
\begin{equation}
    \label{eq:muNdefs}
    \mu_n = \frac1{|\mathbf{S}_n|} \sum_{|\gamma| = n} \delta_\gamma,
\end{equation}

be the $\ell_1$-normalized sum over the coronas of radius $n$, denoted by $\mathbf{S}_n = \{|\gamma| =n \}$. These elements lie naturally within the group algebra $\CC[\Fr]$ and span the linear subspace of radial functions. An elementary computation yields that
\begin{equation}
    \label{mun}
    \mu_1 \ast \mu_n 
    =
    \frac{1}{2r} \mu_{n-1} + \frac{2r-1}{2r} \mu_{n+1}, 
    \quad \text{ for every } \, n \ge 1.
\end{equation}
Therefore, using an elementary recurrence argument, the radial subalgebra of $\CC[\Fr]$ coincides with the unital $\ast$-algebra generated by $\mu_1$ and, as such, is Abelian. As usual, we will denote by $\lambda \colon \Fr\to B(\ell_2\Fr)$ the left regular representation of $\Fr$ on the Hilbert space $\ell_2(\Fr)$. Let $\cL \Fr^\rad \subset B(\ell_2 \Fr)$ be the weak-$\ast$ closure of the radial algebra. Since it is generated by $\mu_1$, its Gel'fand spectrum coincides with the spectrum of $\mu_1$ and the restriction of the usual trace $\tau$ to $\cLFr^\rad$ gives a measure over it. Explicitly, we have 
\begin{equation}
    \label{eq:Gelfand}
    \big(\cL \Fr^\rad,\tau{|}_{\cL \Fr^\rad} \big) \xrightarrow{\quad \Gamma \quad} L_\infty(I,m_r) 
\end{equation}
where $I := \sigma(\lambda(\mu_1))$ is a symmetric interval depending on the rank $r$ and $dm_r$ is an absolutely continuous measure over $I$, see Section \ref{Preliminaries} for the details. Both the formula for the spectrum of $\mu_1$ and the spherical Plancherel measure $dm_r$ go back to the pioneering work of Kesten \cite{Kesten1959}, see as well \cite{Cartier1973Geometrie}, while their connection to the reduced free algebras was studied in \cite{Cohen1982Operator, CohenDeMichele1982, Pytlik1981}. The Gel'fand transform in \eqref{eq:Gelfand} has a particularly nice expression given by sending its generator $\mu_1$ to the identity function $s \mapsto s$ on $I$. The recurrence law \eqref{mun} then implies that that each $\mu_n$ is mapped by $\Gamma$ to an orthogonal polynomial $P_n \colon I \to \RR$ satisfying the recurrence relations in \eqref{Pn}. 

The unique trace-preserving conditional expectation $\EE_\rad \colon \cL \Fr \to \cL \Fr^\rad \subset \cL \Fr$ is given by averaging the group algebra elements over coronas and, since those coronas are finite, $\EE_\rad$ maps the reduced $C^\ast$ of the free group onto its radial part. Composing this expectation with evaluation on $s \in I$ gives states on the reduced $C^\ast$-algebra, ie:
\begin{equation}
    \label{def: state phi_s}
    \begin{tikzcd}[column sep=large]
        C_{\lambda}^{*} \Fr
        \arrow[r, "\EE_{\rad}"]
        \arrow[rrr, bend right=22, "\phi_s"]
        &
        C_{\lambda}^{*} \Fr^{\rad}
        \arrow[r, "\Gamma"]
        &
        C(I)
        \arrow[r, "\text{ev}_s"]
        &
        \CC
    \end{tikzcd}
\end{equation}
The state $\phi_s \colon C_\lambda^\ast \Fr \to \CC$ induces a radial positive definite function $\varphi_s \colon \Fr \to \CC$ on the group given by $\varphi_s(\gamma) = \phi_s(\lambda_\gamma)$. In the above case, we have $\phi_s(\lambda_\gamma) = P_{|\gamma|}(s)$ and that those are all the extremal point of the convex set of radial normalized positive definite functions that also factor through $\lambda$, see \cite{Arnaud1980, CohenDeMichele1982, Pytlik1981} or \cite{HaagerupKnudby2015} for a recent account. 
%
%
Now, taking an element $x \in \CC[\Fr]$ and expanding the $L_p$-norm on the definition of the mixed norm of $L_p(\cL \Fr; \EE_\rad)$ in \eqref{mixed-norms} gives an integral against $dm_r$. Thus for every $1 \leq p < \infty$, we have that $\Gamma$ can be extended to an isomorphism 
\begin{equation}
    \label{eq: mixed space as a field of hilberts}
    L_p(\cL \Fr; \EE_\rad^c) \; \cong \; L_p \big( I, \, dm_r \, ; \big(\EuScript{H}_s\big)_{s \in I} \big), 
    \qquad \text{ where } 
    \EuScript{H}_s = L_2(C_\lambda^* \Fr; \phi_s).
\end{equation}
In the expression above $L_2(C_\lambda^\ast \Fr; \phi_s)$ is the Gel'fand-Neumark-Segal construction associated to the state $\phi_s$ and the right hand side of the isomorphism is a mixed $L_p(L_2)$ associated with the measurable field of Hilbert spaces $s \mapsto \EuScript{H}_s$, see \cite[part II.1]{Dix1981} for the specific definitions on measurable fields.
%
As a key ingredient, we employ the following \emph{Plancherel theory} of free groups \cite{FP1983} in order to transform all the GNS representations $\EuScript{H}_s$ into boundary representations of $\Fr$. Indeed, let $\partial \Fr$ be the hyperbolic boundary of the group, given by all possible infinite geodesic paths starting at $e$. The boundary has a natural Borel structure generated by the sets of infinite words that start with any given $m$-word $\gamma$, ie:
\begin{equation*}
    \label{eq: Neighboorhoodgamma}
    E(\gamma) := \big\{ \omega \in \partial \Fr : \omega \; \text{ starts with } \gamma \big\}
\end{equation*}
and a natural probability measure $\nu$ given by setting $\nu(E(\gamma))$ to be the fraction of $\mathbf{S}_m$ covered by $\gamma$, see Section \ref{Preliminaries}. The natural action $\FF_r \curvearrowright \partial \Fr$ is clearly not measure-preserving. Nevertheless, it preserves the measure class and we have a family of unitary representations on $L_2(\partial \Fr)$ parametrized by $s \in I$ and given by
\begin{equation}
    \label{def:pseries_rep}
    \pi_s(\gamma) \xi (\omega) = \Big( \frac{d \nu_\gamma}{d \nu}(x)\Big)^{\frac12 + \imagi t } \xi(\gamma^{-1} \omega)
    \qquad \text{ where } \; s = \frac1{2 r}\big[(2r - 1)^{\frac12 + \imagi t} + (2r - 1)^{\frac12 - \imagi t}\big].
\end{equation}
In the above expression $\nu_\gamma(E) = \nu(\gamma^{-1}E)$ and $t$ is set to be in the brach of the inverse given by $[0, {\pi}/{\log(2r-1)}]$. It is also well known that $\varphi_s(\gamma)$ is actually a matrix coefficient of $\pi_s$, ie: $\varphi_s(\gamma) = \langle \pi_s(\gamma) \1, \1 \rangle$ and, by the uniqueness of the GNS representation, each $\mathcal{H}_s$ is isomorphic to $L_2(\partial \Fr)$. This is already useful since it allows us to define an explicit isomorphism
\[
    L_p(\cL \Fr; \EE_\rad^c) \xrightarrow{\quad \Phi \quad} L_p\big(I, dm_r; L_2(\partial \Fr)\big),
    \quad \text{ given by } \; \Phi(\lambda(f)) = \big( \pi_s(f) \1 \big)_{s \in I}.
\]
So far, the facts described are well known. Our contribution is the decomposition of $L_2(\partial \Fr)$ using a natural martingale sequence in a way that simultaneously decomposes $L_p(\cL \Fr; \EE_\rad^c)$ into a direct sum of right $\cL \Fr^\rad$-modular factors that are preserved by every radial Fourier multiplier. First, notice that the sets $E(\gamma)$ for $|\gamma| = m$ generate an atomic (finite) sub-$\sigma$-algebra $\cF_m$ of the boundary with the sequence $(\cF_m)_{m \geq 0}$ being an increasing filtration approximating the Borel $\sigma$-algebra of $\partial \Fr$. This filtration has a natural counterpart on the free group that has already been explored in \cite{GJP2017}. Namely, $\cA_m \subseteq \cP(\mathbf{S}_n)$ is the collection of equivalence classes under the following equivalence relation
\[
    \gamma \, \sim_{m} \, \eta
    \quad \text{ if and only if they shared starting words of length } m,
\]
where, as usual, if either word has length at most \(m\), we require the words to be equal. The finite filtration $(\cA_m)_{0 \leq m \leq n}$ intuitively approximates $(\cF_m)_{m \geq 0}$ as $n \to \infty$. This finite filtration induces an orthogonal decomposition into martingale difference spaces as follows
\[
    \ell_2(\mathbf{S}_n) \cong \bigoplus_{0 \leq m \leq n} V_{m,n},
    \quad \text{ where } \; 
    V_{m, n} = \ell_2(\mathbf{S}_n; \cA_m) \ominus \ell_2(\mathbf{S}_n; \cA_{m-1}),
\]
with the usual convention that $\ell_2(\cA_{-1}) = \{0\}$. Similarly $L_2(\partial \Fr)$ can be decomposed as a direct sum of finite dimensional martingale difference spaces $\cH_m = L_2(\partial \Fr; \cF_{m}) \ominus L_2(\partial \Fr; \cF_{m-1})$. Thus, we have a decomposition
\begin{equation}
    \label{eq:directDecomposition}    
    L_p(\cL \Fr; \EE_\rad^c) \xrightarrow{\; \Phi \; } L_p(I, dm_r; L_2(\partial \Fr)) \xrightarrow{ \; \cong \; } 
    L_p\bigg(I, dm_r; \bigoplus_{m \geq 0} \cH_m\bigg),
\end{equation}
and, crucially, the map $\Phi$, after decomposing it as a direct sum, has an explicit formula. Indeed, the preimage of $L_2(I; \cH_m)$ would be the direct sum over $n \geq m$ of all the finite dimensional subspaces $V_{n,m}$. Furthermore, there are explicit sequences of orthogonal functions on $(I,dm_r)$, that we will denote by $(Q_{m,n}(x))_{n \geq m}$ such that
\[
    \Phi(\lambda(f)) = \frac{Q_{m,n}(x)}{\| Q_{m,n}\|_2} \, \otimes \, \cW_n(f) \in L_2(I, dm_r) \otimes_2 \cH_m, 
    \qquad 
    \text{ for every } \; f \in V_{m,n},
\]
where $\cW_n$ is the natural isometry $\cW_n\colon \ell_2(\mathbf{S}_n) \to L_2(\partial \Fr; \cF_n)$ that sends $\delta_\gamma$, where $|\gamma| = n$, to $\sqrt{|\mathbf{S}_n|} \, \chi_{E(\gamma)}$. The recurrence formula for $Q_{m,n}$, as well as the details of the proof of the identities above, can be read in Section \ref{sec: kernels}. Notice that, if we fix $m=0$, then $V_{0,n}$ is just the span of the radius $n$ coronas. Thus, the first summand is just the radial functions $L_p(\cL \Fr^\rad) \cong L_p(I;dm_r; \cH_0)$ and our technique recovers a similar spectral decomposition and orthogonal basis for functions that depend on both the radius and the starting $m$-subword. We have that the direct summand in \eqref{eq:directDecomposition} satisfies the following
\begin{enumerate}
    \item Each of the summands $L_p(I, dm_r; \cH_m)$ is a right $\cL \Fr^\rad$-module.
    \item If $a: \Fr \to \CC$ is a radial symbol, meaning that there is a profile $\mathring{a} \colon \NN_0 \to \CC$ such that $a(\gamma) = \mathring{a}(|\gamma|)$, then the Fourier multiplier $T_a$ leaves the summands on \eqref{eq:directDecomposition} invariant, ie: $T_a \big[ L_2 (I) \otimes \cH_m \big] \subseteq L_2 (I) \otimes \cH_m$.
    
	\item Furthermore, the restriction $T_a|_{L_2 (I) \otimes \cH_m }$ takes the form $R_{a}^{(m)} \otimes \id_{\cH_m}$, with the map $R_{a}^{(m)}$ satisfying
	\begin{equation}
	   \label{eq:SpectralQmultiplier}	
       R_{a}^{(m)}( Q_{m,n} ) = \mathring{a}(n) \, Q_{m,n}.
    \end{equation}
\end{enumerate}
Thus, each radial Fourier multiplier is composed of a direct sum of spectral multipliers with respect to different bases in each of the $m$-numerated blocks of the martingale difference on $L_2(\partial \Fr)$. We can summarize this in the following technical proposition. Recall that a family of operators
\(
    \EuScript{T}
    \subset
    B (L_p(\Omega,d\nu))
\)
is \(R\)-bounded if there exists a constant \(C\geq 0\) such that
\[
    \left\|
        \bigg(
            \sum_{j = 1}^N
            |T_{j}(f_j)|^2
        \bigg)^{\frac12}
    \right\|_{L_p(\Omega,d\nu)}
    \leq
    C \, 
    \left\|
        \bigg(
            \sum_{j = 1} |f_j|^2
        \bigg)^{\frac12}
    \right\|_{L_p(\Omega,d\nu)},
\]
for any finite collection with $T_j \in \EuScript{T}$, admitting thus potential repetitions. We denote the infimum of all such constants \(C > 0\) by \( \cR_p(\EuScript{T})\). Notice as well that our definition of the $R$-boundedness constant is not standard, since the inequality above is normally defined via randomization with Rademacher variables, see \cite[Definition 8.1.1(1), p.~165]{Hyt2017BookII}. Nevertheless, for $L_p$-spaces, the two definitions are equivalent up to constants due to Khinchin's inequality. 

\begin{lproposition}
    \label{Prop_Rbounded_kernels}
    Let \(r \in \{2, 3, \dots, \} \), \(1 \leq p <\infty\) and $a(\gamma) = \mathring{a}(|\gamma|)$ be a bounded radial function on the free group. For each $m \ge 0$, there exists a bounded operator $R_{a}^{(m)}: L_2(I,dm_r)\to L_2(I,dm_r)$ such that the restriction of $\widetilde{T}_a = \Phi T_a \Phi^{-1}$ to $L_2(I,dm_r;\cH_m)$ equals $R_{a}^{(m)} \otimes \id_{\cH_m}$. Moreover, the following assertions hold.
    \begin{enumerate}[label={\rm(\roman*)}]
        \item We have that $\displaystyle \big\| T_a: L_p(\cLFr; \Erad^c) \to L_p(\cLFr; \Erad^c) \big\|
            \; = \; \cR_p\big\{ R_{a}^{(m)}  : m \geq 0 \big\}$.
        \item Furthermore if $\mathring{a}$ is finitely supported, each of the operators $R_{a}^{(m)}$ have (scalar) kernel representations of the form
        \[
            \begin{split}
                R_{a}^{(m)}& (h)(s) = \int_I K_{a}^{(m)}(s,t) \, h(t) \, dm_r(t),
                \qquad h \in L_p(I,dm_r), 
            \end{split}
        \]
        where 
        \[
            K_{a}^{(m)}(s,t) = \sum_{n \geq m} \mathring{a}(n) \, \frac{\, Q_{m,n}(s) \, \overline{Q_{m,n}(t)}}{\|Q_{m,n}\|_{2}^2}.
        \]    
        This kernel identity is also true for general symbols $\mathring{a}$ if the sums above are interpreted in the distributional sense.
    \end{enumerate}
\end{lproposition}

In the case of ball multipliers with $\mathring{a}_N(n) =\chi_{\{0,\dots, N\}}(n)$, we may carry out a delicate analysis of the spectral multipliers $(R_{N}^{m})_{m \geq 0}$ in Section~\ref{subsec: generalization of Pollard's arguments}, by means of a decomposition into an operator behaving like a weighted Hilbert transform and another that behaves like a regular integral operator. This allows us to perform the desired weighted $L_2$ analysis of $R_{\varphi,i}$ using singular integral theory with $A_p$ weights. 

The remainder of the paper is organized as follows. Section \ref{Preliminaries} reviews the main background material on free groups needed throughout the paper, including the radial subalgebra $\LFrrad$, the Poisson boundary and the Poisson kernel, as well as the necessary Plancherel theory. In Section \ref{sec: kernels}, we establish a decomposition of $L_2(\cLFr)$ into common eigenspaces for arbitrary radial Fourier multipliers, under which these multipliers act by scalar-valued kernels at each level. We then derive an explicit formula for the kernel of the ball multiplier, Lemma~\ref{AN-CN}, while in Section \ref{sec: proof of main theorem}, we decompose the ball multiplier kernel into regular and singular components, estimate each part using the theory of singular integrals and special functions, and thereby complete the proof of the main results. Lastly, in Section \ref{sec: q-balls} we show the non-uniform bounded of the $q$-balls.

\section{Preliminaries}\label{Preliminaries}

%

\subsection{Free groups and the radial subalgebra}\label{subsec: radial subalgebra definition}
Let $\Fr$ be the free group on $r$ generators, with $r \in \{2, 3, \cdots \}$. As usual, we will denote by $\lambda \colon \Fr\to B(\ell_2\Fr)$ the left-regular representation on its space of $\ell_2$-summable sequences $\ell_2(\Fr)$. Meanwhile, its associated left regular von Neumann algebra is set to be 
\[
    \cL \Fr := \big\{ \lambda_\gamma : \gamma \in \Gamma \big\}'' \subset B(\ell_2 \Fr).
\]
We will identify the group algebra $\CFr$, that is the set of complex-valued finitely supported functions in $\Fr$, with the subalgebra of $\cL \Fr$ spanned by a finite number of unitaries $\lambda_\gamma$. The von Neumann algebra $\cL \Fr$ is a natural example of a $\mathrm{II}_1$ factor and thus comes equipped with a unique tracial state $\tau(x)=\langle x\delta_e,\delta_e\rangle$, which allows us to define their noncommutative $L_p$-spaces $L_p(\cLFr)$, $p\in [1,\infty)$, as the completions of $\cLFr$ with respect to the norms
\begin{equation*}
    \|x\|_{L_p(\cLFr)} :=
        \tau(|x|^p)^{\frac1{p}} ,\quad 1\leq p<\infty.
\end{equation*}
Meanwhile $L_\infty (\cLFr)$ would denote $\cLFr$. The construction of noncommutative $L_p$-spaces is standard by now \cite{PX2003, GoldsteinLabuschagne2025, Hiai2021}. We also recall that it is also possible to describe them as subsets of linear (unbounded) $\tau$-measurable operators $L_0(\tau)$, see \cite{Terp1981}. Given a function $f: \Fr \to \CC$ and a unitary representation $\pi: \Fr \to B(\cH)$, we will denote the linear extension of $\pi$ by
\[
    \pi(f) := \sum_{\gamma\in\Fr}f(\gamma)\pi(\gamma).
\]
This is a well-defined and bounded operator if $f \in \ell_1(\Fr)$ but, in the case of $\pi = \lambda$, much more is true. Namely, since $\tau$ satisfies the Plancherel identity---ie $\| \lambda(f) \|_{L_2(\cL \Fr)} = \| f \|_2$---the expression $\lambda(f)$ gives a well defined unbounded operator for every $f \in \ell_2(\Fr)$.  

We will also denote by $|\cdot|$ the canonical word-length function. A function $f$ on $\Fr$ is called \emph{radial} if $f(\gamma)$ depends only on $|\gamma|$, that is, $f(\gamma)=f(\gamma')$ as soon as $|\gamma|=|\gamma'|$. Clearly, when that is the case, there exists a \emph{radial profile} $\mathring{f}: \NN_0 \to \CC$ such that $f(\gamma) = \mathring{f}(|\gamma|)$. As discussed in the introduction, the radial (finitely) supported functions are $\ast$-algebra generated by $\mu_1$, and as such, from an Abelian $\ast$-algebra $\CFrrad \coloneqq \{f\in\CFr\,\colon\, f\text{ is radial} \}$. Its weak-$\ast$ closure inside the group von neumann algebra would be denoted by $\LFrrad \subset \LFr$.


The radius $n$ corona on the free group would be denoted by $\mathbf{S}_n = \{ \gamma \in \Fr : |\gamma| = n \}$ while its cardinality is 
\begin{equation}\label{cn-def}
	c(n) \coloneqq |\mathbf S_n|= \left\{
	\begin{aligned}
		&\ {2r(2r-1)^{n-1}}, & \hbox{ $n\geq1$}, \\
		&\ 1, & \hbox{ $n=0$}.
	\end{aligned}
	\right.
\end{equation}
As in \eqref{eq:muNdefs}, $\mu_n$ will be the $\ell_1$-normalized characteristic function of the corona of radius $n$ in $\Fr$. The family $(\mu_n)_{n\geq 0}$ clearly constitutes a linear basis of the space $\CFrrad$. The well-known convolution identity \eqref{mun} follows by noting that each element $\gamma$ of length $|\gamma| = n$ has exactly 1 adjacent element on the Cayley graph at distance $n-1$ and $2r -1$ elements at distance $n+1$. Since $\cL \Fr^\rad$ is the weak-$\ast$ closure of $\CC[\Fr]^\rad$, it is an Abelian von Neumann algebra and, by spectral theory, we have a trace-preserving $\ast$-isomorphism
\[
    \Gamma\,\colon\,(\LFrrad,\tau) \longrightarrow L_\infty(I,m_r),
\]
where $I := \sigma(\lambda(\mu_1))$ is the spectrum of the self-adjoint generator $\mu_1$ and $dm_r$ is the restriction of the Plancherel trace. Both $I$ and $dm_r$ are well known by now
\begin{equation*}
    I=\Big[-\frac{\sqrt{2r-1}}{r},\frac{\sqrt{2r-1}}{r}\Big],\qquad dm_r(s)=\frac{1}{\pi}\frac{\sqrt{2r-1-r^2s^2}}{1-s^2}ds,
\end{equation*}
see \cite{FP1982}. In this context, the Gel'fand transform $\Gamma$ is sometimes referred to as the \emph{spherical transform}, and $dm_r$ as the \emph{spherical Plancherel measure}. Observe that the spherical transform has an explicit expression as follows. Let $P_n: I \to \CC$ be the image under $\Gamma$ of $\lambda(\mu_n)$. Then, using that $\Gamma$ is a homomorphism, applying it to \eqref{mun} and noticing that $P_1(s) =s$ gives the following recursive identities for $P_n$
\begin{equation}
    \label{Pn}
    \begin{cases}
        P_0(s)=1,\; P_1(s)=s, & \\[7pt]
        \displaystyle{P_{n+1}(s)=\frac{2r}{2r-1}sP_n(s)-\frac{1}{2r-1}P_{n-1}(s)}, & \text{ when } n \geq 1.
    \end{cases}
\end{equation}
The $(P_n)_{n \geq 0}$ are therefore a sequence of orthonormal polynomials with respect to the measure $d m_r$.

We note that the following explicit expression for $P_n(s)$ can be derived from the recurrence above. Fix $s \in I$. The characteristic equation associated with the recursive sequence \eqref{Pn} is
\[
    \sigma^2-\frac{2r \, s}{2r-1} \, \sigma + \frac{1}{2r-1} = 0. 
\]
Their two roots, depending on $s$, would be denoted by $\sigma(s)$ and $\overline{\sigma}(s)$, where 
\begin{equation}
    \label{sigma-def}
    \sigma(s)
    =
    \frac{rs-\mathrm{i}\sqrt{2r-1-r^2s^2}}{2r-1}.
\end{equation}
These two roots are different when $s$ is in the interior of the interval $I$. Hence, using that $\sigma(s)^n$ and $\overline{\sigma}(s)^n$ satisfy the recurrence law in \eqref{Pn} and solving for the initial data $P_0(s)=1$ and $P_1(s)=s$, it follows that
\begin{equation}
    \label{Pn-sigma}
    P_n(s)
    =
    \frac{s-\sigma(s)}{\overline\sigma(s)-\sigma(s)}\,\overline\sigma(s)^n
    +
    \frac{\overline\sigma(s)-s}{\overline\sigma(s)-\sigma(s)}\,\sigma(s)^n,
    \qquad s\in \bigg(-\frac{\sqrt{2r-1}}{r},\frac{\sqrt{2r-1}}{r} \bigg).
\end{equation}

Let us obtain an explicit formula for the unique trace-preserving conditional expectation $\EE_\rad$ onto $\cL \Fr^\rad$. Notice that it naturally extends to an orthogonal projection from $L_2(\cL \Fr)$ onto $L_2(\cL \Fr^\rad) \subset L_2(\cL \Fr)$. This projection is the adjoint of the inclusion $L_2(\cL \Fr^\rad) \into L_2(\cL \Fr)$. Thus, using that 
\[
    \tau \big( \lambda_\gamma^\ast \, \lambda(\mu_n) \big) = \frac1{c(n)} \, \delta_{\{|\gamma| = n\}}
\]
and taking adjoins, it is obtained that $\EE^\rad$ is given by averaging the frequencies of each corona $\mathbf{S}_n$, that is
\[
    \EE^\rad [\lambda_\gamma] = \lambda( \mu_n), \quad \text{ where } \; |\gamma| = n.
\]
  



\subsection{Poisson boundary and principal series representations}\label{subsec: poisson boundary}

Let $\partial \Fr$ be the set of all infinite reduced words $\omega = a_{j_0}a_{j_1}a_{j_2}\dots$ without cancellation, where each $a_{j_n}$ is a generator or the inverse of a generator. For each $\gamma\in \Fr$, we define $E(\gamma)\subseteq \partial\Fr$ as the set of infinite words starting with $\gamma$. We equip \(\partial\Fr\) with its natural topology, the coarser topology that makes all $E(\gamma)$ open. This makes $\partial \Fr$ into a totally disconnected space. We denote by
\[
	\cB =\sigma\{E(\gamma ):\gamma \in\Fr\}
\]
the associated Borel \(\sigma\)-algebra. Note also that for each integer $n\geq 0$, the family $\{E(\gamma)\,\colon\,|\gamma|=n\}$ is a partition of $\partial\Fr$ into $c(n)$ pieces. Thus, there is a unique Radon probability measure \(\nu\) on \((\partial\Fr,\cB)\) characterized by 
\[
    \nu(E(\gamma))=\frac{1}{c(|\gamma|)}
    \quad \text{ for every cylinder set } E(\gamma).
\]
The resulting Borel probability space \((\partial\Fr,\cB,\nu)\) is the \emph{Poisson boundary} of \(\Fr\), see \cite{FP1982}. The measurable cylinders $\cF_n := \{ E(\gamma) : |\gamma| = n\}$ are finite (atomic) $\sigma$-algebras and form a filtration $\cF_n \subset \cF_{n+1}$ on the Poisson boundary that approximates its Borel $\sigma$-algebra \(\cB\). Their associated conditional expectations and martingale differences would be denoted by
\begin{equation}
    \label{eq:def poisson martingale diff}
    \EE_n \colon L_\infty(\partial \Fr,\cB,d\nu) \to L_\infty(\partial\Fr,\cF_n,d\nu),\quad \Delta_n = \EE_n - \EE_{n-1} \ \text{(with } \EE_{-1}\coloneqq 0),
\end{equation}
while, as mentioned in the introduction, the $L_2$-range of $\Delta_n$, will be denoted by 
\[
    \cH_n := \text{ran}(\Delta_n) = L_2(\partial\Fr,\cF_n,d\nu) \ominus L_2(\partial\Fr,\cF_{n-1},d\nu).
\]

The free group \(\Fr\) acts on its boundary \(\partial\Fr\) by left multiplication: for \(g\in\Fr\) and \(\omega\in\partial\Fr\), we define \(g\cdot\omega\) to be the infinite reduced word obtained by concatenating \(g\) and \(\omega\), followed by reduction.  This action does not preserve the probability measure \(\nu\), however, for each \(\gamma\in\Fr\), the push-forward measure \(\nu_\gamma(E):=\nu(\gamma^{-1}E)\) satisfies that \(\nu_\gamma\ll\nu\). Thus its Radon--Nikodym derivative defines a multiplicative cocycle of the action $\Fr \acts \partial \Fr$ given by
\[
    P(\gamma,\omega):=\frac{d\nu_\gamma}{d\nu}(\omega),\qquad \gamma\in\Fr,\ \omega\in\partial\Fr.
\]
The $P(\gamma, \omega)$ would be referred to as the \emph{Poisson kernel}. This name comes from a well-known geometric interpretation where $P(\gamma, \omega)$ is the probability density that a random walk starting at $\gamma \in \Fr$, and proceeding with homogeneous probability along the Cayley graph, would diverge to the boundary in the direction given by $\omega$. The Poisson kernel plays a central role in the representation theory of free groups, although we will only make use of some key results. First, recall that if $\gamma \in \Fr$ and $\omega \in \overline{\FF}_r = \Fr \cup \partial \Fr$, then we will denote by $\gamma \wedge \omega \in \Fr$ the largest common initial subword of the 
(reduced) words corresponding to $\gamma$ and $\omega$. Since $\gamma$ is finite, it follows that $\gamma \wedge \omega \in \Fr$. In geometric terms, one can see $\gamma \wedge \omega$ as the last common vertex of the Cayley graph geodesic joining $e$ with $\gamma$ and the geodesic representing $\omega$. The following explicit formula:
\begin{equation}
    \label{eq:poisson-kernel-gromov}
    P(\gamma,\omega)
    =
    (2r-1)^{2 \, |\gamma \wedge \omega| - |\gamma| },
\end{equation}
follows by computing $\nu(\gamma \, E(\eta))$ for every $\eta \in \Fr$, see \cite{FP1982}. Throughout the paper, for each \(\gamma\in\Fr\), let  \(\gamma^{(n)}\) denote the prefix of \(\gamma\) of length \(n\). The following telescoping rearrangement of the formula above for complex powers of $P(\gamma, \omega)$  would  play an important role in subsequent computations. 
\begin{lemma}[{\cite[Ch. 3, Proposition 2.5]{FP1983}}] 
    \label{Possion-lem} 
    Let \(\gamma \in \Fr\), \(\omega\in \partial \Fr\), and \(z\in \CC\), with \(|\gamma|=n\). Then
    \begin{equation*}
        \label{eq:poisson-kernel-cylinders}
        P(\gamma,\omega)^z
        =
        (2r-1)^{nz}\chi_{E(\gamma)}(\omega)
        +
        \sum_{k=0}^{n-1} (2r-1)^{(2k-n)z}
        \left[
            \chi_{E(\gamma^{(k)})}(\omega)
            -
            \chi_{E(\gamma^{(k+1)})}(\omega)
        \right].
    \end{equation*}
\end{lemma}

As mentioned in the introduction, the maps $\pi_s$ in \eqref{def:pseries_rep} constitute unitary representations of $\Fr$ on $L_2(\partial \Fr)$. This follows easily from the definition of the measure cocycle. The choice of the parameter $t$ in the exponent is taken because, due to \eqref{eq:poisson-kernel-gromov}, the imaginary powers of $P(\gamma, \omega)$ have period $2 \pi/\log(2 r - 1)$ and the map sending $t$ to $-t$ in \eqref{def:pseries_rep} giving representations that are equivalent after conjugation. This parametrized family of representations $(\pi_s ) _{s\in I}$ constitutes the \emph{principal series} of \(\Fr\), in analogy with the principal series for semisimple Lie groups.
Let $\varphi_s(\gamma) = \big(\Gamma \circ \EE_\rad(\lambda_\gamma)\big)(s)$, for $s \in I$, be the normalized positive definite function on $\Fr$ coming from the state $\phi_s: C_\lambda^\ast \Fr \to \CC$, defined in \eqref{def: state phi_s}. The following lemma establishes the equality between the GNS construction of $\varphi_s$ and the principal series. 

\begin{lemma}[{\cite{FP1983}}]
    \label{phi-rep-con}
    Let $s\in I$ and \(\mathds{1}(\omega)\equiv 1\) be the constant function in \(L_2(\partial\Fr)\) and let $F \colon \big[0, \tfrac{\pi}{\log(2r-1)}\big] \to I$ be the bijection given by
    \begin{equation}
        \label{eq: defining F}
        F(t) := \frac1{2 r}\big[(2r - 1)^{\frac12 + \imagi t} + (2r - 1)^{\frac12 - \imagi t}\big].
    \end{equation}
    Then, it holds that
    \[
        \varphi_s(\gamma)
        =
        \langle \pi_s(\gamma)\mathds{1},\mathds{1}\rangle_{L_2(\partial\Fr)}
        =
        \int_{\partial\Fr} P(\gamma,\omega)^{\frac12 + \imagi \, t}\,d\nu(\omega),
        \; \text{ where } \; 
        s = F(t).
    \]
\end{lemma}
An alternative description of the spherical functions $(\varphi_s)_{s \in I}$ states that those are precisely the extremal points of the convex set of normalized positive definite functions on $\Fr$ that are also radial and such that the map $\lambda_\gamma \mapsto \varphi(\gamma)$ extends boundedly to the reduced $C^\ast$-algebra $C_\lambda^\ast \Fr$. This last condition is equivalent to requiring the associated GNS representation of $\varphi$ to be weakly contained in $\lambda$; see \cite[Theorem~3.4.4 and Section~18.1]{Di1977}. Observe that, if we where to drop the last condition, the resulting convex body of radial and positive definite functions becomes larger as per \cite[Theorem~3.1]{HaagerupKnudby2015}. 

\subsection{Plancherel Theory}\label{subsec: Phi definition}
The aforementioned theory of principal series representations provides a Plancherel theory on free groups, analogous to the one on type I semisimple Lie groups. This yields an \(L_2\)-isomorphism $\ell_2( \Fr) \cong L_2\big(I; L_2(\partial\mathbb F_r)\big)$, decomposing the left regular representation \(\lambda\) into a direct integral of principal series representations. We refer to \cite{FP1983} for more details on this Plancherel theory, and to \cite{Di1977,Ta2002} for the theory of direct integrals. Consider the map \(\Phi: \CC[\Fr] \to L_0(I;L_2(\partial \Fr))\) defined by
\begin{equation}
    \label{phi-def}
	\Phi(\lambda(f))(s)= {\pi}_s(f)\mathds{1},\quad f\in \CC[\Fr],\ s\in I.
\end{equation}
By the following result $\Phi$ explicitly intertwines $\lambda$ and $\pi_s$.


\begin{lemma}[{\cite[Chapter 5]{FP1983}}]\label{Plancherelthm}
The map
\[
    U:\ell_2(\Fr) \longrightarrow  L_2(I,dm_r;L_2(\partial \Fr)),\quad \delta_\gamma \mapsto (s\mapsto \pi_s (\gamma) \mathds 1 )
\]
defines a unitary operator, and we have
    \[
    \lambda(\gamma ) = U^*\bigg(\int_{I}^{\oplus} {\pi}_{s}(\gamma ) \, dm_r(s)\bigg) U,\quad \gamma \in \Fr.
    \]
Moreover, for each \(f\in \CC [\Fr]\), we have
\begin{equation*}\label{trace-pri-se}
    \tau(\lambda(f))=\int_I \phi_s(\lambda(f))dm_r(s)=\int_I\langle {\pi}_s(f)\mathds{1},\mathds{1}\rangle_{L_2 (\partial \Fr)} dm_r(s).
    \end{equation*} 
In particular, the map \(\Phi\) extends to an isometric isomorphism from \(L_2(\cLFr)\) to \(L_2(I;L_2(\partial\Fr))\).
\end{lemma}

\section{Mixed-norm spaces and integral kernels of radial multipliers}\label{sec: kernels}

Mixed-norm spaces associated with conditional expectations have already been studied in full generality in \cite{JP2010}. Let us recall some basic facts of the definition on the particular case of our interest. The expression for the norm has already been defined in the Introduction; see \eqref{mixed-norms}. Regarding a description of the spaces $L_p(\cL \Fr; \EE_\rad^c)$ as sets, we will explicitly describe them as follows
\begin{align*}
    L_p(\cL \Fr; \EE_\rad^c) := \big\{ x \in L_p(\cL \Fr) : \| x \|_{L_p(\cL \Fr; \EE_\rad^c)} < \infty \big\} & \text{ when } 1 \leq p < 2\\[6pt]
    L_p(\cL \Fr; \EE_\rad^c) := \big\{ x \in L_2(\cL \Fr) : \| x \|_{L_p(\cL \Fr; \EE_\rad^c)} < \infty \big\} & \text{ when } 2 \leq p \leq \infty.
\end{align*}
The above definition gives complete Banach spaces. To see that, we just need to notice that we have the following a priori inequalities between mixed norms and $L_p$-norms. First, when $2 \leq p$, we have 
\[
    \| x \|_{L_p(\cL \Fr; \EE_\rad^c)} = \big\| \EE_\rad[x^*x] \big\|_{\frac{p}{2}}^\frac1{2} \geq \big\| \EE_\rad[x^\ast x] \big\|_{1}^{\frac12} = \| x \|_2.
\]
For $1 \leq p < 2$, we can use the operator concavity inequality in \cite[Theorem~2.1]{Choi1974}, obtaining the following:
\[
    \| x \|_{L_p(\cL \Fr; \EE_\rad^c)}^p = \big\| \EE_\rad[x^*x]^\frac{p}{2} \big\|_{1} \geq \big\| \EE_\rad\big[|x|^p\big] \big\|_1 = \| x \|_{p}^p.
\]
Recall that the mixed-norm spaces inherit an interpolation compatible operator space structure from $L_p$, see \cite{JP2010} for the details. An immediate consequence of the definition above is that, for $1 \leq p < \infty$, the group algebra $\CC[\Fr]$ is dense inside the mixed-norm spaces $L_p(\cL \Fr; \EE_\rad^c)$. The same inequalities and definitions work in the row case.  

\begin{proposition}
    \label{Phi}
	Let \(1\leq p<\infty\) and \(\Phi\) be the linear map defined in \eqref{phi-def}. The following holds
    \begin{enumerate}[label={(\roman*)}, ref={(\roman*)}, font=\normalfont]
        \item \label{itm:Phi1}
        The map $\Phi$ is right $\cL \Fr^\rad$-modular, ie:
        \[
            \Phi( x \, y )(s) = \Phi(x)(s) \cdot (\Gamma y)(s), \quad \text{ for every } \; x \in L_p(\cL; \EE_\rad), \, y \in \cL \Fr^\rad.
        \]
        \item \label{itm:Phi2}
        $\Phi$ extends to an isometric isomorphism between Banach spaces
	    \begin{equation*}
		      \Phi :{L_p}(\cL(\Fr);\EE_{\rad}^c) \longrightarrow {L_p}(I,m_r;L_2(\partial \Fr)).
	    \end{equation*}
    \end{enumerate}
    The analogous assertions also holds for \(L_p(\cL(\Fr);\EE_{\rad}^r)\).
\end{proposition}
\begin{proof}
    For \ref{itm:Phi1}, we will start by showing the elementary identity 
    \begin{equation}
        \label{eq: elementary modular id}
        \pi_s(\mu_1) \, \mathds{1} = s \, \mathds{1}.
    \end{equation}
    A computation shows that
    \begin{eqnarray*}
         \big\| \big( \pi_s(\mu) - s \, \id \big) \1 \big\|_{L_2(\partial \Fr)}^2
            & = & \big\langle \pi_s(\mu_1 \ast \mu_1) \1, \1 \big\rangle + s^2 
                  - 2 \, \Re \big\langle \pi_s(\mu_1) \1, \1 \big\rangle \\
            & = & \phi_s( \lambda(\mu_1 \ast \mu) ) + s^2 - 2 \, \phi_s( \lambda(\mu_1) ) \; = \; s^2 + s^2 - 2 \, s^2 \; = \; 0.
    \end{eqnarray*}
    In the above equations, we have used that $\varphi_s(\gamma)$ coincides with $\langle \pi_s(\gamma) \1, \1\rangle$ and that 
    $\varphi_s$ induces a state on $C_\lambda^\ast \Fr$ coinciding with spectral evaluation at $s \in I$. 
    It follows that $\Phi( x \, \lambda(\mu_1))(s) = \Phi(x)(s) \, s$, for every $x \in \CC[\Fr]$. Similarly, iteratively multiplying by powers of $\lambda(\mu_1)$ we get \eqref{eq: elementary modular id} by $y \in \CC[\Fr]^\rad$. By norm density of $\CC[\Fr]$ inside $L_p(\cL \Fr; \EE_\rad)$ we can extend it to general $x$. For $y$, we use that the right action of $\cL \Fr$ on the mixed-norm space $L_p(\cL \Fr; \EE_\rad^c)$ is continuous with respect to the strong operator topology. 
    
	For \ref{itm:Phi2}, by a density argument, it suffices to prove that \(\Phi\) is an isometry on \(\CC[\Fr]\). Let \(f\in \CC[\Fr]\). Using Lemma \ref{Plancherelthm},  we obtain
	\[
	   \left \| \lambda(f) \right \|_{{L_p}(\cL(\Fr);\EE_{\rad}^c)}^p 
       = \tau\big( \EE_{\rad}[\lambda(f)^*\lambda(f)]^{\frac{p}{2}}\big) 
       =\int_I \phi_s\big(\EE_{\rad}[\lambda(f)^*\lambda(f)]^{\frac{p}{2}}\big) \, dm_r(s).
	\]
	Observe that, since $\varphi_s$ is radial, $\phi_s \circ \EE_\rad = \phi_s$. Now, using that $\phi_s(\lambda(f))$ equals $\langle \pi_s(f) \mathds{1}, \mathds{1} \rangle$ from  Lemma~\ref{phi-rep-con} gives that
	\begin{align*}
		\left \| \lambda(f) \right \|_{{L_p}(\cL(\Fr);\EE_{\rad}^c)}^p 
		      & = \int_I \left |\langle {\pi}_s(f^**f)\mathds{1},\mathds{1}\rangle \right |^{\frac{p}{2}}dm_r(s) 
              = \int_I\left| \langle {\pi}_s(f) \mathds{1}, {\pi}_s(f)\mathds{1}\rangle \right|^{\frac{p}{2}}dm_r(s)\\
		      & = \int_I \left\| {\pi}_s(f) \mathds{1} \right\|_{L^2(\partial \Fr)}^p dm_r(s)  = \left \|\Phi(\lambda(f)) \right\|_{{L_p}(I,dm_r;L^2(\partial \Fr))}^p. 
	\end{align*} 
    To see that it is an isomorphism, it remains to show that $\Phi$ is surjective. Assume it is not, fix $p'$ as the conjugate exponent satisfying $\tfrac1{p} + \tfrac1{p'} = 1$, there is a nonzero $h \in L_{p'}(I, dm_r; L_2(\partial \Fr))$ annihilating the range of $\Phi$. In particular, $h$ satisfies that
    \[
        \int_I \langle h(s), \Phi(x)(s) \rangle_{L_2(\partial \Fr)} \, d m_r(s) = 0, \quad \text{ for all } x \in \CC[\Fr],
    \]
    using the modular identity in \ref{itm:Phi1}, and the density of polynomials on $C(I)$, we get that $h(s)  =0$ for almost every $s \in I$, giving a contradiction and thus concluding the proof.
\end{proof}
In the proof of point \ref{itm:Phi2}, we have chosen to give a proof of the key identity \eqref{eq: elementary modular id} using only the equality between the state $\phi_s$ and the vector state associated with the representations $\pi_s$. However, a direct proof using explicit formulas for $\pi_s$ is also available and would be provided in Lemma \ref{prop1} and Remark \ref{Pn-def2}.


\subsection{The kernels of radial multiplier operators}\label{subsec: kernel decomposition of radial multipliers}
Now, we will obtain the explicit kernel formulas for general radial multiplier operators.  These expressions serve as the basis for the study of partial sum operators.

In the introduction, we defined sub-$\sigma$ algebras $\cA_m$ of $\cP(\mathbf{S}_n)$ given by the equivalence classes in $\mathbf{S}_m$ that have the same starting $m$-subword.  When $n \leq m$ this $\sigma$-algebra is generated by the singletons $\{ \gamma : |\gamma|= n\}$, while, for $n \geq m$ we have an atomic $\sigma$-algebra generated by the $c(m)$ arcs
\[
    S(m,\gamma)\coloneqq \{\gamma' : |\gamma ' |=|\gamma |,\ {\gamma'} ^{(m)}= \gamma^{(m)}   \}, \quad |\gamma|>m,
\]
where, we recall, \(\gamma^{(m)}\) denotes the starting $m$-subword of \(\gamma\). Observe as well that \(|S(m,\gamma)|=\tfrac{c(|\gamma|)}{c(m)}\). In a slight abuse of notation, we will denote by $\cA_m$ both the $\sigma$-algebra of $\mathbf{S}_n$ described above, as well as the $\sigma$-algebra of the full group $\Fr$ generated by all of the arcs $S(m,\gamma)$ of any length. Equivalent, the latter can be seen as the equivalence classes of elements that have the same length and the same starting $m$-subword. Note that $\cA _m \subset \cA_{m+1}$ for all $m$ and thus, they form a filtration whose final $\sigma$-algebra is $\cP(\Fr)$. Denote the associated conditional expectations by $\cE_m: \ell_2 (\Fr) \to \ell_2 (\Fr;\cA_m)$ and their martingale differences by $\cD_m = \cE_m - \cE_{m-1}$ where we set \(\cE_{-1}=0\). We have the following explicit formula
\begin{equation}
    \label{Em-delta}
    \cE _m (\delta_\gamma) 
    = 
    \begin{cases}
	   \delta_\gamma & \text{ when } |\gamma|\leq m , \\[5pt]
	   \displaystyle{\frac1{|S(m,\gamma)|} \chi_{S(m,\gamma )}} & \text{ when } |\gamma|>m.
    \end{cases}
\end{equation}
The spaces of $\cA_m$-measurable functions admit natural interpretations. For instance, when $m=0$, we have that $\ell_2 (\Fr, \cA _0)$ consists of all radial functions in $\ell_2 (\Fr)$, and hence $\EE_\rad (\lambda(f)) = \lambda (\cE_0 (f))$ for all $f\in \mathbb C [\Fr]$. Similarly, for larger $m$, we get spaces of functions depending on the radial length and the $c(m)$ possible starting subwords. With this definitions at hand, we may obtain the following description of invariant subspaces under radial Fourier multiplier maps. Denote
\[
    V_{m,n}
    \coloneqq
    \ell_2(\mathbf{S}_n;\cA_m) \ominus \ell_2(\mathbf{S}_n;\cA_{m-1})
    =
    \operatorname{span}
    \left\{ 
        \cD_m(\delta_\gamma ):  \gamma\in \mathbf S_n
    \right\},
    \quad m,n\in \NN.
\]
We have the following lemma.
\begin{lemma}
    \label{l2-Vmn} Let \(r\geq 2\). 
    The following statements hold.
    \begin{enumerate}[label={\rm (\roman*)}, ref={(\roman*)}]
        \item \label{itm:l2-Vmn1}
        For any $m\in \mathbb N$, $n\geq m$ and $\gamma \in \mathbf S_n$, the function $\cD _m (\delta_\gamma ) $ is supported in $\mathbf S_n$. In particular, the subspaces $(V_{m,n})_{n\geq m}$ are orthogonal in $\ell_2 (\Fr)$ and 
        \[\ell_2 (\Fr) = \bigoplus_{m=0}^\infty \bigoplus_{n= m}^\infty V_{m,n} .\]
        
        \item \label{itm:l2-Vmn2} 
        Let  $a:\Fr \to \mathbb C$ be a radial function. Then $T_a \big[\lambda (V_{m,n})\big] \subseteq \lambda (V_{m,n})$, and  
        \[
            T_a(\lambda (f))
		      =
		      \mathring{a}(n)
		      \lambda (f),\quad f\in V_{m,n}.
        \]
    \end{enumerate}
\end{lemma}
\begin{proof} For any $m\in \mathbb N$, $n\geq m$ and $\gamma \in \mathbf S_n$,  \eqref{Em-delta} yields
\[
    \cD _m (\delta_\gamma ) =|S(m,\gamma)|^{-1}\chi_{S(m,\gamma)}-|S(m-1,\gamma)|^{-1}\chi_{S(m-1,\gamma)},
\]
where \(S(m,\gamma)\), \(S(m-1,\gamma)\subseteq \mathbf S_n\). Consequently, the function $\cD _m (\delta_\gamma ) $ is supported in $\mathbf S_n$,  and the  subspaces $(V_{m,n})_{n\geq m}$ are orthogonal in $\ell_2 (\Fr)$. The Hilbert space \(\ell_2 (\Fr)\) admits the orthogonal decomposition 
\[
    \ell_2 (\Fr) = \bigoplus_{n=0}^\infty \ell_2(\mathbf{S}_n) = \bigoplus_{n=0}^\infty \; \bigoplus_{m=0}^n V_{m,n},
\]
switching the order of the direct sums proves \ref{itm:l2-Vmn1}. Lastly, for each \(f\in V_{m,n}\), we have \(\operatorname{supp}(f)\subseteq \mathbf S_n\). Since \(a\) is radial, \ref{itm:l2-Vmn2} follows.
\end{proof}
For elements in the subspaces \(\{V_{m,n}\}_{n\geq m}\), we derive a recurrence relation analogous to that given in \eqref{mun}.

\begin{lemma}
    \label{Dmn-re}
    Let \(r \in \{2,3, \dots\}\) and \(m\ge 0\). For each \(n\ge m+1\) and \(\gamma\in\Fr\) with \(|\gamma|\ge n+1\), one has
    \begin{equation}
        \label{relation_Dmmu1}
        \cD_m(\delta_{\gamma^{(n)}})*\mu_1
        =
        \frac{1}{2r}\cD_m(\delta_{\gamma^{(n-1)}})
        +
        \frac{2r-1}{2r}\cD_m(\delta_{\gamma^{(n+1)}}).
    \end{equation}
\end{lemma}

\begin{proof} For the case \(m=0\), it can be directly obtained from \eqref{mun}. For \(m\geq1\) and   \(n\geq m\),
	\begin{equation}\label{Dmn}
			\cD_m(\delta_{\gamma^{(n)}})*\mu_1=\frac{c(m)}{c(n)}\chi_{S(m,\gamma^{(n)})}*\mu_1
			-
			\frac{c(m-1)}{c(n)}\chi_{S(m-1,\gamma^{(n)})}*\mu_1,
	\end{equation}
    where and in the following proof, \(c(n)\) is defined as in \eqref{cn-def}.
	For each \(\rho \in\Fr\), we have
	\begin{equation}\label{chi-mu-1}
		\begin{aligned}
        \chi_{S(m,\gamma^{(n)})}*\mu_1(\rho)&=\sum_{|\eta|=1}	\chi_{S(m,\gamma^{(n)})}(\rho \eta^{-1})\mu_1(\eta)\\
			&=\frac{1}{2r}\chi_{S(m,\gamma^{(n+1)})}(\rho)+\frac{2r-1}{2r}\chi_{S(m,\gamma^{(n-1)})}(\rho).
		\end{aligned}
	\end{equation}
	An analogous equality holds for \(\chi_{S(m-1,\gamma^{(n)})}\).
	Combining \eqref{Dmn} and \eqref{chi-mu-1}, we obtain
	\begin{align*}
		\cD_m(\delta_{\gamma^{(n)}})*\mu_1
            & = \frac{c(m)}{c(n)}\Big[\frac{1}{2r}\chi_{S(m,\gamma^{(n+1)})}+\frac{2r-1}{2r}\chi_{S(m,\gamma^{(n-1)})}\Big]\\
		    & \hspace*{50pt} -\frac{c(m-1)}{c(n)}\Big[\frac{1}{2r}\chi_{S(m-1,\gamma^{(n+1)})}+\frac{2r-1}{2r}\chi_{S(m-1,\gamma^{(n-1)})}\Big]\\
		      & =\frac{1}{2r}\Big[\frac{c(m)}{c(n-1)}\chi_{S(m,\gamma^{(n-1)})}-\frac{c(m-1)}{c(n-1)}\chi_{S(m-1,\gamma^{(n-1)})}\Big]\\
		      & \hspace*{50pt} +\frac{2r-1}{2r}\Big[\frac{c(m)}{c(n+1)}\chi_{S(m,\gamma^{(n+1)})}-\frac{c(m-1)}{c(n+1)}	\chi_{S(m-1,\gamma^{(n+1)})}\Big]. \qedhere
	\end{align*}
\end{proof}
We will show that the space $\Phi(V_{m,n})$ can be identified with the tensor product of a one-dimensional space of functions on the interval $I$ and a martingale difference subspace on the Poisson boundary $\partial\Fr$. Before establishing this identification, we shall use the identity in Lemma~\ref{Possion-lem} to derive an explicit formula for the map $\Phi$ as follows.
\begin{lemma}
    \label{prop1}
	Let \(r \in \{2, 3, \dots\}\). For each \(\gamma  \in \Fr\) with  \(|\gamma |=n\) and \(s\in I\), we have
	\begin{equation}\label{phi-alpha}
		\Phi(\lambda(\gamma ))(s)=\pi_s(\gamma ) \mathds{1}= \sum_{j=0}^n \alpha_{j,n}(s)\chi_{E(\gamma^{(j)})} ,
	\end{equation}
	where the coefficients are given by
	\begin{equation}\label{alpha-jn}
		\alpha_{j,n}(s):= \left\{
		\begin{aligned}
			&\ \sigma(s)^{n}\big(\sigma(s)^{-2j}-\sigma(s)^{-2j+2}\big), & \hbox{ $j\geq1$}, \\
			&\ \sigma(s)^{n}, & \hbox{ $j=0$},
		\end{aligned}
		\right.
	\end{equation}
	and \(\sigma(s)\) is  defined as in \eqref{sigma-def}.
\end{lemma}
\begin{proof}
	Let \(t = F^{-1}(s)\) and denote \(z =\frac{1}{2}+\imagi \, t\). By \eqref{def:pseries_rep} and Lemma~\ref{Possion-lem}, we obtain
	\begin{align*}
		\quad\pi_s(\gamma )\mathds{1} & = P(\gamma ,\cdot)^{\frac{1}{2}+\imagi \, t}\\
			&= (2r-1)^{-n \, z}\mathds{1}+\sum_{j=1}^n\Big[(2r-1)^{-(n-2j) z}-(2r-1)^{-(n-2j+2) z}\Big]\chi_{E(\gamma^{(j)})}\\
			&= (2r-1)^{-z n}\mathds{1} \\ 
            & \qquad +\sum_{j=1}^n \big((2r-1)^{-z}\big)^n\Big[\big((2r-1)^{-z}\big)^{-2j} -\big((2r-1)^{-z}\big)^{-2j+2}\Big]\chi_{E(\gamma^{(j)} )}.
	\end{align*}
	We claim that 
	\begin{equation}\label{sigma_s'}
		(2r-1)^{-z}=\sigma(s).
	\end{equation}
	Then, together with \(E(\gamma^{(0)})=\partial \Fr\), \eqref{phi-alpha} follows. It therefore suffices to prove \eqref{sigma_s'}. By \eqref{eq: defining F}, it follows that 
    \[
        s=F(t)= \frac{\sqrt{2r-1}}{r}\cos{\big(\log{(2r-1)}t\big)}.
    \]
	Hence,
	\[
	    \cos{(\log{(2r-1)}t)}=\frac{rs}{\sqrt{2r-1}} 
        \quad \hbox{and} \quad 
        \sin{(\log{(2r-1)}t)}=\frac{\sqrt{2r-1-r^2s^2}}{\sqrt{2r-1}}.
	\]
	Hence, by \eqref{sigma-def}, we have
	\begin{align*}
			&=\frac {1} {\sqrt{2r-1}}\Big(\frac{rs}{\sqrt{2r-1}}-\mathrm{i}\frac{\sqrt{2r-1-r^2s^2}}{\sqrt{2r-1}}\Big)\\
			&=\frac{rs-\mathrm{i}\sqrt{2r-1-r^2s^2}}{2r-1}=\sigma(s). \qedhere
	\end{align*}
\end{proof}
\begin{remark} \normalfont
    \label{Pn-def2} 
    \
	\begin{enumerate}[label=(\roman*), ref=(\roman*), font=\normalfont]
		\item \label{itm:Pn-def2.1}
        By the fact \(\varphi_s(\gamma)=P_{n}(s)\) with \(|\gamma|=n\) and Lemma~\ref{phi-rep-con}, integrating both sides of \eqref{phi-alpha} with respect to \(\omega\), we obtain
		\[
            P_n(s)=\sum_{j=0}^n\frac{\alpha_{j,n}(s)}{c(j)},\quad s\in I,\quad n\ge 0,
        \]
		where the coefficients $\alpha_{j,n}(s)$ are defined as in Lemma~\ref{prop1}.
		\item \label{itm:Pn-def2.2} 
        For each  \(x\in \cL \Fr^{\rad}\), we have that $\Phi$ extends the spherical transform $\Gamma$, ie 
		\[
            \Phi(x)(s)=\Gamma (x)(s)\mathds1,\quad s\in I.
        \]
		While this follows from the definitions, it can also be seen as a consequence of the computations in Lemma~\ref{prop1}. Indeed, we obtain
		\begin{align*}
			\pi_s(\mu_n)\mathds 1&=\frac{1}{c(n)}\sum_{|\gamma |=n}\pi_s(\gamma )\mathds 1=\frac{1}{c(n)}\sum_{|\gamma |=n}\sum_{j=0}^n\alpha_{j,n}(s)\chi_{E(\gamma^{(j)})}\\
			&=\frac{1}{c(n)}\sum_{j=0}^n\sum_{|\gamma |=n}\alpha_{j,n}(s)\chi_{E(\gamma^{(j)})}=\frac{1}{c(n)}\sum_{j=0}^n\alpha_{j,n}(s) \frac{c(n)}{c(j)}\mathds 1=P_n(s)\mathds1,
		\end{align*}
		where the last equality is due to \ref{itm:Pn-def2.1}
	\end{enumerate}
\end{remark}

We now proceed to explicitly establish the identification of \(V_{m,n}\) with the tensor product of a one-dimensional function space on the interval \(I\) and a martingale difference subspace on the Poisson boundary \(\partial\Fr\). More precisely, for each \(m\ge 0\), the function space on \(I\) corresponds to the deformed Bożejko-Fendler orthogonal system \(\{Q_{m,n}\}_{n\ge m}\) defined by
\begin{equation}
    \label{defQmn}
	Q_{m,n}(s)=\sum_{j=m}^n \frac{\alpha_{j,n}(s)}{c(j)}, \quad s\in I.
\end{equation}

\begin{proposition}
    \label{Phi_QmnHm}
	Let \(r \in \{2, 3, \dots, \} \) and \(m\ge 0\). The following statements hold: 
	\begin{enumerate}[label=(\roman*), ref=(\roman*), font=\normalfont]
		\item \label{itm:Phi_QmnHm1}
        For \(n\geq m\), we have
        \[
            \Phi\big(\lambda(V_{m,n})\big)=Q_{m,n} \otimes \cH_m,
        \]
        where \(\cH_m\) is the range of the martingale difference \(\Delta_m\) on the Poisson boundary defined in \eqref{eq:def poisson martingale diff}. More precisely,
		\begin{equation}
            \label{Dm-gamma-Qmn}
		      \Phi\big(\lambda(\cD_m(\delta_\gamma ))\big)=Q_{m,n}
	        \left[
	           c(m)\chi_{E(\gamma^{(m)})}
	           -
	           c(m-1)\chi_{E(\gamma^{(m-1)})}
	        \right],\quad |\gamma|=n.
		\end{equation}
		
		\item \label{itm:Phi_QmnHm2}
        The family \(\{Q_{m,n}\}_{n\ge m}\) forms a complete orthogonal system on \(L_2(I,dm_r)\) with
		\[
		      \|Q_{m,n}\|_{L_2(I,dm_r)}^2=\frac {1} {c(n)}.
		\]
		\item \label{itm:Phi_QmnHm3}
        The family
		\(
		\{Q_{m,n}\}_{n\ge m}
		\) satisfy the  recurrence relation
		\begin{equation}\label{Qnm}
			Q_{m,n+1}(s)=\frac{2r}{2r-1} \, s \, Q_{m,n}(s)-\frac{1}{2r-1}Q_{m,n-1}(s),\quad n\ge m+1,\ s\in I.
		\end{equation}
	\end{enumerate}
\end{proposition}

\begin{proof}

In the case \(m=0\), \ref{itm:Phi_QmnHm1}, \ref{itm:Phi_QmnHm2} and \ref{itm:Phi_QmnHm3} follow from \eqref{defQmn} and Remark~\ref{Pn-def2} \ref{itm:Pn-def2.1}, \ref{itm:Pn-def2.2}.

We now assume that \(m\ge 1\). Let \(\gamma \in\Fr\) with \(|\gamma |=n\geq m\). 
Applying \(\Phi\) to \(\lambda(\mathcal{D}_m(\delta_\gamma ))\), \eqref{Em-delta} and Lemma \ref{prop1}, we get
\begin{align*}
	\quad \Phi(\lambda(\mathcal{D}_m(\delta_\gamma ))) & = \pi_s\Big(\frac{c(m)}{c(n)}\chi_{S(m,\gamma)}
	-
	\frac{c(m-1)}{c(n)}\chi_{S(m-1,\gamma)}\Big)\mathds{1}\\
	& = \frac{c(m)}{c(n)}\sum_{\eta\in S(m,\gamma)}\sum_{j=0}^n \alpha_{j,n}(s)\chi_{E(\eta^{(j)})}-
	\frac{c(m-1)}{c(n)}\sum_{\rho \in S(m-1,\gamma)}\sum_{j=0}^n \alpha_{j,n}(s)\chi_{E(\rho^{(j)})}\\
	& =  \sum_{j=0}^n\alpha_{j,n}(s) \Big[ \frac{c(m)}{c(n)} \sum_{\eta\in S(m,\gamma)}\chi_{E(\eta^{(j)})}- \frac{c(m-1)}{c(n)} \sum_{\rho \in S(m-1,\gamma)}\chi_{E(\rho^{(j)})}\Big].
\end{align*}
For \(j<m\), for any \(\eta\in S(m,\gamma)\) and \(\rho\in S(m-1,\gamma)\), \(\chi_{E(\eta^{(j)})}=\chi_{E(\rho^{(j)})}=\chi_{E(\gamma^{(j)})}\) hold. From \(|S(m,\gamma)|=c(n)/c(m)\) and \(|S(m-1,\gamma)|=c(n)/c(m-1)\), we directly see that the term inside the bracket is zero. For \(j\ge m\), we have
\[
\sum_{\eta\in S(m,\gamma)}\chi_{E(\eta^{(j)})}
=
\frac{c(n)}{c(j)}\chi_{E(\gamma^{(m)})}
\]
and
\[
\sum_{\rho \in S(m-1,\gamma)}\chi_{E(\rho^{(j)})}
=
\frac{c(n)}{c(j)}\chi_{E(\gamma^{(m-1)})}.
\]
Therefore,
\begin{align*}   
	\Phi\bigl(\lambda(\cD_m(\delta_\gamma ))\bigr)(s)
	&=
	\sum_{j=m}^n
	\frac{\alpha_{j,n}(s)}{c(j)}
	\left[
	c(m)\chi_{E(\gamma^{(m)})}
	-
	c(m-1)\chi_{E(\gamma^{(m-1)})}
	\right] \\
	&=
	Q_{m,n}(s)
	\left[
	c(m)\chi_{E(\gamma^{(m)})}
	-
	c(m-1)\chi_{E(\gamma^{(m-1)})}
	\right].
\end{align*}
The boundary term
\[
c(m)\chi_{E(\gamma^{(m)})}
-
c(m-1)\chi_{E(\gamma^{(m-1)})}
\]
belongs to the martingale difference subspace \(\cH_m\). Thus concluding the proof of \ref{itm:Phi_QmnHm1}.

We now turn to prove the orthogonality of the functions \(\{Q_{m,n}\}_{n\ge m}\).
Fix \(n,n'\ge m\). Choose \(\gamma ,\eta \in\Fr\) such that
\[
|\gamma |=n,\ |\eta|=n'\  \text{and}\  \gamma^{(m)}=\eta^{(m)}.
\]
Then the corresponding terms coincide
\[
h_m
\coloneqq
c(m)\chi_{E(\gamma^{(m)})}
-
c(m-1)\chi_{E(\gamma^{(m-1)})}
=
c(m)\chi_{E(\eta^{(m)})}
-
c(m-1)\chi_{E(\eta^{(m-1)})}.
\]
By the fact that \(\Phi: L_2(\cL \Fr) \to L_2(I,dm_r; L_2(\partial \Fr))\) is an isomorphism between Hilbert space and \(\Phi\big(\lambda(\cD_m(\delta_\gamma))\big)=Q_{m,n}\otimes h_m\), we obtain
\[
\int_I Q_{m,n}(s) \, \overline{Q}_{m,n'}(s)\,dm_r(s)
	\cdot
	\|h_m\|_{L_2(\partial\Fr)}^2  =
	\tau\big(
	\lambda(\cD_m(\delta_\gamma ))^*
	\lambda(\cD_m(\delta_\eta))\big)
    =
    \langle \cD_m(\delta_\gamma ), \cD_m(\delta_\eta)\rangle_{\ell_2(\Fr)}.
\]
On the left hand side,
\[
    \|h_m\|_{L_2(\partial\Fr)}^2
    =
    c(m)-c(m-1).
\]
On the right side, the inner product is zero if
\(n\ne n'\) due to the orthogonality of the subspaces \(V_{m,n}\). And if \(n=n'\), then \(\cD_m(\delta_\gamma)=\cD_m(\delta_\eta)\), and  a direct computation gives
\[
    \|\lambda(\mathcal{D}_m(\delta_\gamma))\|^2_{L_2(\cL \Fr)}=\left\|
    \frac{c(m)}{c(n)}\chi_{S(m,\gamma)}
    -
    \frac{c(m-1)}{c(n)}\chi_{S(m-1,\gamma)}
    \right\|_{\ell_2(\Fr)}^2
=
\frac{c(m)-c(m-1)}{c(n)}.
\]
Consequently,
\[
\int_I Q_{m,n}(s)\overline{Q}_{m,n'}(s)\,dm_r(s)
=
\frac{1}{c(n)}\delta_{n,n'}.
\]

Applying \(\Phi\) to both sides of \eqref{relation_Dmmu1}, we obtain the recurrence relation stated in (iii).
This completes the proof.
\end{proof}
Combining the definition of the above orthogonal system \(\{Q_{m,n}\}_{n\ge m}\) with Lemma \ref{prop1}, we identify its connection with the polynomials \(\{P_n\}_{n\ge 0}\).
\begin{lemma}\label{QnmPnm-lemma} Let \(m\ge 1,\ n\ge m\) and \(s\in I\). We have
\[Q_{m,n}(s)=P_n(s)-\sigma(s)^{n-(m-1)}P_{m-1}(s).\]
\end{lemma}
\begin{proof} By \eqref{alpha-jn} and \eqref{defQmn}, we obtain
\begin{align*}
    \label{Qnk-alpha}
		Q_{m,n}(s)&=\sum_{j=m}^n\frac{\alpha_{j,n}(s)}{c(j)}=\sum_{j=0}^n\frac{\alpha_{j,n}(s)}{c(j)}-\sum_{j=0}^{m-1}\frac{\alpha_{j,n}(s)}{c(j)}\\
		&=\sum_{j=0}^n\frac{\alpha_{j,n}(s)}{c(j)}-\sum_{j=0}^{m-1}\frac{\sigma(s)^{n-(m-1)}\alpha_{j,m-1}(s)}{c(j)}\\
		&=P_n(s)-\sigma(s)^{n-(m-1)}P_{m-1}(s). \qedhere
\end{align*}
\end{proof}
We now proceed to prove Proposition~\ref{Prop_Rbounded_kernels}.
%

\begin{proof}[Proof of Proposition \ref{Prop_Rbounded_kernels}]
	Due to Proposition~\ref{Phi_QmnHm} and Lemma \ref{l2-Vmn} , we have
	\[\Phi\big(\lambda(V_{m,n})\big)=Q_{m,n}\otimes \cH_m,\]
	and
	\[
	\widetilde{T}_a(\Phi(\lambda(f))) = \mathring{a}(n) \, \Phi(\lambda(f)),\quad f\in V_{m,n}.
	\]
	Combining these identities gives
	\begin{equation}\label{T_Qmn}
		\begin{aligned}
			\widetilde{T}_a(Q_{m,n}\otimes v)=\mathring{a}(n) \, Q_{m,n}\otimes v,\quad v\in \cH_m.
		\end{aligned}
	\end{equation}
	 Since \((Q_{m,n})_{n\geq m}\) is a complete orthogonal system in
	 \(L_2(I,dm_r)\), there exists a unique bounded operator
	 \(R_{a}^{(m)}\) on \(L_2(I,dm_r)\) satisfying
	 \[
	   R_{a}^{(m)}(Q_{m,n})
	   =
	   \mathring{a}(n) \, Q_{m,n},
	   \qquad n\geq m.
	 \]
	 Indeed,
	 \(
	   \|R_{a}^{(m)} L_2(I) \to L_2(I) \|
	   \leq \|\mathring{a}\|_{\ell_\infty}.
	 \)
	 It follows from \eqref{T_Qmn} that
	 \begin{equation}
        \label{T_R}
	 	\widetilde T_a
	 	{\big|}_{L_2(I,dm_r;\cH_m)}
	 	=
	 	R_{a}^{(m)}\otimes\operatorname{id}_{\cH_m}.
	 \end{equation}
	Using the Proposition~\ref{Phi}, we have
	\[
	\big\| T_a: L_p(\cLFr; \Erad^c) \to L_p(\cLFr; \Erad^c) \big\|
	= \big\| \widetilde{T}_a: L_p(I,dm_r; L_2(\partial \Fr)) \to L_p(I,dm_r; L_2(\partial \Fr)) \big\|.
	\]
	 Choose an orthonormal basis
	 \((e_{m,i})_{1\leq i\leq n_m}\) of each \(\cH_m\), where
	 \(n_m\coloneqq\dim\cH_m\). For any \(\xi\in L_p(I,dm_r; L_2(\partial \Fr))\), we can write
	 \[
	 \xi(s)=\sum_{m=0}^\infty\sum_{i=1}^{n_m} \xi_{m,i}(s)e_{m,i},\quad s\in I.
	 \] Through \eqref{T_R}, we have
	\begin{equation*}
		\begin{aligned}
			\|\widetilde{T}_a(\xi)\|_{L_p(I,dm_r; L_2(\partial \Fr))}^p    
            & = \bigg \|\widetilde{T}_a\big(\sum_{m=0}^\infty\sum_{i=1}^{n_m} \xi_{m,i}e_{m,i}\big) \bigg\|_{L_p(I,dm_r; L_2(\partial \Fr))}^p\\
			& = \bigg\|\sum_{m=0}^\infty\sum_{i=1}^{n_m}R_{a}^{(m)}(\xi_{m,i})e_{m,i}\bigg\|_{L_p(I,dm_r; L_2(\partial \Fr))}^p\\
			& = \int_I \bigg\|\sum_{m=0}^\infty\sum_{i=1}^{n_m}R_{a}^{(m)}(\xi_{m,i})(s)e_{m,i}\bigg\|_{L_2(\partial \Fr)}^pdm_r(s)\\
			& = \bigg\|\Big(\sum_{m=0}^\infty\sum_{i=1}^{n_m}|R_{a}^{(m)}(\xi_{m,i})|^2\Big)^{\frac{1}{2}}\bigg\|_{L_p(I,dm_r)}^p.
 		\end{aligned}
	\end{equation*}
	Similarly, we have
	\[
	   \|\xi\|_{L_p(I,dm_r; L_2(\partial \Fr))}^p 
       =
       \bigg\|\Big(\sum_{m=0}^\infty\sum_{i=1}^{n_m}|\xi_{m,i}|^2\Big)^{\frac{1}{2}}\bigg\|_{L_p(I,dm_r)}^p.
	\]
	Through these two identities and the definition of R-boundedness, we have
	\[
	\big\| T_a: L_p(\cLFr; \Erad^c) \to L_p(\cLFr; \Erad^c) \big\| = \cR_p\big\{ R_{a}^{(m)}  : m \geq 0 \big\},
	\]
	which proves (i). Finally, when \(\mathring{a}\) is finitely supported, the orthogonality of \(\{Q_{m,n}\}_{n\ge m}\) and definition of \(R_{a}^{(m)}\) gives
	\[
	   R_{a}^{(m)}(h)(s)
	   =
	   \sum_{n\geq m}
	   \mathring{a}(n)
	   \frac{\langle h,Q_{m,n}\rangle}
	   {\|Q_{m,n}\|_2^2}
	   Q_{m,n}(s) = \int_I K_{a}^{(m)}(s,t) h(t) \, dm_r(t),\quad h\in L_p(I,dm_r).
	\]
	For general bounded \(\mathring{a}\), the same computation is understood in the sense of distributions. This completes the proof.
\end{proof}

Observe as well that, in the proof above, we have intentionally given a definition of the $R$-boundedness constant that allows for repetitions on the operators present in the collection $\big\{ R_a^{(m)} : m \leq 0 \big\}$ in order to account for those operators appearing after tensoring with the identity on $\cH_m$. However, that difference is superfluous up to constants by \cite[Lemma 4.3]{KVW2016}. 

\begin{corollary}
    \label{tmn-vec-kernel} Let \(r \in \{2, 3, \dots\}\). For partial sum operators \(S_N\) on $\Fr$, we obtain
    \begin{equation}
        \label{phi-tmn-ker}
        \Phi\big(S_N(x)\big)(s)=
        \sum_{m=0}^N\int_IK_{N,m}(s,t)\Delta_m(\Phi(x)(t))dm_r(t), \quad x\in \cL \Fr,\  s\in I,
    \end{equation}
    where
    \begin{equation}
        \label{KmN}
        K_{N,m}(s,t)= \sum_{n=m}^N c(n)Q_{m,n}(s)\overline{Q}_{m,n}(t),\quad m\geq 0,\ s,t\in I.
    \end{equation}
\end{corollary} 

\subsection{The case of partial sum operators}\label{subsec: generalization of Pollard's arguments}
We now consider the partial sum operators. The resulting kernel is decomposed into weighted Hilbert-transform kernels and a regular term. Furthermore, we establish asymptotic estimates for the family $\{Q_{m,n}\}_{n\ge m}$, which will be crucial for the subsequent proof.
In order to more conveniently study these estimates, we normalize the orthogonal functions rescaling the spectral variable to $[-1,1]$. Specifically, for \(s\in [-1,1]\), we define
\[
U_{m,n}(s)\coloneqq\sqrt{c(n)} Q_{m,n}\Big(\frac{\sqrt{2r-1}}{r}s\Big).
\]
Then by Proposition~\ref{Phi_QmnHm} (iii), for each \(m\ge 1\), the sequences \(\{U_{m,n}(s)\}_{n\ge m}\) satisfy the recurrence relation
\begin{equation}\label{Qmn-rec}
	U_{m,n+1}(s)=2sU_{m,n}(s)-U_{m,n-1}(s),\quad n\ge m+1.
\end{equation}
For \(m=0\),
\begin{equation}\label{Pn-rec}
	\begin{aligned}
		U_{0,0}(s)=1,\ U_{0,1}(s)=2\sqrt{\frac{2r-1}{2r}}s,\ U_{0,2}(s)=2sU_{0,1}(s)-\sqrt{\frac{2r}{2r-1}}U_{0,0}(s),\\
		U_{0,n+1}(s)=2sU_{0,n}(s)-U_{0,n-1}(s),\quad n\ge 2.
	\end{aligned}
\end{equation}
In addition, by Lemma \ref{QnmPnm-lemma}, they satisfy
\begin{equation}
    \label{tildePwithtildeQ}
	\begin{aligned}
		U_{m,n}(s) &=U_{0,n}(s)-(s-\imagi \sqrt{1-s^2})^{n-(m-1)}U_{0,m-1}(s),\quad \ m\ge 2,\ n\ge m;\\
		U_{1,n} & =U_{0,n}(s)-(s-\imagi \sqrt{1-s^2} \, \mathrm i)^{n}\sqrt{\frac{2r}{2r-1}}U_{0,0}(s),\quad \ n\ge 1.
	\end{aligned}
\end{equation}
The corresponding probability measure becomes: \[\frac{1}{\pi}\frac{r(2r-1)\sqrt{1-s^2}}{r^2-(2r-1)s^2}ds=:w(s)ds.\] 
Now, we let
\begin{equation}\label{KmN-new}
	\widetilde K_{N,m}(s,t)=\sum_{n=m}^N{U}_{m,n}(s)\overline{U}_{m,n}(t),\quad m\geq0,s,t\in[-1,1].
\end{equation}

For the normalized system of functions $\{U_{m,n}\}_{m\geq0,n\geq m}$, we have an alternative representation as follows:
\begin{lemma}\label{Qmn-sin-lemma} For \(r\ge 2\) and \(s\in[-1,1]\),  let \(s=\cos{\theta}\) with \(\theta \in [-\pi,0]\). Then, we have
	\[{U}_{0,0}=1,\quad {U}_{0,n}(s)=\sqrt{\frac{2r}{2r-1}}\frac{\frac{{2r-1}}{r}\cos\theta\sin n\theta-\sin(n-1)\theta}{\sin\theta},\quad n\ge 1,\]
	and
		\begin{equation}\label{qnm-sin}
		{U}_{m,n}(s)=	\sqrt{\frac{2r}{2r-1}} \Big(\frac{{2r-1}}{r} \cos{\theta}-e^{\mathrm i \theta }\Big)e^{-\mathrm{i}(m-1)\theta}\frac{\sin\big((n-(m-1))\theta\big)}{\sin \theta},\quad m\geq1,\ n\ge m.
	\end{equation} 
\end{lemma}
\begin{proof}  Let \(s=\cos{\theta}\) with \(\theta \in [-\pi,0]\). Then by \eqref{sigma-def}, we have
	\begin{equation}\label{sigma-n-theta}
		\sigma \Big(\frac{\sqrt{2r-1}}{r}s \Big)
            =\frac{1}{\sqrt{2r-1}}e^{\mathrm i \theta}.
	\end{equation}
	For \(m=0\) and \(n\geq 1\), substituting these formulae to \eqref{Pn-sigma}, we have
	\begin{align*}
		U_{0,n}(s)&=\sqrt{c(n)}{P}_{n}\Big(\frac{\sqrt{2r-1}}{r}s\Big)\\
        &=\sqrt{\frac{2r}{2r-1}} \Big[\frac{\frac{{2r-1}}{r} \cos{\theta}-e^{\mathrm i \theta }}{e^{-\mathrm{i}\theta}-e^{\mathrm{i}\theta}}e^{-\mathrm{i}n\theta}+\frac{e^{-\mathrm i \theta }-\frac{{2r-1}}{r} \cos{\theta}}{e^{-\mathrm{i}\theta}-e^{\mathrm{i}\theta}}e^{\mathrm{i}n\theta}\Big]\\
		&=\sqrt{\frac{2r}{2r-1}}\frac{\frac{{2r-1}}{r}\cos\theta(e^{-\mathrm in\theta}-e^{\mathrm in\theta})-(e^{\mathrm i\theta-\mathrm in\theta}-e^{-\mathrm i\theta+\mathrm in\theta})}{-2\mathrm i\sin\theta}\\
			&=\sqrt{\frac{2r}{2r-1}}\frac{\frac{{2r-1}}{r}\cos\theta\sin n\theta-\sin(n-1)\theta}{\sin\theta}.
	\end{align*}
	For \(m\geq2\) and \(n\geq m\), substituting these formulae and using Lemma \ref{QnmPnm-lemma} and \eqref{Pn-sigma}, we have
	\begin{equation*}
		\begin{aligned}
			{U}_{m,n}(s)&={U}_{0,n}(s)-(s-\sqrt{1-s^2}i)^{n-m+1}{U}_{0,m-1}(s)\\
			&=\sqrt{\frac{2r}{2r-1}} \bigg[\frac{\frac{{2r-1}}{r} \cos{\theta}-e^{\mathrm i \theta }}{e^{-\mathrm{i}\theta}-e^{\mathrm{i}\theta}}e^{-\mathrm{i}n\theta}+\frac{e^{-\mathrm i \theta }-\frac{{2r-1}}{r} \cos{\theta}}{e^{-\mathrm{i}\theta}-e^{\mathrm{i}\theta}}e^{\mathrm{i}n\theta}\\
			&\quad\quad\quad\quad\quad \quad-e^{\mathrm i(n-(m-1))\theta}\Big(\frac{\frac{{2r-1}}{r} \cos{\theta}-e^{\mathrm i \theta }}{e^{-\mathrm{i}\theta}-e^{\mathrm{i}\theta}}e^{-\mathrm{i}(m-1)\theta}+\frac{e^{-\mathrm i \theta }-\frac{{2r-1}}{r} \cos{\theta}}{e^{-\mathrm{i}\theta}-e^{\mathrm{i}\theta}}e^{\mathrm{i}(m-1)\theta}\Big)\bigg],\\
			&=\sqrt{\frac{2r}{2r-1}} \frac{\frac{{2r-1}}{r} \cos{\theta}-e^{\mathrm i \theta }}{-2\mathrm i\sin \theta}\Big[e^{-\mathrm{i}n\theta}-e^{\mathrm{i}(n-2(m-1))\theta}\Big]\\
			&=\sqrt{\frac{2r}{2r-1}} \frac{\frac{{2r-1}}{r} \cos{\theta}-e^{\mathrm i \theta }}{-2\mathrm i\sin \theta}e^{-\mathrm{i}(m-1)\theta}\Big[e^{-\mathrm{i}(n-(m-1))\theta}-e^{\mathrm{i}(n-(m-1))\theta}\Big].
		\end{aligned}
	\end{equation*}
	The same conclusion holds in the case \(m=1\). The proof is complete.
\end{proof}

We introduce an auxiliary family \(\{q_{m,n}\}_{n\ge m}\) defined on $[-1,1]$ by
		\[q_{0,0}(s)=2\big(\frac{2r-1}{2r}\big)^{\frac14}\frac{1-s^2}{\sqrt{U_{0,2}(1)}}\]
        and
        \[
		q_{m,n}(s)=-\frac{1}{2}\sqrt{\frac{D_n}{-D_{n+1}}}\frac{D_n(s)}{D_n}
		,\quad m\geq0,n\geq \max\{1,m\},\]
		where
		\[
		\ D_n= \begin{vmatrix}
			{U}_{m,n+1}(1) & {U}_{m,n}(1) \\
			{U}_{m,n+1}(-1)&  {U}_{m,n}(-1)
		\end{vmatrix}\ \quad\mbox{and}\quad D_n(s)= \begin{vmatrix}
			U_{m,n+2}(s) & {U}_{m,n+1}(s) & {U}_{m,n}(s)\\
			U_{m,n+2}(1) & {U}_{m,n+1}(1) & {U}_{m,n}(1)\\
			U_{m,n+2}(-1)&  {U}_{m,n+1}(-1)&{U}_{m,n}(-1)
		\end{vmatrix}.
		\]
The family \(\{q_{m,n}\}_{n\ge m}\)  plays an essential role in the kernel decomposition of partial sum operators. And we collect estimates for the  family \(\{U_{m,n}\}_{n\ge m}\)  and \(\{q_{m,n}\}_{n\ge m}\).
\begin{lemma}\label{Qmn-qmn-w}
	Let \(r\ge 2\) and \(m\ge 0\). There exists a universal constant \(A>0\) such that the following assertions hold.
	\begin{enumerate}[label=(\roman*), font=\normalfont]
		\item  For all \(n\ge m\) and \(s\in [-1,1]\),
		\[
		(1-s^2)^{1/4}w(s)^{1/2}|{U}_{m,n}(s)|\le A.
		\]
		\item 
		The family
        \[
        \left\{(1-s^2)^{-1}q_{m,n}(s)\right\}_{n\ge m}
        \]
        is orthonormal with respect to the measure \((1-s^2)w(s)\,ds\). And for all \(n\ge m\) and \(s\in[-1,1]\), we have
        \[
        (1-s^2)^{-1/4}w(s)^{1/2}
        \left|q_{m,n}(s)\right|
        \le A.
        \]
	\end{enumerate}
\end{lemma}
\begin{proof}
	 Let \(s=\cos{\theta}\) with \(\theta \in [-\pi,0]\). We have
	 \begin{equation}\label{w-theta}
	 w(s)=\frac{-1}{\pi} \frac{r(2r-1)\sin\theta}{r^2-(2r-1)\cos^2\theta}.
	 \end{equation}
	 For \(m=0\),  by Lemma \ref{Qmn-sin-lemma}, we have
	 \begin{align*}
	 	(1-s^2)^{1/4}w(s)^{1/2}|{U}_{0,n}(s)|
	 	&=(-\sin\theta)^{1/2}\big(\frac{-1}{\pi} \frac{r(2r-1)\sin\theta}{r^2-(2r-1)\cos^2\theta}\big)^{1/2}\\
	 	&\quad\cdot\Big|\sqrt{\frac{2r}{2r-1}}\frac{\frac{{2r-1}}{r}\cos\theta\sin n\theta-\sin(n-1)\theta}{\sin\theta}\Big|\leq 6.
	 \end{align*}	 
Similarly, 	using Lemma \ref{Qmn-sin-lemma} with \(m\ge 1\), we get :
	\begin{equation*}
		\begin{aligned}
			(1-s^2)^{1/4}w(s)^{1/2}|{U}_{m,n}(s)|&=(-\sin\theta)^{1/2}\big(\frac{-1}{\pi}\frac{r(2r-1)\sin\theta}{r^2-(2r-1)\cos^2\theta}\big)^{1/2}\\
			&\quad\cdot\Big|\sqrt{\frac{2r}{2r-1}} \big(\frac{{2r-1}}{r} \cos{\theta}-e^{\mathrm i \theta }\big)e^{-\mathrm{i}(m-1)\theta}\frac{\sin\big((n-(m-1))\theta\big)}{\sin \theta}\Big|\le 6.
		\end{aligned}
	\end{equation*}
	
	We now turn to prove (ii). For $m=0$, condition~(ii) reduces to Hypothesis~4 in \cite{P1948},
which was verified in Section~6 (Reduction of Hypotheses) of that paper.
We note that, although Hypothesis~4 as stated in \cite{P1948} does not
explicitly require the implied constant to be independent of $r$,
the arguments given in Section~6 in fact yield a constant that is
universal, hence independent of $r$. For \(m\geq1\), by Lemma \ref{Qmn-sin-lemma}, we have
	\[
	{U}_{m,n}(1)=\sqrt{\frac{2r}{2r-1}}\frac{r-1}{r}(n-m+1),\quad n\ge m,
	\]
	and
	\[
	{U}_{m,n}(-1)=(-1)^n\sqrt{\frac{2r}{2r-1}}\frac{r-1}{r}(n-m+1),\quad n\ge m.
	\]
	By the definition  of \(D_n\) and $D_{n}(s)$ in (ii), we have
	\[D_n=2(-1)^n\frac{2r}{2r-1}\Big(\frac{r-1}{r}\Big)^2(n-m+2)(n-m+1),\]
	\[D_{n+1}=2(-1)^{n+1}\frac{2r}{2r-1}\Big(\frac{r-1}{r}\Big)^2(n-m+3)(n-m+2),\]
	and
	\[D_n(s)=D_{n+1}U_{m,n}(s)+D_nU_{m,n+2}(s),\quad s\in[-1,1].\]
	Then, we can express  $q_{m,n}$ as
	\begin{equation}\label{AnCn}
		\begin{aligned}
			q_{m,n}(s)&=-\frac{1}{2}
			\sqrt{\frac{D_n}{-D_{n+1}}} \frac{D_n(s)}{D_n}\\
			&=\frac12\sqrt{\frac{D_{n+1}}{-D_n}}
			U_{m,n}(s)-
			\frac12\sqrt{\frac{D_n}{-D_{n+1}}}U_{m,n+2}(s)\\
			&=\frac{1}{2}\sqrt{\frac{n-m+3}{n-m+1}}U_{m,n}(s)-\frac12\sqrt{\frac{n-m+1}{n-m+3}}U_{m,n+2}(s).
		\end{aligned}
	\end{equation}
	Substituting \eqref{qnm-sin} into the above formula, we obtain
	\begin{equation}\label{qnm-theta}
		\begin{aligned}
			q_{m,n}(s)=&\frac{1}{2}\sqrt{\frac{2r}{2r-1}} e^{-i(m-1)\theta}\Big(e^{i\theta}-\frac{2r-1}{r}\cos{\theta}\Big)\\
			&\cdot\Big( \sqrt{\frac{n-m+1}{n-m+3}}\frac{\sin\big((n-m+3)\theta\big)}{\sin{\theta}} - \sqrt{\frac{n-m+3}{n-m+1}}\frac{\sin\big({(n-m+1)\theta}\big)}{\sin{\theta}} \Big).
		\end{aligned}
	\end{equation}
	Now, we check the orthonormality of \(\{(1-s^2)^{-1} q_{m,n}(s)\}_{n\ge m}\) with respect to the measure \((1-s^2) w(s) ds\), by \eqref{w-theta}, we have
	\begin{equation*}
		\begin{aligned}
			&\quad\int_{-1}^1(1-s^2)^{-1}q_{m,n}(s)(1-s^2)^{-1}\overline{q}_{m,\ell}(s)(1-s^2)w(s)ds\\
			&=\frac{1}{4}\frac{2r}{2r-1}\int_{-\pi}^0\Big|e^{i\theta}-\frac{2r-1}{r}\cos{\theta}\Big|^2 \Big( \sqrt{\frac{n-m+1}{n-m+3}}\frac{\sin{(n-m+3)\theta}}{\sin{\theta}} - \sqrt{\frac{n-m+3}{n-m+1}}\frac{\sin{(n-m+1)\theta}}{\sin{\theta}} \Big)\\
			&\quad\quad\quad\cdot\Big( \sqrt{\frac{\ell-m+1}{\ell-m+3}}\frac{\sin{(\ell-m+3)\theta}}{\sin{\theta}} - \sqrt{\frac{\ell-m+3}{\ell-m+1}}\frac{\sin{(\ell-m+1)\theta}}{\sin{\theta}} \Big)\cdot\frac{1}{\sin^2\theta}\frac{1}{\pi} \frac{r(2r-1)\sin^2\theta}{r^2-(2r-1)\cos^2\theta}d\theta\\
			&=\frac{1}{2\pi}\int_{-\pi}^0 \Big( \sqrt{\frac{n-m+1}{n-m+3}}\frac{\sin{(n-m+3)\theta}}{\sin{\theta}} - \sqrt{\frac{n-m+3}{n-m+1}}\frac{\sin{(n-m+1)\theta}}{\sin{\theta}} \Big)\\
			&\quad\quad\quad\quad\cdot\Big( \sqrt{\frac{\ell-m+1}{\ell-m+3}}\frac{\sin{(\ell-m+3)\theta}}{\sin{\theta}} - \sqrt{\frac{\ell-m+3}{\ell-m+1}}\frac{\sin{(\ell-m+1)\theta}}{\sin{\theta}} \Big)d\theta.
		\end{aligned}
	\end{equation*}
	Notice that for  the trigonometric system\(\{\frac{\sin(\ell\theta)}{\sin\theta}\}_{\ell\ge 1}\), we have 
	\begin{equation*}\int_{-\pi}^{0}
		\frac{\sin(\ell\theta)}{\sin\theta}\,
		\frac{\sin(\ell'\theta)}{\sin\theta}\,d\theta
		=\left\{
		\begin{aligned}
			&\pi\min\{\ell,\ell'\}, & \hbox{\(\ell\equiv \ell'\pmod 2\)},\\[4pt]
			&0, & \hbox{\(\ell\not \equiv \ell' \pmod 2\)}.
		\end{aligned}
		\right.
	\end{equation*}
	Hence, we obtain that \(\{(1-s^2)^{-1}q_{m,n}(s)\}_{n\geq m}\) is orthonormal. 
	Using \eqref{qnm-theta}, we have
	\begin{equation*}
		\begin{aligned}
			|q_{m,n}(s)|=\frac{1}{2}&\Big| \sqrt{\frac{2r}{2r-1}} e^{-i(m-1)\theta}\Big(e^{i\theta}-\frac{2r-1}{r}\cos{\theta}\Big) \Big|\\
			\cdot &\Big| \sqrt{\frac{n-m+1}{n-m+3}}\frac{\sin{(n-m+3)\theta}}{\sin{\theta}} - \sqrt{\frac{n-m+3}{n-m+1}}\frac{\sin{(n-m+1)\theta}}{\sin{\theta}} \Big|\\
			\le 3&\Big| \sqrt{\frac{n-m+1}{n-m+3}}\frac{\sin{(n-m+3)\theta}}{\sin{\theta}} - \sqrt{\frac{n-m+3}{n-m+1}}\frac{\sin{(n-m+1)\theta}}{\sin{\theta}} \Big|\\
			=3&\Big| \frac{(n-m+1)\sin{(n-m+3)\theta}-(n-m+3)\sin{(n-m+1)\theta}}{\sqrt{(n-m+1)(n-m+3)}sin\theta} \Big|.
		\end{aligned}
	\end{equation*}
	Using the identities
	\begin{equation*}
		\begin{aligned}
			&\sin{(n-m+3)\theta}=\sin{(n-m+2)\theta}\cos{\theta}+ \cos{(n-m+2)\theta}\sin{\theta},\\
			&\sin{(n-m+1)\theta}=\sin{(n-m+2)\theta}\cos{\theta}- \cos{(n-m+2)\theta}\sin{\theta},
		\end{aligned}
	\end{equation*}
	we derive
	\begin{equation*}
		\begin{aligned}
			\big| q_{m,n}(s)\big| &\le 3\Big| \frac{-2\sin{(n-m+2)\theta}\cos{\theta}+(2n-2m+4)\cos{(n-m+2)\theta}\sin{\theta}}{\sqrt{(n-m+1)(n-m+3)}sin\theta}  \Big|\\
			&\le 3 \Big( \Big|\frac{2\cos{\theta}\sin{(n-m+2)\theta}}{\sqrt{(n-m+1)(n-m+3)}sin\theta}\Big|+ \Big| \frac{2n-2m+4}{\sqrt{(n-m+1)(n-m+3)}} \Big| \Big)\\
			&\le 6\Big|  \frac{2n-2m+4}{\sqrt{(n-m+1)(n-m+3)}} \Big|\le 12,
		\end{aligned}
	\end{equation*}
	where the last second inequality uses the bound \(|\sin(\ell\theta)| \le \ell|\sin\theta|\) for \(\ell \geq 1\). Finally, we get
	\begin{align*}
		(1-s^2)^{-1/4}w(s)^{1/2}| q_{m,n}(s)| &= (\sin^2{\theta})^{-1/4}\Big( \frac{1}{\pi}\frac{r(2r-1)| \sin{\theta} |}{r^2-(2r-1)\cos^2{\theta}}\Big)^{1/2} | q_{m,n}(s)|\\
		&\le 2| q_{m,n}(s)|\le 24. \qedhere
	\end{align*}
\end{proof}
The following lemma provides a fundamental identity for the partial sum kernel employed in the subsequent decomposition.
\begin{lemma}\label{KmN-new-expression} For each \(m\geq0, N\geq1\), we have
	\begin{equation*}
		\widetilde K_{N,m}(s,t)= \frac{U_{m,N+1}(s)\overline{U}_{m,N}(t)-U_{m,N}(s)\overline{U}_{m,N+1}(t)}{2(s-t)},\quad s,t\in [-1,1].
	\end{equation*}
\end{lemma}
\begin{proof} For the case $m=0$,  we can refer to \cite[Page 357]{P1948}.  For each fixed \(m\geq1\),  let \(U_{m,m-1}(s)=0\). Then   \eqref{tildePwithtildeQ} and \eqref{Pn-rec} yield
	\[U_{m,m-1}(s)=2sU_{m,m}(s)-U_{m,m+1}(s).\]
	 For \(n\geq m\), it can be easily derived from the above equation and  the recurrence formula \eqref{Qmn-rec}
	\begin{equation*}
		\begin{aligned}
			&{U}_{m,n+1}(s)\overline{U}_{m,n}(t)-{U}_{m,n}(s)\overline{U}_{m,n+1}(t)\\
            &=\Big(2s{U}_{m,n}(s)-{U}_{m,n-1}(s)\Big)\overline{U}_{m,n}(t)-{U}_{m,n}(s)\Big(2t\overline{U}_{m,n}(t)-\overline{U}_{m,n-1}(t)\Big)\\
			&=(2s-2t){U}_{m,n}(s)\overline{U}_{m,n}(t)+\Big[{U}_{m,n}(s)\overline{U}_{m,n-1}(t)-{U}_{m,n-1}(s)\overline{U}_{m,n}(t)\Big].
		\end{aligned}
	\end{equation*}
	The above equation can be reformulated as
    \begin{align*}
        {U}_{m,n}(s)\overline{U}_{m,n}(t)&=\frac{{U}_{m,n+1}(s)\overline{U}_{m,n}(t)
        -{U}_{m,n}(s)\overline{U}_{m,n+1}(t)}{2(s-t)}\\
    &\quad\quad-\frac{{U}_{m,n}(s)\overline{U}_{m,n-1}(t)-{U}_{m,n-1}(s)\overline{U}_{m,n}(t)}{2(s-t)}.
    \end{align*}
	Consequently, for $m\geq1$
	\begin{equation*}
		\begin{aligned}
			\widetilde K_{N,m}(s,t)&=\sum_{n=m}^N\Big(\frac{{U}_{m,n+1}(s)\overline{U}_{m,n}(t)-{U}_{m,n}(s)\overline{U}_{m,n+1}(t)}{2(s-t)}\\
            &\quad\quad\quad\quad-\frac{{U}_{m,n}(s)\overline{U}_{m,n-1}(t)-{U}_{m,n-1}(s)\overline{U}_{m,n}(t)}{2(s-t)}\Big)\\
			&=\frac{U_{m,N+1}(s)\overline{U}_{m,N}(t)-U_{m,N}(s)\overline{U}_{m,N+1}(t)}{2(s-t)}\\
            &\quad-\frac{U_{m,m}(s)\overline{U}_{m,m-1}(t)-U_{m,m-1}(s)\overline{U}_{m,m}(t)}{2(s-t)},
		\end{aligned}
	\end{equation*}
where the initial term vanishes. 
	The proof is complete.
\end{proof}

Using the preceding lemma, we derive the following key expression for the kernel of the partial sum operator. Specifically, the kernel decomposes into three terms, consisting of two weighted Hilbert transform kernels and a regular component.
\begin{proposition}\label{AN-CN} Let \(\ N\geq 1,1\leq m\leq N\). For \(s,t \in [-1,1]\), we have
	\begin{equation}\label{KmN-3part}
		\begin{aligned}
			\widetilde K_{N,m}(s,t)&= \left(-\big( 1 - \frac{A_{m,N}}{C_{m,N}} \big)^{-1} \frac{1}{2C_{m,N}}\right)\frac{ U_{m,N+1}(s)\overline{q}_{m,N}(t)}{s - t}\\
			&\quad+\left(\big( 1 - \frac{A_{m,N}}{C_{m,N}} \big)^{-1} \frac{1}{2C_{m,N}}\right) \frac{q_{m,N}(s)\overline{U}_{m,N+1}(t)}{s-t}\\
			&\quad+\left(- \big( 1 -  \frac{A_{m,N}}{C_{m,N}} \big)^{-1}\right)  U_{m,N+1}(s) \overline{U}_{m,N+1}(t)\\
			&=:\sum_{i=1}^3a_{m,N}^i K_{N,m}^i(s,t)
		\end{aligned}
	\end{equation}
	with $|a_{m,N}^i|\leq C$, where \(A_{m,N}=\frac{1}{2}\sqrt{\frac{N-m+3}{N-m+1}}\), \(C_{m,N}=-\frac12\sqrt{\frac{N-m+1}{N-m+3}}\)
and \(C\)	is a constant independent  of \(m\) and \(N\).
\end{proposition}
\begin{proof}
	By \eqref{AnCn}, we can express  $q_{m,N}$ as
	\begin{equation*}
		q_{m,N}(s)=A_{m,N} U_{m,N}(s)+C_{m,N} U_{m,N+2}(s).
	\end{equation*}
	Hence,
	\[U_{m,N+2}(s)=\frac{1}{C_{m,N}}
	q_{m,N}(s)-\frac{A_{m,N}}{C_{m,N}} U_{m,N}(s).\]
	Based on  \eqref{KmN-new}, Lemma \ref{KmN-new-expression} and the above equation, a direct computation yields that
	\begin{equation*}
		\begin{aligned}
			\widetilde K_{N,m}(s,t)&=	\widetilde  K_{N+1,m}(s,t)-U_{m,N+1}(s) \overline{U}_{m,N+1}(t)\\
			&=\frac{1}{2(s-t)}\Big[U_{m,N+2}(s) \overline{U}_{m,N+1}(t)	- U_{m,N+1}(s)\overline{U}_{m,N+2}(t)\Big]-U_{m,N+1}(s)\overline{U}_{m,N+1}(t) \\
			&=\frac{1}{2(s-t)}	\bigg[
			\Big(\frac{1}{C_{m,N}}
			q_{m,N}(s)-\frac{A_{m,N}}{C_{m,N}} U_{m,N}(s)\Big)
			\overline{U}_{m,N+1}(t)\\
            & \qquad\qquad\qquad -U_{m,N+1}(s)\Big(\frac{1}{C_{m,N}}
			\bar{q}_{m,N}(t)-\frac{A_{m,N}}{C_{m,N}} \overline{U}_{m,N}(t)\Big)\bigg]-U_{m,N+1}(s)\overline{U}_{m,N+1}(t).
		\end{aligned}
	\end{equation*}
	Then, by Lemma \ref{KmN-new-expression}, we have
	\begin{equation*}
		\begin{aligned}
			\widetilde K_{N,m}(s,t)&=\frac{1}{2C_{m,N}}\frac{q_{m,N}(s)\overline{U}_{m,N+1}(t)}{s-t}-\frac{1}{2C_{m,N}}\frac{U_{m,N+1}(s)\bar{q}_{m,N}(t)}{s-t}\\
			&\qquad\qquad+\frac{A_{m,N}}{C_{m,N}}\frac{U_{m,N+1}(s)\overline{U}_{m,N}(t)-U_{m,N}(s)\overline{U}_{m,N+1}(t)}{2(s-t)}-U_{m,N+1}(s)\overline{U}_{m,N+1}(t) \\
			&=\frac{1}{2C_{m,N}}\frac{q_{m,N}(s)\overline{U}_{m,N+1}(t)}{s-t}-\frac{1}{2C_{m,N}}\frac{U_{m,N+1}(s)\bar{q}_{m,N}(t)}{s-t}\\
            &\qquad\qquad+\frac{A_{m,N}}{C_{m,N}}\widetilde K_{N,m}(s,t)-U_{m,N+1}(s)\overline{U}_{m,N+1}(t).
		\end{aligned}
	\end{equation*}
	Then,
	\begin{align*}
        \Big( 1-\frac{A_{m,N}}{C_{m,N}} \Big) \, \widetilde K_{N,m}(s,t) =-\frac{1}{2C_{m,N}}&\frac{U_{m,N+1}(s)\bar{q}_{m,N}(t)}{s-t}+\frac{1}{2C_{m,N}}\frac{q_{m,N}(s)\overline{U}_{m,N+1}(t)}{s-t}\\ & -U_{m,N+1}(s)\overline{U}_{m,N+1}(t).
	\end{align*}
	Hence, we have the desired decomposition  \eqref{KmN-3part} of the kernel.
	
	It follows easily from  the definition of \(A_{m,N}\) and \(C_{m,N}\)  that
	\[
        \Big|\frac{1}{2C_{m,N}}\Big|\leq\sqrt{3}
        \quad \text{and} \quad 
        \Big|\big(1-\frac{A_{m,N}}{C_{m,N}}\big)^{-1}\Big| \leq \frac12
    \]
	 and the proof is complete.
\end{proof}
\begin{remark}\label{KN0-dec} For $N\ge 1$ and $s,t\in[-1,1]$, the kernel
$\widetilde K_{N,0}(s,t)$ admits a decomposition analogous to that
in Proposition~3.12, with coefficients bounded by a universal
constant independent of $r$ and $N$; see \cite{P1948}.
\end{remark}

\section{The proofs of the mixed-norm estimates}\label{sec: proof of main theorem}
In this section, we present a detailed proof of Theorem \ref{Tm-pp}. By duality arguments and the \(A_p\)-weighted boundedness of the Hardy–Littlewood maximal operator \(M\), we reduce the problem to the weighted \(L_2\) inequality for partial sum operators, namely Theorem \ref{mainthm}. This conclusion is proved with the help of the auxiliary lemmas presented below. In this section, for every integrable function \(g\) on \([-1,1]\), 
we denote
\[
    \cM_\alpha(g)(s)=\sup_{h>0}\Big(\frac{1}{2h}\int_{s-h}^{s+h}|g(t)|^\alpha dt\Big)^\frac{1}{\alpha}.
\]


\begin{lemma}\label{kernel1}
	For each \(\alpha>2\) and \(1<\rho<\frac{2\alpha}{\alpha+2}\), there exists \(C_{\rho,\alpha}>0\) such that for every integrable function \(g\), we have
    \[
        \kappa(s) \,
        \cM_\rho\!\left(\frac{g}{\kappa}\right)(s)
        \leq C_{\rho,\alpha}\,\cM_\alpha(g)(s),
        \qquad 
        \text{ where } \kappa(s) = (1-s^2)^{\frac12}, \; s \in [-1,1].
    \]
\end{lemma}
\begin{proof}
	 Fix \(h>0\) and denote
	\[
	   \mathcal{I} = \kappa(s) \Big( \frac{1}{2h}\int_{s-h}^{s+h}\bigg| \frac{g(t)}{\kappa(t)} \bigg|^\rho \chi_{[-1,1]}(t)\, dt \Big)^\frac{1}{\rho}.
	\]
	Then it suffices to show 
	\(\mathcal{I}\le C \cM_\alpha(g)(s)\). 
    We shall use the elementary inequalities 
    \begin{equation}\label{ks}
       1-|s|\leq \kappa(s)^2\leq 2(1-|s|),\quad s\in[-1,1].
    \end{equation}
	
	\noindent\textbf{Case 1: \(h\le \frac{1-|s|}{2}\).} This case means that the integration interval is away from the singular points. For \(t\in [s-h,s+h]\cap[-1,1]\), we have  $|t-s|\leq h\leq \frac{1-|s|}{2}$, which implies 
    \[1-|t|\ge 1-|s|-|t-s|\geq\frac{1-|s|}{2}\]  
    Together with \eqref{ks}, gives
	\[\frac{\kappa(s)}{\kappa(t)}\leq \big(\frac{2(1-|s|)}{1-|t|}\big)^{\frac12}\le 2.
	\]
	Consequently, since \(\rho<\alpha\), we obtain
	\[
	\mathcal{I}\le 2\Big( \frac{1}{2h}\int_{s-h}^{s+h} |g(t)|^\rho dt\Big)^\frac{1}{\rho}\le 2\cM_\alpha(g)(s).
	\]
	
	\noindent\textbf{Case 2: \(h > \frac{1-|s|}{2}\).}  Set \(J=[s-h,s+h]\cap[-1,1]\) and \(\beta=\frac{\rho\alpha}{2(\alpha-\rho)}\). Our assumptions imply that \(0<\beta<1\). By H\"{o}lder's inequality with \(\frac1\rho=\frac1\alpha+\frac{1}{2\beta}\), we have 
	\begin{equation}\label{case2}
		\begin{aligned}
			\mathcal{I} 
			\le \kappa(s) \,  \Big(\frac{1}{2h} \int _J |g(t)|^\alpha dt \Big)^\frac{1}{\alpha} \Big( \frac{1}{2h}\int_J\kappa(t)^{-2\beta}dt\Big)^{\frac{1}{2\beta}}.
		\end{aligned}
	\end{equation}
To estimate the last factor, observe that, \(1-t^2\geq (1+t)\) for \(t\in[-1,0]\),  and   \(1-t^2\geq 1-t\) for \(t\in[0,1]\), since \(\beta>0\), it follows that
	\begin{align*}
	   \int_J\kappa(t)^{-2\beta}dt&\leq \int_{J\cap [-1,0]} (1+t)^{-\beta}dt+\int_{J\cap [0,1]} (1-t)^{-\beta}dt\\
       &\leq 2\int_0^{2h}u^{-\beta}du=\frac{2(2h)^{1-\beta}}{1-\beta}.
	\end{align*}
Substituting this into \eqref{case2}, we obtain
	\begin{equation*}
		\begin{aligned}
			\mathcal{I} & \lesssim_{\rho,\alpha} \kappa(s)h^{-\frac{1}{2}}\cM_\alpha(g)(s).
		\end{aligned}
	\end{equation*}
By \eqref{ks} and the assumption \(h>(1-|s|)/2\),
\[\kappa(s)^2
    \le 2(1-|s|)
    <4h.\]
Consequently,
\[
    \mathcal{I}
    \lesssim_{\rho,\alpha}\cM_\alpha(g)(s).
\]
	This completes the proof.
\end{proof}

In what follows, \(A_p\) denotes the class of Muckenhoupt \(A_p\)-weights on \([-1,1]\). For an interval \(J\subset[-1,1]\), we denote the average integral \(\frac{1}{|J|}\int_J\cdot ds\) by \(\avgint_{J} \). Let \(1<p<\infty\). We say that an  integrable function \(w\) on \([-1,1]\), with \(w>0\) a.e., belongs to \(A_p\) if
\[
[w]_{A_p}
\coloneqq
\sup_{J\subset[-1,1]}
\left(\avgint_J w(s)\right)
\left(\avgint_J w(s)^{-\frac{1}{p-1}}\right)^{p-1}
<\infty,
\]
where the supremum is taken over all intervals \(J\subset[-1,1]\). For \(p=1\), an integrable function \(w\) on \([-1,1]\), with \(w>0\) a.e., belongs to \(A_1\) if
\[
[w]_{A_1}
\coloneqq
\sup_{J\subset[-1,1]}
\frac{\avgint_J w(s)}
{\operatorname*{ess\,inf}_{s\in J} w(s)}
<\infty,
\]
where the supremum is taken over all intervals \(J\subset[-1,1]\). For further details on the theory of weights, we refer the reader to \cite{Ga2014}.
\begin{lemma}\label{kernel2}
	For each \(v\in A_1\), we have
    \[
        \big[\kappa \, v\big]_{A_2}
        \lesssim
        [v]_{A_1}, 
        \quad \text{ where } \kappa(s) = (1-s^2)^\frac12.
    \]
    In particular, for every \(\alpha>1\) and every integrable function \(g\), we have
    \[
        \big[\kappa(s) \, \cM_\alpha g\big]_{A_2}
        \le C_\alpha.
    \]
\end{lemma}
\begin{proof}
    Suppose \(v \in A_1\). Notice that \(v\in A_1\subset A_q\) with \(1<q<\frac{3}{2}\).  We take an arbitrary interval \(J\subset [-1,1]\). Then
	\begin{equation*}
		\begin{aligned}
			a& := \Big( \avgint_J(1-s^2)^{\frac{1}{2}}v \Big)\Big( \avgint_J\big((1-s^2)^\frac{1}{2}v\big)^{-1} \Big)\\
			&\lesssim \Big( \avgint_J (1-s^2)^\frac{1}{2}v\Big)\Big( \avgint_Jv^{-\frac{1}{q-1}} \Big)^{q-1}\Big( \avgint_J (1-s^2)^{-\frac{1}{2(2-q)}}\Big)^{2-q}.
		\end{aligned}
	\end{equation*}
	where the last inequality follows from H\"{o}lder's inequality with respect to the index \((\frac{1}{q-1},\frac{1}{2-q})\).  Observe that \((1-s^2)^\ell \in A_1\) whenever \(-1 < \ell \le 0\). Since \(q < \frac{3}{2}\), we have \(-\frac{1}{2(2-q)} \in (-1,0]\) and  \((1-s^2)^{-\frac{1}{2(2-q)}} \in A_1\). By the definition of  \(A_1\) weight, we obtain
	\[
	\avgint_J (1-s^2)^\frac{1}{2(q-2)} \le [(1-s^2)^{-\frac{1}{2(2-q)}}]_{A_1}\operatorname*{ess\,inf} \big( (1-s^2)^{-\frac{1}{2(2-q)}} \big).
	\]
	Substituting this estimate into the bound for \(a\) gives
	\[
	a\lesssim_q \Big( \avgint_J \Big[(1-s^2)^\frac{1}{2}\operatorname*{ess\,inf}\big((1-s^2)^{-\frac{1}{2}}\big)\Big]v\Big)\Big( \avgint_Jv^{-\frac{1}{q-1}} \Big)^{q-1}\lesssim[v]_{A_q}\leq[v]_{A_1}.
	\]
	The last second inequality uses the fact that \(v \in A_q\).
	Therefore, \((1-s^2)^{\frac{1}{2}} v(s) \in A_2\).
	
	For each  integrable function \(g\), by the Coifman–Rochberg Theorem \cite[Theorem 8.33]{KLV2021}, we have \(\cM_\alpha(g)\in A_1\)  and \([\cM_\alpha(g)]_{A_1}\lesssim \frac{\alpha}{\alpha-1}\). This completes the proof.
\end{proof}

We now proceed to prove  the weighted \(L^2\) inequality in Theorem~\ref{mainthm}.

\begin{proof}[Proof of Theorem \ref{mainthm}]
	\medskip
	\noindent\textbf{Step 1:} We first reduce the weighted \(L_2\) inequality to its vector valued version by  Plancherel theory on  free groups.
	We may assume that
	\(x=\lambda(f)\) and \(y\in (\cL \Fr^{\rad})_+\)
	for some \(f\in\mathbb C[\Fr]\). 
	 By Lemma \ref{Plancherelthm} and Lemma \ref{phi-rep-con}, we have
	 \begin{align*}
    \tau\big(|x|^2\cdot y\big) 
        & = \tau\big(\EE_\rad(|x|^2)\cdot y \big)
          =\int_I \phi_s\big(\EE_\rad(|x|^2)\cdot y\big)(s) \, dm_r(s)\\
        & = \int_I\phi_s\big(\EE_\rad(|x|^2)\big)\cdot \phi_s(y) \, dm_r(s)
          =\int_I\phi_s\big(\lambda(f)^*\lambda(f)\big)\cdot \Gamma(y)(s) \, dm_r(s).
	 \end{align*}
And Lemma \ref{Plancherelthm} gives
	\begin{equation}\label{x2-y}
	\begin{aligned}
		\tau\big(|x|^2\cdot y\big)
		&=\int_I\|\widetilde\pi_s (f)\mathds 1\|_{L_2(\partial \Fr)}^2\Gamma (y)(s)dm_r(s)\\
		&=\int_I\|\Phi (x)(s)\|_{L_2(\partial \Fr)}^2\Gamma (y)(s)dm_r(s),
	\end{aligned}
\end{equation}
	 Substituting \(S_N(x)\) for \(x\) in \eqref{x2-y},  we obtain
	\begin{equation*}
			\tau\big(|S_N(x)|^2\cdot y\big)=\int_I\|\Phi (S_N(x))(s)\|_{L_2(\partial \Fr)}^2\Gamma (y)(s)dm_r(s).
	\end{equation*}
Similarly, we have
	\begin{equation*}
			\tau\big(|x|^2\cdot \cM_\alpha(y)\big)=\int_{I}\|\Phi(x)(s)\|_{L_2(\partial \Fr)}^2\Gamma (\mathcal M_\alpha(y))(s)dm_r(s)=\int_{I}\|\Phi(x)(s)\|_{L_2(\partial \Fr)}^2M_\alpha(\Gamma (y))(s)dm_r(s).
	\end{equation*}
	In addition, Corollary~\ref{tmn-vec-kernel}, the definition of \(\widetilde{K}_{N,m}\) \eqref{KmN-new}, and the orthogonality of martingale difference sequences give
	\begin{align*}
	    &\int_I\|\Phi (S_N(x))(s)\|_{L_2(\partial \Fr)}^2\Gamma (y)(s)dm_r(s)\\
        &
	=\int_{-1}^1\sum_{m=0}^N\bigg\|\int_{-1}^1\widetilde K_{N,m}(s,t)\Delta_m \Big(\xi\Big(\frac{\sqrt{2r-1}}{r}t\Big)\Big)w(t)dt\bigg\|_{L_2(\partial \Fr)}^2
                \, \Gamma (y)\Big(\frac{\sqrt{2r-1}}{r}s\Big) \, w(s) \, ds,
	\end{align*}
	and
	\begin{align*}
	     \int_{I}\|\Phi(x)\|_{L_2(\partial \Fr)}^2&\cM_\alpha(\Gamma (y))(s)dm_r(s)\\
         &=\int_{-1}^1\sum_{m=0}^\infty\bigg\|\Delta_m \Big(\xi\Big(\frac{\sqrt{2r-1}}{r}s\Big)\Big)\bigg\|_{L_2(\partial \Fr)}^2
                      \cM_\alpha(\Gamma (y))\Big(\frac{\sqrt{2r-1}}{r}s\Big) \, w(s) \, ds,
	\end{align*}
	where \(\xi(t)=\Phi(x)(\frac{\sqrt{2r-1}}{r}t).\)
	As a result, combining the above equalities, we see that the desired
weighted $L_2$ inequality in the theorem amounts to the following estimate
\[
  \int_{-1}^1 \Big\| \int_{-1}^1 \widetilde{K}_{N,m}(s,t) f(t) w(t) \, dt \Big\|_{L_2(\partial \Fr)}^2 g(s) w(s) \, ds
  \leq C_\alpha \int_{-1}^1 \|f(s)\|_{L_2(\partial \Fr)}^2 \cM_\alpha(g)(s) w(s) \, ds,
\]
where $0 \leq m \leq N$ and, for simplicity, we set
$f(t) = \Delta_m\!\bigl( \xi( \frac{\sqrt{2r-1}}{r} t ) \bigr)$ and
$g(s) = \Gamma(y)\!\bigl( \frac{\sqrt{2r-1}}{r} s \bigr)$.

\medskip
\noindent\textbf{Step 2:} In this step, we establish the desired
vector-valued weighted $L_2$ inequality reduced to at the end of
the previous step. By Proposition~\ref{AN-CN}, we have
\[
  \widetilde{K}_{N,m}(s,t) = \sum_{i=1}^3 a_{m,N}^i \, K_{N,m}^i(s,t), \quad s,t \in [-1,1],
\]
where $|a_{m,N}^i| < C$ for a universal constant $C > 0$.
It therefore suffices to show that, for each $i \in \{1,2,3\}$,
\begin{equation}\label{Ki-l2weighted}
        \int_{-1}^1\Big\|\int_{-1}^1K^i_{N,m}(s,t)f(t)w(t)dt\Big\|_{L_2(\partial \Fr)}^2g(s)w(s)ds\leq C_\alpha \int_{-1}^1\|f(s))\|_{L_2(\partial \Fr)}^2 \cM_\alpha(g)(s)w(s)ds.
    \end{equation}
    
    We claim that for every $f,g\colon [-1,1]\to\mathbb{C}$ with $f\in L_2([-1,1],\cM_\alpha(g)wds)$ and $g\geq 0$,  we have
	\begin{equation}\label{mainthmeq}   
		\int_{-1}^1\Big| \int_{-1}^1K_{N,m}^i(s,t) f(t)w(t)dt \Big|^2g(s)w(s)ds\le C_{\alpha}\int_{-1}^1 | f(s) |^2 \cM_\alpha(g)(s)w(s)ds,\quad i\in\{1,2,3\}.
	\end{equation}
    Once the claim is proved, by Fubini’s theorem, we have
    \begin{align*}
        \mathrm{LHS}\,\eqref{Ki-l2weighted}&=\int_{-1}^1\int_{\partial \Fr}\big|\int_{-1}^1K^i_{N,m}(s,t)f(t)w(t)dt\big|^2d\nu\ g(s)w(s)ds\\
        &=\int_{\partial \Fr}\int_{-1}^1\big|\int_{-1}^1K^i_{N,m}(s,t)f(t)w(t)dt\big|^2\ g(s)w(s)ds\ d\nu\\
        &\leq C_\alpha\int_{\partial \Fr}\int_{-1}^1\big|f(s)\big|^2\ \cM_\alpha(g)(s)w(s)ds\ d\nu\\
        &= C_\alpha\int_{-1}^1\|f(s)\|_{L_2(\partial \Fr)}^2\ \cM_\alpha(g)(s)w(s)ds,
    \end{align*}
	which proves \eqref{Ki-l2weighted}.   Now we prove the claim \eqref{mainthmeq}.

	\smallskip
    
	\noindent\textit{Case 1:} We first consider the case \(i=1\). By Proposition \ref{AN-CN}, we have
\[K_{N,m}^1(s,t)=\frac{ U_{m,N+1}(s)\bar{q}_{m,N}(t)}{s - t},\quad s,t\in[-1,1].\]
Then,
	\[
	\text{LHS }\eqref{mainthmeq}=\int_{-1}^1| H(\bar{q}_{m,N}fw)(s)|^2 |U_{m,N+1}(s)|^2g(s)w(s)ds,
	\] 
    here and in what follows, 
    \(H\) is the Hilbert transform defined by
    \[H(f)(s)=\operatorname{p.v.}\int_{-1}^1\frac{f(t)}{s-t}dt,\quad s\in[-1,1].\]
	For  \(\alpha>2\) and \(1<\rho\leq2\alpha/(\alpha+2)\), by the Coifman–Rochberg Theorem \cite[Theorem 8.33]{KLV2021}, we have \(\cM_\rho((U_{m,N+1})^2gw)\in A_2\) for any \(\rho>1\)  and \([\cM_\rho((U_{m,N+1})^2gw)]_{A_2}\lesssim \frac{\rho}{\rho-1}\). Hence, Lemma \ref{Qmn-qmn-w} (i) (ii) give
	\begin{equation*}
		\begin{aligned}
			\text{LHS }\eqref{mainthmeq}&\le \int_{-1}^1\left| H(q_{m,N}fw) (s)\right|^2\cM_\rho(|U_{m,N+1}|^2g \, w)(s) \, ds\\
			&\lesssim_\alpha \int_{-1}^1|q_{m,N}fw|^2(s) \, \cM_\rho(|U_{m,N+1}|^2 g \, w)(s) \, ds\\
			& \lesssim \int_{-1}^1\left| f(s) \right|^2 \, \kappa(s) \, \cM_\rho \Big( \frac{g}{\kappa} \Big)(s) \, w(s) \, ds\\
			& \lesssim_\alpha \int_{-1}^1\left| f(s) \right|^2 \, \cM_\alpha(g)(s) \, w(s) \, ds,
		\end{aligned}
	\end{equation*}
where the last inequality follows from Lemma \ref{kernel1} with \(\rho=(1+2\alpha/(\alpha+2))/2\).  
	
	\smallskip
    
	\noindent\textit{Case 2:} For the case \(i = 2\), by Proposition \ref{AN-CN}, we have
\[K_{N,m}^2(s,t)= \frac{q_{m,N}(s)\overline{U}_{m,N+1}(t)}{s-t},\quad s,t\in[-1,1].\]
    We similarly have
	\[
	\text{LHS }\eqref{mainthmeq}=\int_{-1}^1\left| H(\bar U_{m,N+1}fw) (s)\right|^2|q_{m,N}(s)|^2g(s)w(s)ds.
	\]
	Using Lemma \ref{Qmn-qmn-w} (ii) and Lemma \ref{kernel2}, we obtain
	\begin{equation*}
		\begin{aligned}
			\text{LHS }\eqref{mainthmeq}&\lesssim\int_{-1}^1\left|H(\bar U_{m,N+1}fw)(s)\right|^2(1-s^2)^{\frac{1}{2}}g(s)ds\\
			& \le \int_{-1}^1\left|H(\bar U_{m,N+1}fw)(s)\right|^2(1-s^2)^{\frac{1}{2}}\cM_\alpha(g)(s)ds\\
			&\lesssim_\alpha \int_{-1}^{1}|\bar U_{m,N+1}(s)|^2|f(s)|^2 w^2(s)(1-s^2)^{\frac{1}{2}}M_\alpha(g)(s)ds\\
			& \lesssim \int_{-1}^1 \left| f(s) \right|^2\cM_\alpha(g)(s)w(s)ds,
		\end{aligned}
	\end{equation*}
	where the last inequality is due to Lemma \ref{Qmn-qmn-w} (i).
	
	\smallskip
    
	\noindent\textit{Case 3:} For the case \(i=3\), by Proposition \ref{AN-CN}, we have
\[K_{N,m}^3(s,t)= U_{m,N+1}(s) \overline{U}_{m,N+1}(t),\quad s,t\in[-1,1].\]
    By the H\"older inequality,
	\begin{equation*}
		\begin{aligned}
			\text{LHS }\eqref{mainthmeq}&=\int_{-1}^1\Big| \int_{-1}^1\overline{U}_{m,N+1}(t)f(t)w(t)dt \Big|^2|U_{m,N+1}(s)|^2(s)g(s)w(s)ds\\
			&\leq \int_{-1}^1|\overline{U}_{m,N+1}(t)|^2w(t)dt\int_{-1}^1|f(t)|^2w(t)dt \int_{-1}^1|U_{m,N+1}(s)|^2(s)g(s)w(s)ds.
		\end{aligned}
	\end{equation*}
Proposition~\ref{Phi_QmnHm} (ii) and  Lemma \ref{Qmn-qmn-w} (i) give
	\begin{equation*}
\text{LHS }\eqref{mainthmeq}
			\lesssim \int_{-1}^1|f(t)|^2w(t)dt \int_{-1}^1g(s) (1-s^2)^{-\frac{1}{2}}ds.
	\end{equation*}
	Applying the H\"older inequality once more together with the fact the evenness of the function \(1-s^2\) and \(1-s^2>1-s\) for \(s\in[0,1]\), we deduce that
		\begin{equation*}
		\begin{aligned}
			\text{LHS }\eqref{mainthmeq}
			&\leq \int_{-1}^1|f(t)|^2 w(t)dt \Big(\int_{-1}^1|g(s)|^{\alpha}ds\Big)^{\frac{1}{\alpha}}\Big(\int_{-1}^1(1-s^2)^{-\frac{ \alpha'}{2}}ds\Big)^{\frac 1{\alpha'}}\\
			&\leq \int_{-1}^1|f(t)|^2 w(t)dt \Big(\int_{-1}^1|g(s)|^{\alpha}ds\Big)^{\frac{1}{\alpha}}\Big(2\int_{0}^1(1-s)^{-\frac{ \alpha'}{2}}ds\Big)^{\frac 1{\alpha'}}\\
			&\lesssim_\alpha \int_{-1}^1|f(t)|^2 \cM_\alpha g(t)w(t)dt,
		\end{aligned}
	\end{equation*}
		where the last inequality is due to \(-\frac {\alpha'}{2}>-1\).
	Moreover, the above inequality allows us to extend the assignment $x=\lambda(f)$ to all $x\in \cL \Fr$. The proof is complete.
\end{proof}

We are now positioned to complete the proof of Theorem~\ref{Tm-pp}.

\begin{proof}[Proof of Theorem \ref{Tm-pp}]
We only prove the column case, since the row case follows by replacing $x$ with $x^*$. By density and duality, it suffices to establish the estimate for \(x\in \cL \Fr\) in the range \(2<p<3\). We have
\begin{equation*}
	\begin{aligned}
	 \|S_N(x)\|_{{L_p}(\cL \Fr;\mathbb{E}_{\rad}^c)}^2
	 &=\tau\Big(\EE_\rad\big[S_N(x)^*S_N(x)\big] ^{\frac{p}{2}} \Big)^{\frac{2}{p}}\\
	&=\sup_y \tau\Big(\EE_\rad\big[S_N(x)^*S_N(x)\big] \, y \Big)\\
	&=\sup_y \tau\Big(|S_N(x)|^2 \, \EE_\rad(y) \Big),
	\end{aligned}
\end{equation*}
where the supremum is taken over all positive operators \(y\in L_{(\frac{p}{2})'}(\cL \Fr)\) of unit norm. Fix \(2<\alpha<\frac{2}{3}(\frac{p}{2})'\).
By Theorem B and H\"older's inequality,
\begin{align*}
    \|S_N(x)\|_{L_p(\mathcal L\Fr;\EE_\rad^c)}^2
        & \lesssim_\alpha \sup_y \tau\left( |x|^2 \, \cM_\alpha(\EE_\rad[y] ) \right) \\
        & = \sup_y \tau\left( \EE_\rad(|x|^2) \, \cM_\alpha(\EE_{\rad}[y]) \right) \\
        & \le \sup_y \|x\|_{L_p(\cL \Fr;\EE_\rad^c)}^2 \, \big\|\cM_\alpha(\EE_\rad(y))\|_{L_{(\frac{p}{2})'}(\cL \Fr)}.
\end{align*}
Since $\EE_\rad$ extends as a contraction to every $L_p$-space and
\[
    \cM_\alpha = \Gamma^{-1}\circ \cM_\alpha\circ\Gamma: L_{(\frac{p}{2})'}(\cL \Fr^\rad) \longrightarrow L_{(\frac{p}{2})'}(\cL \Fr^\rad)
\]
is bounded, as \(w\in A_s\) for every \(s>3/2\), we obtain the desired estimate. This completes the proof.
\end{proof}

\section{Non-uniform bounds for \texorpdfstring{$q$}{q}-balls}\label{sec: q-balls}

Here we will proceed with the proof of Theorem \ref{thm:qBall}. To that end, we recall some facts on the relationship between Fourier and Herz-Schur multipliers. When $(\Omega, \mu)$ is a measure space and $1 \leq p < \infty$, we will denote by $S_p(\Omega) \subseteq B(L_2 \Omega)$ its Schatten $p$-class, that is the noncommutative $L_p$-spaces associated with $B(L_2 \Omega)$ endowed with its natural semifinite trace $\mathrm{Tr}$. When $p=2$, this process gives the so-called Hilbert-Schmidt matrices and in that case we have an isometric isomorphism $S_2(\Omega) \cong L_2(\Omega \times \Omega)$ that associates to each operator $T \in S_2$ its integral kernel $k_T$. Given a \textit{Schur symbol} $M: \Omega \times \Omega \to \CC$, this allows us to define its associated Schur multiplier $S_M: S_2(\Omega) \to S_2(\Omega)$ as the operator sending $T \in S_2(G)$ to the operator associated with the kernel
\[
    M(\omega,\eta) \, k_T(\omega, \eta). 
\]
Observe that $S_M$, a priori, is only bounded over $S_2(\Omega)$. When it extends to a bounded operator from $S_2(\Omega) \cap S_p(\Omega)$ to $S_p(\Omega)$, we have a $S_p$-bounded Schur multiplier. In the same fashion, complete boundedness for $S_M$ can be defined. The interested reader can look up \cite{P2025} or the introductory sections of \cite{LD2011,Parcet2024Idempotentschur} for more information on $S_p$-bounded Schur multipliers. 

In the case of $\Omega = G$ being a locally compact second countable group and $m: G \to \CC$ being a bounded function, its associated Herz-Schur multiplier $H_m$ is the Schur multiplier whose symbol is given by $M(\gamma, \eta) = m(\gamma \, \eta^{-1})$. In our proof we will use the now well-known fact that Herz-Schur multipliers provide a lower bound for (complete bounded) Fourier multipliers, ie
\begin{equation}
    \label{eq:lower transference}
    \big\| H_m: S_p(G) \to S_p(G) \big\|_\cb \; \leq \; \big\| T_m: L_p(\cL G) \to L_p(\cL G) \big\|_\cb, \quad \text{ for every } 1 \leq p \leq \infty.
\end{equation}
There $L_p(\cL G)$ represents the noncommutative $L_p$-space of the regular von Neumann algebra of $G$. The above inequality was shown by \cite{NeuRic2011} in the discrete case and by \cite{CasSall2015} in the locally compact setting. Recall as well that the reverse inequality to \eqref{eq:lower transference} has been the subject of a lot of research. Initially it was only known for amenable groups \cite{NeuRic2011}, now there are recent counterexamples on the free group \cite{GPPW2026}. 

\begin{proof}[Proof of Theorem~\ref{thm:qBall}]
Fix $1<q<\infty$. By duality and interpolation with the completely contractive $L_2$-bounds, it suffices to consider the case of $2<p<\infty$. Similarly, by Fourier restriction it is enough to prove in the case of $r = 2$, that is $\FF_2 = \langle a, b\rangle$. Assume, by contradiction, that
\[
    \sup_{N\geq 1} \|S_N^q \colon L_p(\cL \Fr) \to L_p(\cL \Fr) \|_\cb = C < \infty.
\]
Using the lower bound of the Fourier--Schur transference \eqref{eq:lower transference} yields that the following symbols
\[
    K_N(\gamma,\eta)=\chi_{[0,1]} \Big( \frac1{N^q} \big|\gamma \, \eta^{-1} \big|_q^q \Big)
\]
induce uniformly bounded Schur multiplier in $S_p(\Fr)$. Restricting to the two-parameters subsets $X \subset \Fr$ given by $X = \big\{ a^n \, b^m : n, \, m \in \ZZ \big\}$ gives a uniformly bounded family of Schur multipliers
\begin{align*}
    \sup_{N} \big\| S_{F_N} \colon &S_p(X) \to S_p(X) \big\| \leq C,\\[-3pt]
        & \text{ where }\; F_N((n_1, m_1), (n_1, m_2)) = \chi_{[0,1]} \Big( \frac1{N^q} \big| a^{n_1}b^{m_1}(a^{n_2}b^{m_2})^{-1}\big|_q^q \Big).
\end{align*}
Now, let $\psi_N: \RR^2 \to \ZZ^2$ be the measurable transformation given by $\psi_N(x_1, x_2) = \big( {\lfloor N \, x_2\rfloor}, \, {\lfloor N \, x_1\rfloor} \big)$. Now, we get that $G_N(x, y) = F_N(\psi_N(x), \psi_N(y))$, where $x, y \in \RR^2$, are also uniformly $S_p(\RR^2)$-bounded Schur multipliers. To see that, just notice that there is a natural unitary isometry splitting $L_2(\RR^2)$ as $\ell_2(X) \otimes_2 L_2([0,\tfrac1{N}]^2)$ and that, conjugating under these unitaries just intertwines $S_{G_N}$ and $S_{F_N} \otimes \id_{S_p([0,\frac1{N}]^2)}$, therefore:
\[
    \big\| S_{G_N}\colon S_p(\RR^2) \to S_p(\RR^2) \big\| = \big\| S_{F_N} \colon S_p(X) \to S_p(X) \big\|_\cb.
\]
Notice as well that the left hand side is automatically completely bounded by \cite[Theorem~1.18]{LD2011}. Now, we have that the limit of the sequence $G_N$ as $N \to \infty$ converges to
\[
    G(x,y) = \chi_{[0,1]} \big( |x_1|^q+|x_2-y_2|^q+|y_1|^q \big),
\]
where $x = (x_1,x_2)$ and $y = (y_1, y_2)$. Thus, by the pointwise limit stability of $S_p$-bounded Schur multipliers, it follows that $S_G\colon S_p(\RR^2) \to S_p(\RR^2)$ is a (completely) bounded idempotent Schur multiplier supported in the domain
\[
    \Omega := \{ (x,y) \in \RR^2 \times \RR^2 :  \underbrace{|x_1|^q+|x_2-y_2|^q+|y_1|^q}_{H(x,y)} \leq 1 \big\},
\]
whose boundary is of $C^1$ class. As such, it is in the hypotheses of the characterization of idempotent Schur multipliers in \cite[Theorem~A]{Parcet2024Idempotentschur}. Seeing that it fails such characterization is enough to obtain a contradiction. To that end, fix $0<\varepsilon<1$ and put
\[
    y=(\varepsilon,0),
    \qquad
    x(t)=\big(t,v(t)\big),
    \qquad
    v(t)=(1-\varepsilon^q-t^q)^{1/q},
\]
for $0<t<(1-\varepsilon^q)^{\frac{1}{q}}$. Then $H(x(t),y)=1$ and
\[
    \nabla_x \, H(x(t),y)=q\big(t^{q-1},v(t)^{q-1}\big),
    \qquad
    \nabla_y \, H(x(t),y)=q\big(\varepsilon^{q-1},-v(t)^{q-1}\big).
\]
Near these points the boundary is smooth and transverse, since both partial gradients are nonzero. The cited theorem forces
$\nabla_y \, H(x(t),y)$ and $\nabla_y \, H(x(t'),y)$ to be parallel for sufficiently close $t,t'$. However, their first coordinates are
identical and nonzero, whereas their second coordinates differ whenever $t\neq t'$, because $q>1$ and $v$ is strictly decreasing.
This completes the proof.
\end{proof}

\textbf{Use of AI statement. } The mathematical contents of this article were completed by the early spring of 2026---before the advent of widely available IA systems capable of doing research-level mathematics---. Indeed, the results were presented in May 2026 in the \textit{Harmonic Analysis and von Neumann Algebras (HAVNA)}, Madrid. All the results have been written by humans and the use of LLM's has been circumscribed to help in polishing the language, editing and reference keeping. However, interactions around the present article with IA tools have produced two follow up results. The first is a clever use of the filtration $(\cA_m)_m$ in $\Fr$ to produce the first examples of Herz-Schur multipliers in $S_p(\Fr)$ that are not Fourier multipliers in $L_p(\cL \Fr)$; see \cite{GPPW2026}. The second \cite{GLPPWZ2026} would be a succeeding text in which we will develop, among other results
 \begin{enumerate}[label={\small (\alph*)}]
     \item That being bounded between mixed-norm spaces is, for radial symbols $a: \Fr \to \CC$, a strictly weaker notion than being (completely) $L^p$-bounded, ie
     \[
        \big\| T_a: L_p(\cL \Fr;\EE_\rad^c) \to L_p(\cL \Fr;\EE_\rad^c) \big\|_\cb
        \; \leq \; \big\| T_a: L_p(\cL \Fr) \to L_p(\cL \Fr) \big\|_\cb,
     \]
    and the above result fails for non-radial symbols. As a consequence, using complex interpolation techniques, Bochner-Riesz multipliers are bounded on mixed-norm spaces for the range of exponents $p$ corresponding to dimension $3$.

    \item We also obtain a very striking result neatly connecting the mixed norm theory with (central) spectral multipliers in $SU(2)$. Indeed, if $a(\gamma) = \mathring{a}(|\gamma|)$ is a radial multiplier we have that
    \[
        \big\| T_a: L_p(\cL \Fr;\EE_\rad^c) \to L_p(\cL \Fr;\EE_\rad^c) \big\|
        \simeq_{(p)} \cR_p \big\{ S_{\mathring{a}}^{(m)}: L_p(SU(2)) \to L_p(SU(2)) \big\},
    \]
    where $S_{\mathring{a}}^{(m)}$ is the $m$-shifted spectral multiplier sending the characters of each irreducible representation $\chi_k$ to $\mathring{a}(m+k) \, \chi_k$. This connection between the free group and $SU(2)$ extends earlier results due to Biane \cite{Biane2008RandomWalks}.
    
    \item Lastly, we show that in the Herz range $p \in (\tfrac{3}{2}, 3)$, a combination of $R$-boundedness and $A_p$-weight techniques gives that it is enough to take the single, unshifted, multiplier $S_a^{(0)}$.
 \end{enumerate}

\textbf{Acknowledgment.} A. González-Pérez was supported by the Ramón y Cajal grant \texttt{RYC2022}\texttt{-}\texttt{037045-I} (Ministry of Science, Spain) while in a research stay at Xiamen University. J. Parcet was supported by the European Research Council: Advanced Grant SMARP \texttt{101263125} and  Severo Ochoa Grant \texttt{CEX2023-001347-S}. The first three-named authors were also supported in part by Spanish Grant \texttt{PID2022-141354NB-I00}. J. Pérez-García and S. Wang were supported by the NSF of China: \texttt{No.12671160}, \texttt{No.W2611005}, \texttt{No.W2441002}, \texttt{No.12031004}.  

\bibliographystyle{alpha}

\bibliography{refs}

@book{FP1983,
  author        = {Fig\`a-Talamanca, A. and Picardello, M.},
  title         = {Harmonic analysis on free groups},
  series        = {Lecture Notes in Pure and Applied Mathematics},
  volume        = {87},
  pages         = {viii+145},
  publisher     = {Marcel Dekker, Inc., New York},
  year          = {1983},
  isbn          = {0-8247-7042-0},
  mrclass       = {43-02 (22D10 46Lxx)},
  mrnumber      = {710827},
  mrreviewer    = {M.\ Bo\.zejko},
}

@article{Co1984,
  author        = {C{\'o}rdoba, A.},
  title         = {An integral inequality for the disc multiplier},
  journal       = {Proc. Amer. Math. Soc.},
  fjournal      = {Proceedings of the American Mathematical Society},
  volume        = {92},
  number        = {3},
  pages         = {407--408},
  year          = {1984},
  doi           = {10.1090/S0002-9939-1984-0759664-6},
}

@article{Pytlik1981,
  author        = {Pytlik, T.},
  title         = {Radial functions on free groups and a decomposition of the regular representation
                  into irreducible components},
  journal       = {J. Reine Angew. Math.},
  fjournal      = {Journal f{\"u}r die reine und angewandte Mathematik},
  volume        = {326},
  pages         = {124--135},
  year          = {1981},
  doi           = {10.1515/crll.1981.326.124},
  url           = {https://doi.org/10.1515/crll.1981.326.124},
  mrnumber      = {0622348},
  zbl           = {0464.22004},
}

@book{CheCowJoJulgVal2001,
  author        = {Cherix, P. A. and Cowling, M. and Jolissaint, P. and Julg, P. and Valette, A.},
  title         = {Groups with the {H}aagerup property},
  volume        = {197},
  pages         = {viii+126},
  publisher     = {Birkh\"auser Verlag, Basel},
  year          = {2001},
  doi           = {10.1007/978-3-0348-8237-8},
  isbn          = {3-7643-6598-6},
  url           = {http://dx.doi.org/10.1007/978-3-0348-8237-8},
  mrreviewer    = {Tullio G. Ceccherini-Silberstein},
}

@book{Hyt2017BookII,
  author        = {Hyt{\"o}nen, T. and van Neerven, J. and Veraar, M. and Weis, L.},
  title         = {Analysis in Banach Spaces. Volume {II}: Probabilistic Methods and Operator
                  Theory},
  series        = {Ergebnisse der Mathematik und ihrer Grenzgebiete. 3. Folge. A Series of Modern
                  Surveys in Mathematics},
  volume        = {67},
  publisher     = {Springer},
  address       = {Cham},
  year          = {2017},
  doi           = {10.1007/978-3-319-69808-3},
  isbn          = {978-3-319-69808-3},
}

@article{ChuahLiuMei2023,
  author        = {Chuah, C. Y. and Liu, Z. and Mei, T.},
  title         = {Characterization of positive definite, radial functions on free groups},
  journal       = {Proc. Amer. Math. Soc. Ser. B},
  fjournal      = {Proceedings of the American Mathematical Society, Series B},
  volume        = {10},
  pages         = {168--181},
  year          = {2023},
  doi           = {10.1090/bproc/158},
  eprint        = {2112.00630},
  archivePrefix = {arXiv},
  primaryClass  = {math.FA},
}

@article{NeuRic2011,
  author        = {Neuwirth, S. and Ricard, {\'E}.},
  title         = {Transfer of {Fourier} multipliers into {Schur} multipliers and sumsets in a
                  discrete group},
  journal       = {Canad. J. Math.},
  fjournal      = {Canadian Journal of Mathematics},
  volume        = {63},
  number        = {5},
  pages         = {1161--1187},
  year          = {2011},
  doi           = {10.4153/CJM-2011-053-9},
  language      = {English},
  reviewer      = {{Vladimir V. Peller (East Lansing)}},
}

@article{Kesten1959,
  author        = {Kesten, H.},
  title         = {Symmetric random walks on groups},
  journal       = {Trans. Amer. Math. Soc.},
  fjournal      = {Transactions of the American Mathematical Society},
  volume        = {92},
  number        = {2},
  pages         = {336--354},
  year          = {1959},
  doi           = {10.1090/S0002-9947-1959-0109367-6},
  url           = {https://doi.org/10.1090/S0002-9947-1959-0109367-6},
}

@misc{GPPW2026,
  author        = {Gonz{\'a}lez-P{\'e}rez, A. and Parcet, J. and P{\'e}rez Garc{\'i}a, J. and Wang,
                  S.},
  title         = {Radial {Schur} multipliers on trees and the failure of {Fourier--Schur}
                  transference},
  year          = {2026},
  note          = {Preprint},
  doi           = {10.48550/arXiv.2610.15881},
  eprint        = {2610.15881},
  archivePrefix = {arXiv},
  primaryClass  = {math.FA},
  url           = {https://arxiv.org/abs/2610.15881},
}

@article{KVW2016,
  author        = {Kwapie{\'n}, S. and Veraar, M. and Weis, L.},
  title         = {{$R$}-boundedness versus {$\gamma$}-boundedness},
  journal       = {Ark. Mat.},
  fjournal      = {Arkiv f{\"o}r Matematik},
  volume        = {54},
  number        = {1},
  pages         = {125--145},
  year          = {2016},
  doi           = {10.1007/s11512-015-0223-1},
  eprint        = {1404.7328},
  archivePrefix = {arXiv},
  primaryClass  = {math.FA},
}

@article{CasSall2015,
  author        = {Caspers, M. and de la Salle, M.},
  title         = {Schur and {F}ourier multipliers of an amenable group acting on non-commutative
                  {$L_p$}-spaces},
  journal       = {Trans. Amer. Math. Soc.},
  fjournal      = {Transactions of the American Mathematical Society},
  volume        = {367},
  number        = {10},
  pages         = {6997--7013},
  year          = {2015},
  doi           = {10.1090/S0002-9947-2015-06281-3},
}

@article{Parcet2024Idempotentschur,
  author        = {Parcet, J. and de la Salle, M. and Tablate, E.},
  title         = {The local geometry of idempotent {Schur} multipliers},
  journal       = {Forum Math. Pi},
  fjournal      = {Forum of Mathematics, Pi},
  volume        = {13},
  pages         = {e14},
  eid           = {e14},
  year          = {2025},
  doi           = {10.1017/fmp.2025.6},
  eprint        = {2312.02895},
  archivePrefix = {arXiv},
}

@article{Bozejko1986,
  author        = {Bo{\.z}ejko, M.},
  title         = {Positive definite functions on the free group and the noncommutative {Riesz}
                  product},
  journal       = {Boll. Un. Mat. Ital. A (6)},
  fjournal      = {Bollettino della Unione Matematica Italiana. A. Serie VI},
  volume        = {5},
  number        = {1},
  pages         = {13--21},
  year          = {1986},
}

@article{Cohen1982Operator,
  author        = {Cohen, J. M.},
  title         = {Operator norms on free groups},
  journal       = {Boll. Un. Mat. Ital. B (6)},
  fjournal      = {Bollettino dell'Unione Matematica Italiana. B. Serie VI},
  volume        = {1},
  number        = {3},
  pages         = {1055--1065},
  year          = {1982},
  zbl           = {0518.46050},
}

@incollection{CohenDeMichele1982,
  author        = {Cohen, J. M. and De Michele, L.},
  title         = {The radial {Fourier--Stieltjes} algebra of free groups},
  booktitle     = {Operator Algebras and {$K$}-Theory},
  series        = {Contemporary Mathematics},
  volume        = {10},
  pages         = {33--40},
  publisher     = {American Mathematical Society},
  address       = {Providence, RI},
  year          = {1982},
  doi           = {10.1090/conm/010/658507},
  url           = {https://doi.org/10.1090/conm/010/658507},
  mrnumber      = {658507},
}

@article{Car1985,
  author        = {Carbery, A.},
  title         = {A weighted inequality for the maximal {Bochner--Riesz} operator on
                  {$\mathbb{R}^{2}$}},
  journal       = {Trans. Amer. Math. Soc.},
  fjournal      = {Transactions of the American Mathematical Society},
  volume        = {287},
  number        = {2},
  pages         = {673--680},
  year          = {1985},
  doi           = {10.1090/S0002-9947-1985-0768732-X},
}

@incollection{RdF1982Weighted,
  author        = {Rubio de Francia, J. L.},
  editor        = {Ricci, F. and Weiss, G.},
  title         = {Weighted norm inequalities and vector valued inequalities},
  booktitle     = {Harmonic Analysis (Minneapolis, Minn., 1981)},
  series        = {Lecture Notes in Mathematics},
  volume        = {908},
  pages         = {86--101},
  publisher     = {Springer-Verlag},
  address       = {Berlin},
  year          = {1982},
  doi           = {10.1007/BFb0093282},
}

@article{Arnaud1980,
  author        = {Arnaud, J.-P.},
  title         = {Fonctions sph{\'e}riques et fonctions d{\'e}finies positives sur l'arbre
                  homog{\`e}ne},
  journal       = {C. R. Acad. Sci. Paris S{\'e}r. A-B},
  fjournal      = {Comptes Rendus Hebdomadaires des S{\'e}ances de l'Acad{\'e}mie des Sciences.
                  S{\'e}ries A et B},
  volume        = {290},
  number        = {2},
  pages         = {A99--A101},
  year          = {1980},
  mrnumber      = {563948},
}

@article{HaagerupKnudby2015,
  author        = {Haagerup, U. and Knudby, S.},
  title         = {A {L{\'e}vy-Khinchin} formula for free groups},
  journal       = {Proc. Amer. Math. Soc.},
  fjournal      = {Proceedings of the American Mathematical Society},
  volume        = {143},
  number        = {4},
  pages         = {1477--1489},
  year          = {2015},
  doi           = {10.1090/S0002-9939-2014-12466-X},
  eprint        = {1304.5763},
  archivePrefix = {arXiv},
  primaryClass  = {math.GR},
  url           = {https://doi.org/10.1090/S0002-9939-2014-12466-X},
}

@incollection{Cartier1973Geometrie,
  author        = {Cartier, P.},
  title         = {G{\'e}om{\'e}trie et analyse sur les arbres},
  booktitle     = {S{\'e}minaire {Bourbaki}, vol. 1971/72, expos{\'e}s 400--417},
  series        = {Lecture Notes in Mathematics},
  volume        = {317},
  pages         = {123--140},
  publisher     = {Springer-Verlag},
  address       = {Berlin, Heidelberg},
  year          = {1973},
  note          = {Expos{\'e} no. 407, February 1972},
  doi           = {10.1007/BFb0069280},
  url           = {https://www.numdam.org/item/SB_1971-1972__14__123_0/},
}

@article{P1948,
  author        = {Pollard, H.},
  title         = {The mean convergence of orthogonal series. {II}},
  journal       = {Trans. Amer. Math. Soc.},
  fjournal      = {Transactions of the American Mathematical Society},
  volume        = {63},
  pages         = {355--367},
  year          = {1948},
  doi           = {10.2307/1990435},
  issn          = {0002-9947,1088-6850},
  url           = {https://doi.org/10.2307/1990435},
  mrclass       = {42.4X},
  mrnumber      = {23941},
  mrreviewer    = {E.\ G.\ Kogbetliantz},
}

@incollection{BF2006,
  author        = {Bo\.zejko, M. and Fendler, G.},
  title         = {A note on certain partial sum operators},
  booktitle     = {Quantum probability},
  series        = {Banach Center Publ.},
  volume        = {73},
  pages         = {117--125},
  publisher     = {Polish Acad. Sci. Inst. Math., Warsaw},
  year          = {2006},
  doi           = {10.4064/bc73-0-7},
  url           = {https://doi.org/10.4064/bc73-0-7},
  mrclass       = {46L51 (42C10)},
  mrnumber      = {2423120},
}

@book{KLV2021,
  author        = {Kinnunen, J. and Lehrbäck, J. and Vähäkangas, A.},
  title         = {Maximal Function Methods for Sobolev Spaces},
  series        = {Mathematical Surveys and Monographs},
  publisher     = {American Mathematical Society},
  year          = {2021},
  isbn          = {9781470465759},
  url           = {https://books.google.co.jp/books?id=-vU8EAAAQBAJ},
  lccn          = {2021013318},
}

@article{JP2010,
  author        = {Junge, M. and Parcet, J.},
  title         = {Mixed-norm inequalities and operator space {$L_p$} embedding theory},
  journal       = {Mem. Amer. Math. Soc.},
  fjournal      = {Memoirs of the American Mathematical Society},
  volume        = {203},
  number        = {953},
  pages         = {vi+155},
  year          = {2010},
  doi           = {10.1090/S0065-9266-09-00570-5},
  isbn          = {978-0-8218-4655-1},
  issn          = {0065-9266,1947-6221},
  url           = {https://doi.org/10.1090/S0065-9266-09-00570-5},
  mrclass       = {46L52 (46B09 46L07 46L09 46L51)},
  mrnumber      = {2589944},
  mrreviewer    = {Christian\ Le Merdy},
}

@article{Co1979,
  author        = {C{\'o}rdoba, A.},
  title         = {A note on {B}ochner-{R}iesz operators},
  journal       = {Duke Math. J.},
  fjournal      = {Duke Mathematical Journal},
  volume        = {46},
  number        = {3},
  pages         = {505--511},
  year          = {1979},
  issn          = {0012-7094,1547-7398},
  url           = {http://projecteuclid.org/euclid.dmj/1077313571},
  mrclass       = {42B15 (42B25)},
  mrnumber      = {544242},
  mrreviewer    = {Richard\ Bagby},
}

@article{CS1972,
  author        = {Carleson, L. and Sj\"{o}lin, P.},
  title         = {Oscillatory integrals and a multiplier problem for the disc},
  journal       = {Studia Math.},
  fjournal      = {Studia Mathematica},
  volume        = {44},
  pages         = {287--299. (errata insert)},
  year          = {1972},
  doi           = {10.4064/sm-44-3-287-299},
  issn          = {0039-3223,1730-6337},
  url           = {https://doi.org/10.4064/sm-44-3-287-299},
  mrclass       = {42A92 (42A18)},
  mrnumber      = {361607},
  mrreviewer    = {B.\ Muckenhoupt},
}

@article{Fe1972,
  author        = {Fefferman, C.},
  title         = {A note on spherical summation multipliers},
  journal       = {Israel J. Math.},
  fjournal      = {Israel Journal of Mathematics},
  volume        = {15},
  pages         = {44--52},
  year          = {1973},
  doi           = {10.1007/BF02771772},
  issn          = {0021-2172},
  url           = {https://doi.org/10.1007/BF02771772},
  mrclass       = {42A92},
  mrnumber      = {320624},
  mrreviewer    = {R.\ A.\ Hunt},
}

@book{Zy2002,
  author        = {Zygmund, A.},
  title         = {Trigonometric series. {V}ol. {I}, {II}},
  series        = {Cambridge Mathematical Library},
  edition       = {Third},
  pages         = {xii; Vol. I: xiv+383 pp.; Vol. II: viii+364},
  publisher     = {Cambridge University Press, Cambridge},
  year          = {2002},
  note          = {With a foreword by Robert A. Fefferman},
  isbn          = {0-521-89053-5},
  mrclass       = {01A75 (42-02)},
  mrnumber      = {1963498},
}

@article{Fe1971,
  author        = {Fefferman, C.},
  title         = {The multiplier problem for the ball},
  journal       = {Ann. of Math. (2)},
  fjournal      = {Annals of Mathematics. Second Series},
  volume        = {94},
  pages         = {330--336},
  year          = {1971},
  doi           = {10.2307/1970864},
  issn          = {0003-486X},
  url           = {https://doi.org/10.2307/1970864},
  mrclass       = {42A18 (42A92 47B99)},
  mrnumber      = {296602},
  mrreviewer    = {R.\ Larsen},
}

@inproceedings{St1979,
  author        = {Stein, E.},
  title         = {Some problems in harmonic analysis},
  booktitle     = {Harmonic analysis in {E}uclidean spaces ({P}roc. {S}ympos. {P}ure {M}ath.,
                  {W}illiams {C}oll., {W}illiamstown, {M}ass., 1978), {P}art 1},
  series        = {Proc. Sympos. Pure Math., XXXV, Part 1},
  pages         = {3--20},
  publisher     = {Amer. Math. Soc., Providence, RI},
  year          = {1979},
  isbn          = {0-8218-1438-9},
  mrclass       = {42B15 (32A35)},
  mrnumber      = {545235},
  mrreviewer    = {R.\ A.\ Hunt},
}

@article{Co1977,
  author        = {C{\'o}rdoba, A.},
  title         = {The {K}akeya maximal function and the spherical summation multipliers},
  journal       = {Amer. J. Math.},
  fjournal      = {American Journal of Mathematics},
  volume        = {99},
  number        = {1},
  pages         = {1--22},
  year          = {1977},
  doi           = {10.2307/2374006},
  issn          = {0002-9327,1080-6377},
  url           = {https://doi.org/10.2307/2374006},
  mrclass       = {42A18},
  mrnumber      = {447949},
  mrreviewer    = {A.\ Schwartz},
}

@article{Mo1990,
  author        = {Mockenhaupt, G.},
  title         = {On radial weights for the spherical summation operator},
  journal       = {J. Funct. Anal.},
  fjournal      = {Journal of Functional Analysis},
  volume        = {91},
  number        = {1},
  pages         = {174--181},
  year          = {1990},
  doi           = {10.1016/0022-1236(90)90051-L},
  issn          = {0022-1236,1096-0783},
  url           = {https://doi.org/10.1016/0022-1236(90)90051-L},
  mrclass       = {42B10 (42B15)},
  mrnumber      = {1054117},
  mrreviewer    = {Kenneth\ F.\ Andersen},
}

@incollection{CHS1991,
  author        = {Carbery, A. and Hern\'{a}ndez, E. and Soria, F.},
  title         = {Estimates for the {K}akeya maximal operator on radial functions in {${\bf
                  R}^n$}},
  booktitle     = {Harmonic analysis ({S}endai, 1990)},
  series        = {ICM-90 Satell. Conf. Proc.},
  pages         = {41--50},
  publisher     = {Springer, Tokyo},
  year          = {1991},
  isbn          = {4-431-70084-6},
  mrclass       = {42B25},
  mrnumber      = {1261427},
  mrreviewer    = {Stephen\ Buckley},
}

@article{Co1989,
  author        = {C{\'o}rdoba, A.},
  title         = {The disc multiplier},
  journal       = {Duke Math. J.},
  fjournal      = {Duke Mathematical Journal},
  volume        = {58},
  number        = {1},
  pages         = {21--29},
  year          = {1989},
  doi           = {10.1215/S0012-7094-89-05802-X},
  issn          = {0012-7094,1547-7398},
  url           = {https://doi.org/10.1215/S0012-7094-89-05802-X},
  mrclass       = {42B15 (47B38)},
  mrnumber      = {1016411},
  mrreviewer    = {Shan\ Zhen\ Lu},
}

@article{Ru1989,
  author        = {Rubio de Francia, J.},
  title         = {Transference principles for radial multipliers},
  journal       = {Duke Math. J.},
  fjournal      = {Duke Mathematical Journal},
  volume        = {58},
  number        = {1},
  pages         = {1--19},
  year          = {1989},
  doi           = {10.1215/S0012-7094-89-05801-8},
  issn          = {0012-7094,1547-7398},
  url           = {https://doi.org/10.1215/S0012-7094-89-05801-8},
  mrclass       = {42B15 (42B25)},
  mrnumber      = {1016410},
  mrreviewer    = {Norman\ J.\ Weiss},
}

@article{Be2014,
  author        = {Bennett, J.},
  title         = {Optimal control of singular {F}ourier multipliers by maximal operators},
  journal       = {Anal. PDE},
  fjournal      = {Analysis \& PDE},
  volume        = {7},
  number        = {6},
  pages         = {1317--1338},
  year          = {2014},
  doi           = {10.2140/apde.2014.7.1317},
  issn          = {2157-5045,1948-206X},
  url           = {https://doi.org/10.2140/apde.2014.7.1317},
  mrclass       = {42B15 (42B20 42B25 42B37)},
  mrnumber      = {3270165},
  mrreviewer    = {Petr\ Honz\'{\i}k},
}

@article{CRS1992,
  author        = {Carbery, A. and Romera, E. and Soria, F.},
  title         = {Radial weights and mixed norm inequalities for the disc multiplier},
  journal       = {J. Funct. Anal.},
  fjournal      = {Journal of Functional Analysis},
  volume        = {109},
  number        = {1},
  pages         = {52--75},
  year          = {1992},
  doi           = {10.1016/0022-1236(92)90011-7},
  issn          = {0022-1236,1096-0783},
  url           = {https://doi.org/10.1016/0022-1236(92)90011-7},
  mrclass       = {42B15 (42B10 42B20 42B25)},
  mrnumber      = {1183604},
  mrreviewer    = {Kenneth\ F.\ Andersen},
}

@article{He1954,
  author        = {Herz, C.},
  title         = {On the mean inversion of {F}ourier and {H}ankel transforms},
  journal       = {Proc. Natl. Acad. Sci. USA},
  fjournal      = {Proceedings of the National Academy of Sciences of the United States of America},
  volume        = {40},
  pages         = {996--999},
  year          = {1954},
  doi           = {10.1073/pnas.40.10.996},
  issn          = {0027-8424},
  url           = {https://doi.org/10.1073/pnas.40.10.996},
  mrclass       = {42.4X},
  mrnumber      = {63477},
  mrreviewer    = {I.\ I.\ Hirschman, Jr.},
}

@article{Co2016,
  author        = {C{\'o}rdoba, A.},
  title         = {Singular integrals, maximal functions and {F}ourier restriction to spheres: the
                  disk multiplier revisited},
  journal       = {Adv. Math.},
  fjournal      = {Advances in Mathematics},
  volume        = {290},
  pages         = {208--235},
  year          = {2016},
  doi           = {10.1016/j.aim.2015.11.039},
  issn          = {0001-8708,1090-2082},
  url           = {https://doi.org/10.1016/j.aim.2015.11.039},
  mrclass       = {42B20 (42B25)},
  mrnumber      = {3451922},
  mrreviewer    = {Norbert\ Ortner},
}

@article{DH1985,
  author        = {De Canni\`ere, J. and Haagerup, U.},
  title         = {Multipliers of the {F}ourier algebras of some simple {L}ie groups and their
                  discrete subgroups},
  journal       = {Amer. J. Math.},
  fjournal      = {American Journal of Mathematics},
  volume        = {107},
  number        = {2},
  pages         = {455--500},
  year          = {1985},
  doi           = {10.2307/2374423},
  issn          = {0002-9327,1080-6377},
  url           = {https://doi.org/10.2307/2374423},
  mrclass       = {43A22},
  mrnumber      = {784292},
  mrreviewer    = {Pierre\ Eymard},
}

@article{CH1989,
  author        = {Cowling, M. and Haagerup, U.},
  title         = {Completely bounded multipliers of the {F}ourier algebra of a simple {L}ie group
                  of real rank one},
  journal       = {Invent. Math.},
  fjournal      = {Inventiones Mathematicae},
  volume        = {96},
  number        = {3},
  pages         = {507--549},
  year          = {1989},
  doi           = {10.1007/BF01393695},
  issn          = {0020-9910,1432-1297},
  url           = {https://doi.org/10.1007/BF01393695},
  mrclass       = {22E27 (22E40 43A22 43A80 46J99)},
  mrnumber      = {996553},
  mrreviewer    = {Christopher\ Meaney},
}

@article{Haa1978,
  author        = {Haagerup, U.},
  title         = {An example of a nonnuclear {$C\sp{\ast} $}-algebra, which has the metric
                  approximation property},
  journal       = {Invent. Math.},
  fjournal      = {Inventiones Mathematicae},
  volume        = {50},
  number        = {3},
  pages         = {279--293},
  year          = {1978/79},
  doi           = {10.1007/BF01410082},
  issn          = {0020-9910,1432-1297},
  url           = {https://doi.org/10.1007/BF01410082},
  mrclass       = {46L05 (22D35 43A35)},
  mrnumber      = {520930},
  mrreviewer    = {Ole\ A.\ Nielsen},
}

@article{LD2011,
  author        = {Lafforgue, V. and de la Salle, M.},
  title         = {Noncommutative {$L^p$}-spaces without the completely bounded approximation
                  property},
  journal       = {Duke Math. J.},
  fjournal      = {Duke Mathematical Journal},
  volume        = {160},
  number        = {1},
  pages         = {71--116},
  year          = {2011},
  doi           = {10.1215/00127094-1443478},
  issn          = {0012-7094,1547-7398},
  url           = {https://doi.org/10.1215/00127094-1443478},
  mrclass       = {46L07 (46B28 46L52)},
  mrnumber      = {2838352},
}

@article{MR2017,
  author        = {Mei, T. and Ricard, {\'E}.},
  title         = {Free {H}ilbert transforms},
  journal       = {Duke Math. J.},
  fjournal      = {Duke Mathematical Journal},
  volume        = {166},
  number        = {11},
  pages         = {2153--2182},
  year          = {2017},
  doi           = {10.1215/00127094-2017-0007},
  issn          = {0012-7094,1547-7398},
  url           = {https://doi.org/10.1215/00127094-2017-0007},
  mrclass       = {46L07 (46L52 46L54)},
  mrnumber      = {3694567},
  mrreviewer    = {Daniele\ Puglisi},
}

@inproceedings{P2025,
  author        = {Parcet, J.},
  title         = {The impact of {S}chur multipliers in harmonic analysis and operator algebras},
  booktitle     = {International Congress of Mathematicians 2026},
  pages         = {730--750},
  organization  = {SIAM},
  year          = {2026},
}

@article{K1984,
  author        = {Kuhn, G.},
  title         = {Convergence of {F}ourier series expansion related to free groups},
  journal       = {Proc. Amer. Math. Soc.},
  fjournal      = {Proceedings of the American Mathematical Society},
  volume        = {92},
  number        = {1},
  pages         = {31--36},
  year          = {1984},
  doi           = {10.2307/2045147},
  issn          = {0002-9939,1088-6826},
  url           = {https://doi.org/10.2307/2045147},
  mrclass       = {43A50 (20E05 33A45 42C10)},
  mrnumber      = {749884},
  mrreviewer    = {B.\ Dreseler},
}

@article{JNRX2004,
  author        = {Junge, M. and Nielsen, N. and Ruan, Z. and Xu, Q.},
  title         = {{$\mathscr C\mathscr O\mathscr L_p$} spaces---the local structure of
                  non-commutative {$L_p$} spaces},
  journal       = {Adv. Math.},
  fjournal      = {Advances in Mathematics},
  volume        = {187},
  number        = {2},
  pages         = {257--319},
  year          = {2004},
  doi           = {10.1016/j.aim.2003.08.010},
  issn          = {0001-8708,1090-2082},
  url           = {https://doi-org-s-224.libdb.csu.edu.cn/10.1016/j.aim.2003.08.010},
  mrclass       = {46L52 (46L07)},
  mrnumber      = {2078339},
  mrreviewer    = {Christian\ Le Merdy},
}

@article{JMP2018,
  author        = {Junge, M. and Mei, T. and Parcet, J.},
  title         = {Noncommutative {R}iesz transforms-dimension free bounds and {F}ourier
                  multipliers},
  journal       = {J. Eur. Math. Soc. (JEMS)},
  fjournal      = {Journal of the European Mathematical Society (JEMS)},
  volume        = {20},
  number        = {3},
  pages         = {529--595},
  year          = {2018},
  doi           = {10.4171/JEMS/773},
  issn          = {1435-9855,1435-9863},
  url           = {https://doi-org-s-224.libdb.csu.edu.cn/10.4171/JEMS/773},
  mrclass       = {43A50 (42B15 46L52)},
  mrnumber      = {3776274},
  mrreviewer    = {Serap\ \"{O}ztop},
}

@article{JMP2014,
  author        = {Junge, M. and Mei, T. and Parcet, J.},
  title         = {Smooth {F}ourier multipliers on group von {N}eumann algebras},
  journal       = {Geom. Funct. Anal.},
  fjournal      = {Geometric and Functional Analysis},
  volume        = {24},
  number        = {6},
  pages         = {1913--1980},
  year          = {2014},
  doi           = {10.1007/s00039-014-0307-2},
  issn          = {1016-443X,1420-8970},
  url           = {https://doi-org-s-224.libdb.csu.edu.cn/10.1007/s00039-014-0307-2},
  mrclass       = {42B30 (42B35 46L52)},
  mrnumber      = {3283931},
  mrreviewer    = {Yasuo\ Komori-Furuya},
}

@article{GJP2017,
  author        = {Gonz\'{a}lez-P\'{e}rez, A. and Junge, M. and Parcet, J.},
  title         = {Smooth {F}ourier multipliers in group algebras via {S}obolev dimension},
  journal       = {Ann. Sci. \'{E}c. Norm. Sup\'{e}r. (4)},
  fjournal      = {Annales Scientifiques de l'\'{E}cole Normale Sup\'{e}rieure. Quatri\`eme
                  S\'{e}rie},
  volume        = {50},
  number        = {4},
  pages         = {879--925},
  year          = {2017},
  doi           = {10.24033/asens.2334},
  issn          = {0012-9593,1873-2151},
  url           = {https://doi-org-s-224.libdb.csu.edu.cn/10.24033/asens.2334},
  mrclass       = {42B15 (42B25 43A25 43A32)},
  mrnumber      = {3679616},
  mrreviewer    = {Li\ Zhong\ Peng},
}

@article{OP2010,
  author        = {Ozawa, N. and Popa, S.},
  title         = {On a class of {${\rm II}_1$} factors with at most one {C}artan subalgebra},
  journal       = {Ann. of Math. (2)},
  fjournal      = {Annals of Mathematics. Second Series},
  volume        = {172},
  number        = {1},
  pages         = {713--749},
  year          = {2010},
  doi           = {10.4007/annals.2010.172.713},
  issn          = {0003-486X,1939-8980},
  url           = {https://doi.org/10.4007/annals.2010.172.713},
  mrclass       = {46L10 (37A20)},
  mrnumber      = {2680430},
  mrreviewer    = {Stuart\ A.\ White},
}

@article{PV2014,
  author        = {Popa, S. and Vaes, S.},
  title         = {Unique {C}artan decomposition for {$\rm II_1$} factors arising from arbitrary
                  actions of free groups},
  journal       = {Acta Math.},
  fjournal      = {Acta Mathematica},
  volume        = {212},
  number        = {1},
  pages         = {141--198},
  year          = {2014},
  doi           = {10.1007/s11511-014-0110-9},
  issn          = {0001-5962,1871-2509},
  url           = {https://doi.org/10.1007/s11511-014-0110-9},
  mrclass       = {46L36 (37A05 46L10 46L37)},
  mrnumber      = {3179609},
  mrreviewer    = {Liguang\ Wang},
}

@incollection{PX2003,
  author        = {Pisier, G. and Xu, Q.},
  title         = {Non-commutative {$L^p$}-spaces},
  booktitle     = {Handbook of the geometry of {B}anach spaces, {V}ol.\ 2},
  pages         = {1459--1517},
  publisher     = {North-Holland, Amsterdam},
  year          = {2003},
  doi           = {10.1016/S1874-5849(03)80041-4},
  isbn          = {0-444-51305-1},
  url           = {https://doi.org/10.1016/S1874-5849(03)80041-4},
  mrclass       = {46L52 (46L51)},
  mrnumber      = {1999201},
  mrreviewer    = {Narcisse\ Randrianantoanina},
}

@article{PRD2022,
  author        = {Parcet, J. and Ricard, {\'E}. and de la Salle, M.},
  title         = {Fourier multipliers in {${\rm SL}_n({\bf R})$}},
  journal       = {Duke Math. J.},
  fjournal      = {Duke Mathematical Journal},
  volume        = {171},
  number        = {6},
  pages         = {1235--1297},
  year          = {2022},
  doi           = {10.1215/00127094-2021-0042},
  issn          = {0012-7094,1547-7398},
  url           = {https://doi.org/10.1215/00127094-2021-0042},
  mrclass       = {46L51 (22E46 43A80)},
  mrnumber      = {4408121},
  mrreviewer    = {Krishnan\ Parthasarathy},
}

@book{S1993,
  author        = {Stein, E.},
  title         = {Harmonic analysis: real-variable methods, orthogonality, and oscillatory
                  integrals},
  series        = {Princeton Mathematical Series},
  volume        = {43},
  pages         = {xiv+695},
  publisher     = {Princeton University Press, Princeton, NJ},
  year          = {1993},
  note          = {With the assistance of Timothy S. Murphy, Monographs in Harmonic Analysis, III},
  isbn          = {0-691-03216-5},
  mrclass       = {42-02 (35Sxx 43-02 47G30)},
  mrnumber      = {1232192},
  mrreviewer    = {Michael\ Cowling},
}

@article{T1975,
  author        = {Tomas, P.},
  title         = {A restriction theorem for the {F}ourier transform},
  journal       = {Bull. Amer. Math. Soc.},
  fjournal      = {Bulletin of the American Mathematical Society},
  volume        = {81},
  pages         = {477--478},
  year          = {1975},
  doi           = {10.1090/S0002-9904-1975-13790-6},
  issn          = {0002-9904},
  url           = {https://doi-org-s-224.libdb.csu.edu.cn/10.1090/S0002-9904-1975-13790-6},
  mrclass       = {42A92},
  mrnumber      = {358216},
  mrreviewer    = {P.\ G.\ Rooney},
}

@phdthesis{Ve1988,
  author        = {Vega, L.},
  title         = {El multiplicador de Schrodinger. La funci{\'o}n maximal y los operadores de
                  restricci{\'o}n},
  school        = {Universidad Aut{\'o}noma de Madrid},
  year          = {1988},
}

@phdthesis{Co1974,
  author        = {C{\'o}rdoba, A.},
  title         = {The Kakeya maximal function and the spherical summation multipliers.},
  school        = {The University of Chicago},
  year          = {1974},
}

@incollection{Biane2008RandomWalks,
  author        = {Biane, P.},
  editor        = {Franz, U. and Sch{\"u}rmann, M.},
  title         = {Introduction to random walks on noncommutative spaces},
  booktitle     = {Quantum Potential Theory},
  series        = {Lecture Notes in Mathematics},
  volume        = {1954},
  pages         = {61--116},
  publisher     = {Springer},
  address       = {Berlin, Heidelberg},
  year          = {2008},
  doi           = {10.1007/978-3-540-69365-9_3},
}

@unpublished{GLPPWZ2026,
  author        = {Gonz{\'a}lez-P{\'e}rez, A. and Lyu, B. and Parcet, J. and P{\'e}rez Garc{\'i}a,
                  J. and Wang, S. and Zhao, T.},
  title         = {Mixed-norm spaces for the free group and {SU(2)}-transference},
  year          = {2026},
  note          = {In preparation},
}

@article{W1995,
  author        = {Wolff, T.},
  title         = {An improved bound for {K}akeya type maximal functions},
  journal       = {Rev. Mat. Iberoamericana},
  fjournal      = {Revista Matem\'{a}tica Iberoamericana},
  volume        = {11},
  number        = {3},
  pages         = {651--674},
  year          = {1995},
  doi           = {10.4171/RMI/188},
  issn          = {0213-2230},
  url           = {https://doi-org-s-224.libdb.csu.edu.cn/10.4171/RMI/188},
  mrclass       = {42B25},
  mrnumber      = {1363209},
  mrreviewer    = {Detlef\ H.\ M\"{u}ller},
}

@incollection{Po2007,
  author        = {Popa, S.},
  title         = {Deformation and rigidity for group actions and von {N}eumann algebras},
  booktitle     = {International {C}ongress of {M}athematicians. {V}ol. {I}},
  pages         = {445--477},
  publisher     = {Eur. Math. Soc., Z\"{u}rich},
  year          = {2007},
  doi           = {10.4171/022-1/18},
  isbn          = {978-3-03719-022-7},
  url           = {https://doi-org-s-224.libdb.csu.edu.cn/10.4171/022-1/18},
  mrclass       = {46L35 (22D25 28D15 37A20 37A55 46L55)},
  mrnumber      = {2334200},
  mrreviewer    = {Bachir\ Bekka},
}

@article{HI2024,
  author        = {Houdayer, C. and Ioana, A.},
  title         = {Asymptotic freeness in tracial ultraproducts},
  journal       = {Forum Math. Sigma},
  fjournal      = {Forum of Mathematics. Sigma},
  volume        = {12},
  pages         = {Paper No. e88, 22},
  year          = {2024},
  doi           = {10.1017/fms.2024.93},
  issn          = {2050-5094},
  url           = {https://doi-org-s-224.libdb.csu.edu.cn/10.1017/fms.2024.93},
  mrclass       = {46L10 (03C66 46L51 46L52 46L54)},
  mrnumber      = {4814543},
  mrreviewer    = {Jorge\ Castillejos},
}

@article{FP1982,
  author        = {Fig\`a-Talamanca, A. and Picardello, M.},
  title         = {Spherical functions and harmonic analysis on free groups},
  journal       = {J. Funct. Anal.},
  fjournal      = {Journal of Functional Analysis},
  volume        = {47},
  number        = {3},
  pages         = {281--304},
  year          = {1982},
  doi           = {10.1016/0022-1236(82)90108-2},
  issn          = {0022-1236},
  url           = {https://doi.org/10.1016/0022-1236(82)90108-2},
  mrclass       = {22E30 (20E05 43A90)},
  mrnumber      = {665019},
  mrreviewer    = {Mogens\ Flensted-Jensen},
}

@book{Ta2002,
  author        = {Takesaki, M.},
  title         = {Theory of operator algebras. {I}},
  series        = {Encyclopaedia of Mathematical Sciences},
  volume        = {124},
  pages         = {xx+415},
  publisher     = {Springer-Verlag, Berlin},
  year          = {2002},
  note          = {Reprint of the first (1979) edition, Operator Algebras and Non-commutative
                  Geometry, 5},
  isbn          = {3-540-42248-X},
  mrclass       = {46Lxx (46-01)},
  mrnumber      = {1873025},
}

@book{Di1977,
  author        = {Dixmier, J.},
  title         = {{$C\sp*$}-algebras},
  series        = {North-Holland Mathematical Library},
  volume        = {Vol. 15},
  pages         = {xiii+492},
  publisher     = {North-Holland Publishing Co., Amsterdam-New York-Oxford},
  year          = {1977},
  note          = {Translated from the French by Francis Jellett},
  isbn          = {0-7204-0762-1},
  mrclass       = {46L05},
  mrnumber      = {458185},
}

@book{GoldsteinLabuschagne2025,
  author        = {Goldstein, S. and Labuschagne, L.},
  title         = {Noncommutative Measures and {$L^p$} and {Orlicz} Spaces, with Applications to
                  Quantum Physics},
  series        = {Oxford Graduate Texts in Mathematics},
  volume        = {33},
  pages         = {xiii+650},
  publisher     = {Oxford University Press},
  address       = {Oxford},
  year          = {2025},
  doi           = {10.1093/oso/9780198950202.001.0001},
  isbn          = {978-0-19-895020-2},
  url           = {https://doi.org/10.1093/oso/9780198950202.001.0001},
}

@article{Choi1974,
  author        = {Choi, M.-D.},
  title         = {A {Schwarz} inequality for positive linear maps on {$C^*$}-algebras},
  journal       = {Illinois J. Math.},
  fjournal      = {Illinois Journal of Mathematics},
  volume        = {18},
  number        = {4},
  pages         = {565--574},
  year          = {1974},
  doi           = {10.1215/ijm/1256051007},
  url           = {https://doi.org/10.1215/ijm/1256051007},
}

@techreport{Terp1981,
  author        = {Terp, M.},
  title         = {{$L^p$} Spaces Associated with von {Neumann} Algebras},
  number        = {3a+3b},
  institution   = {Matematisk Institut, K{\o}benhavns Universitet},
  address       = {Copenhagen},
  year          = {1981},
  month         = jun,
  note          = {Lecture notes},
  url           = {https://www.fuw.edu.pl/~kostecki/scans/terp1981.pdf},
  type          = {Rapport},
}

@book{Hiai2021,
  author        = {Hiai, F.},
  title         = {Lectures on Selected Topics in von {Neumann} Algebras},
  series        = {EMS Series of Lectures in Mathematics},
  pages         = {viii+241},
  publisher     = {EMS Press},
  address       = {Berlin},
  year          = {2021},
  doi           = {10.4171/ELM/32},
  isbn          = {978-3-98547-004-4},
  url           = {https://doi.org/10.4171/ELM/32},
}

@article{DMOS2008,
  author        = {Duoandikoetxea, J. and Moyua, A. and Oruetxebarria, O. and Seijo, E.},
  title         = {Radial {$A_p$} weights with applications to the disc multiplier and the
                  {Bochner--Riesz} operators},
  journal       = {Indiana Univ. Math. J.},
  fjournal      = {Indiana University Mathematics Journal},
  volume        = {57},
  number        = {3},
  pages         = {1261--1282},
  year          = {2008},
  doi           = {10.1512/iumj.2008.57.3282},
}

@book{Dix1981,
  author        = {Dixmier, J.},
  title         = {von {Neumann} Algebras},
  series        = {North-Holland Mathematical Library},
  volume        = {27},
  pages         = {xxxviii+437},
  publisher     = {North-Holland Publishing Co.},
  address       = {Amsterdam--New York},
  year          = {1981},
  note          = {Translated from the second French edition by F. Jellett. With a preface by E. C.
                  Lance},
  isbn          = {0-444-86308-7},
  mrnumber      = {0641217},
}

@misc{De2026,
  author        = {de la Salle, M.},
  title         = {Kakeya conjecture and High-Rank Lattice von {N}eumann algebras},
  year          = {2026},
  note          = {Preprint, arXiv:2602.14623; HAL: hal-05509698},
  eprint        = {2602.14623},
  archivePrefix = {arXiv},
  url           = {https://arxiv.org/abs/2602.14623},
}

@inproceedings{De2023,
  author        = {de la Salle, M.},
  title         = {Analysis on simple {L}ie groups and lattices},
  booktitle     = {International Congress of Mathematicians},
  pages         = {3166--3188},
  organization  = {European Mathematical Society-EMS-Publishing House GmbH},
  year          = {2023},
}

@book{Ga2014,
  author        = {Grafakos, L.},
  title         = {Classical {F}ourier analysis},
  series        = {Graduate Texts in Mathematics},
  volume        = {249},
  edition       = {Third},
  pages         = {xviii+638},
  publisher     = {Springer, New York},
  year          = {2014},
  doi           = {10.1007/978-1-4939-1194-3},
  isbn          = {978-1-4939-1193-6; 978-1-4939-1194-3},
  url           = {https://doi.org/10.1007/978-1-4939-1194-3},
  mrclass       = {42-01 (42Bxx)},
  mrnumber      = {3243734},
  mrreviewer    = {Atanas\ G.\ Stefanov},
}

\vfill

\begingroup
\small
\setlength{\tabcolsep}{7pt}
\renewcommand{\arraystretch}{1.06}
\newcommand{\paperauthorcontact}[3]{%
  \begin{minipage}[t]{0.47\textwidth}\vspace{0pt}%
    \begin{tabular}{@{}|p{\dimexpr\linewidth-2\tabcolsep-\arrayrulewidth\relax}@{}}
      \textbf{#1}\\
      \href{mailto:#2}{\texttt{#2}}\\
      #3
    \end{tabular}%
  \end{minipage}%
}
\par\bigskip\noindent
\paperauthorcontact{Adri\'an Gonz\'alez-P\'erez}
  {adrian.gonzalez@uam.es}
  {Universidad Aut\'onoma de Madrid}
\hfill
\paperauthorcontact{Binbing Lyu}
  {binbinglyu@stu.hit.edu.cn}
  {IASM\newline Harbin Institute of Technology}
\par\medskip\noindent
\paperauthorcontact{Javier Parcet}
  {parcet@icmat.es}
  {ICMAT--CSIC}
\hfill
\paperauthorcontact{Jorge P\'erez Garc\'ia}
  {jorge.perez@icmat.es}
  {IASM\newline Harbin Institute of Technology}
\par\medskip\noindent
\paperauthorcontact{Simeng Wang}
  {simeng.wang@hit.edu.cn}
  {IASM\newline Harbin Institute of Technology}
\hfill
\paperauthorcontact{Tiantian Zhao}
  {zhaotiantian@hit.edu.cn}
  {IASM\newline Harbin Institute of Technology}
\par
\endgroup

\end{document}